\documentclass[10pt]{amsart}

\usepackage{amssymb}
\usepackage{amscd}

\theoremstyle{definition}
\newtheorem{thm}{Theorem}[section]
\newtheorem{lem}[thm]{Lemma}
\newtheorem{prp}[thm]{Proposition}
\newtheorem{dfn}[thm]{Definition}
\newtheorem{cor}[thm]{Corollary}

\newtheorem{rmk}[thm]{Remark}

\newcommand{\beq}{\begin{equation}}
\newcommand{\eeq}{\end{equation}}
\newcommand{\beqr}{\begin{eqnarray*}}
\newcommand{\eeqr}{\end{eqnarray*}}
\newcommand{\bal}{\begin{align*}}
\newcommand{\eal}{\end{align*}}
\newcommand{\bei}{\begin{itemize}}
\newcommand{\eei}{\end{itemize}}

\newcommand{\af}{\alpha}
\newcommand{\bt}{\beta}
\newcommand{\gm}{\gamma}
\newcommand{\dt}{\delta}
\newcommand{\ep}{\varepsilon}

\newcommand{\et}{\eta}
\newcommand{\ch}{\chi}
\newcommand{\io}{\iota}

\newcommand{\ld}{\lambda}
\newcommand{\sm}{\sigma}

\newcommand{\ph}{\varphi}

\newcommand{\rh}{\rho}

\newcommand{\ta}{\tau}

\newcommand{\Z}{{\mathbb{Z}}}
\newcommand{\R}{{\mathbb{R}}}
\newcommand{\C}{{\mathbb{C}}}
\newcommand{\N}{{\mathbb{N}}}

\newcommand{\id}{{\mathrm{id}}}

\newcommand{\spec}{{\mathrm{sp}}}

\newcommand{\card}{{\mathrm{card}}}
\newcommand{\Aut}{{\mathrm{Aut}}}
\newcommand{\Ad}{{\mathrm{Ad}}}

\newcommand{\Zq}[1]{{\Z_{#1}}}
\newcommand{\Zqt}{\Zq{2}}
\newcommand{\Zqh}{\Zq{3}}
\newcommand{\Zqf}{\Zq{4}}
\newcommand{\Zqs}{\Zq{6}}
\newcommand{\Zqn}{\Zq{n}}
\newcommand{\Cs}[3]{C^* (\Zq{#1}, #2, #3)}

\newcommand{\CZnAa}{\Cs{n}{A}{\af}}

\newcommand{\andeqn}{\,\,\,\,\,\, {\mbox{and}} \,\,\,\,\,\,}
\newcommand{\QED}{\rule{0.4em}{2ex}}
\newcommand{\ts}[1]{{\textstyle{#1}}}
\newcommand{\ds}[1]{{\displaystyle{#1}}}
\newcommand{\Ssum}[2]{{\ts{ {\ds{\sum}}_{#1}^{#2} }}}
\newcommand{\ssum}[1]{{\ts{ {\ds{\sum}}_{#1} }}}

\newcommand{\ca}{C*-algebra}

\newcommand{\suca}{simple unital C*-algebra}
\newcommand{\ssuca}{simple separable unital C*-algebra}
\newcommand{\idsuca}{infinite dimensional simple unital C*-algebra}
\newcommand{\idssuca}{infinite dimensional
    simple separable unital C*-algebra}
\newcommand{\idfssuca}{infinite dimensional finite
    simple separable unital C*-algebra}
\newcommand{\uca}{unital C*-algebra}

\newcommand{\ct}{continuous}
\newcommand{\pj}{projection}

\newcommand{\hm}{homomorphism}
\newcommand{\wolog}{without loss of generality}
\newcommand{\Wolog}{Without loss of generality}

\newcommand{\ifo}{if and only if}

\newcommand{\mops}{mutually orthogonal \pj s}

\newcommand{\cfn}{continuous function}
\newcommand{\hsa}{hereditary subalgebra}

\newcommand{\mvnt}{Murray-von Neumann equivalent}
\newcommand{\mvnc}{Murray-von Neumann equivalence}
\newcommand{\tRp}{tracial Rokhlin property}
\newcommand{\sRp}{strict Rokhlin property}

\newcommand{\PSP}{Property~(SP)}

\renewcommand{\S}{\subset}

\newcommand{\SM}{\setminus}
\newcommand{\I}{\infty}

\title[The tracial Rokhlin property]{The tracial Rokhlin property
     for actions of finite groups on C*-algebras}

\author{N.~Christopher Phillips}

\date{25~August 2006}

\address{Department of Mathematics, University of Oregon,
       Eugene OR 97403-1222, USA.}

\email[]{ncp@darkwing.uoregon.edu}

\subjclass[2000]{Primary 46L55; Secondary 46L40.}
\thanks{Research partially supported by
 NSF grants DMS~0070776 and DMS~0302401.}

\begin{document}

\setcounter{section}{-1}

\begin{abstract}
We define ``tracial'' analogs of
the Rokhlin property for actions of finite groups,
approximate representability of actions of finite abelian groups,
and of approximate innerness.
We prove the following four
analogs of related ``nontracial'' results.
\begin{itemize}
\item
The crossed product of an \idssuca\  with tracial rank zero
by an action of a finite group with the tracial Rokhlin property
again has tracial rank zero.
\item
An outer action of a finite abelian group on
an \idssuca\  has the tracial Rokhlin property
\ifo\  its dual is tracially approximately representable,
and is tracially approximately representable
\ifo\  its dual has the tracial Rokhlin property.
\item
If a strongly tracially approximately inner action of a
finite cyclic group on an \idssuca\  has the tracial Rokhlin property,
then it is tracially approximately representable.
\item
An automorphism of an \idssuca\  $A$ with tracial rank zero
is tracially approximately inner
\ifo\  it is the identity on $K_0 (A)$ mod infinitesimals.
\end{itemize}
\end{abstract}

\maketitle

\section{Introduction}\label{Sec:Intro}

\indent
Tracially AF C*-algebras,
now known as C*-algebras with tracial rank zero,
were introduced in~\cite{LnTAF}.
Roughly speaking, a \ca\  has tracial rank zero if
the local approximation characterization of AF~algebras
holds after cutting out a ``small'' approximately central projection.
The term ``tracial'' comes from the fact that,
in good cases, a projection $p$ is ``small'' if $\ta (p) < \ep$
for every tracial state $\ta$ on $A.$
The classification~\cite{Ln15} of simple separable nuclear \ca s
with tracial rank zero and satisfying the Universal Coefficient Theorem
can be regarded as a vast generalization of the classification
of AF~algebras.
This success suggests that one consider ``tracial''
versions of other \ca\  concepts.

In this paper,
motivated by applications to particular crossed products
(see~\cite{PhtRp2} and~\cite{ELP}),
we formulate and prove ``tracial'' versions of the following
theorems:
\begin{itemize}
\item
The crossed product of an AF~algebra
by an action of a finite group with the Rokhlin property
is again~AF (Theorem~\ref{SRokhAF}).
\item
An action of a finite abelian group on
a unital \ca\  has the Rokhlin property
\ifo\  its dual is approximately representable,
and is approximately representable
\ifo\  its dual has the Rokhlin property.
(Lemma~3.8 of~\cite{Iz}).
\item
If an approximately inner action of a
finite cyclic group on a unital \ca\  has the Rokhlin property,
then it is approximately representable
(Proposition~\ref{AppInnRP}).
\item
An automorphism of an AF~algebra $A$
is approximately inner
\ifo\  it is the identity on $K_0 (A)$
(part of Theorem~3.1 of~\cite{Bl1}).
\end{itemize}
Our results are:
\begin{itemize}
\item
The crossed product of an \idssuca\  with tracial rank zero
by an action of a finite group with the tracial Rokhlin property
again has tracial rank zero (Theorem~\ref{RokhTAF}).
\item
An outer action of a finite abelian group on
an \idssuca\  has the tracial Rokhlin property
\ifo\  its dual is tracially approximately representable,
and is tracially approximately representable
\ifo\  its dual has the tracial Rokhlin property
(Theorem~\ref{TRPDualToTAppRep}).
\item
If a strongly tracially approximately inner action of a
finite cyclic group on an \idssuca\  has the tracial Rokhlin property,
then it is tracially approximately representable
(Theorem~\ref{FinOrdTAInn2}).
\item
An automorphism of an \idssuca\  $A$ with tracial rank zero
is tracially approximately inner
\ifo\  it is the identity on $K_0 (A)$ mod infinitesimals
(Theorem~\ref{TAIOnTAF}).
\end{itemize}

The first three of these results were chosen because they
are used in the proof~\cite{PhtRp2} that
every simple higher dimensional noncommutative torus
is an AT~algebra.
(We have not found the ``nontracial'' versions of the first and
third results in the literature.
Therefore they are also proved in this paper.)
The last is related to our effort to find the
``right'' definition of a tracially approximately inner automorphism.
Annoyingly, the strongly tracially approximately inner automorphisms
(as in the third result)
probably don't form a group.
On the other hand, an automorphism of an \idssuca\  with tracial rank zero
which is tracially approximately inner and has finite order
must in fact be strongly tracially approximately inner.
(Combine Proposition~\ref{TAIAndInf} and Theorem~\ref{FOTAIOnTAF}.)

In retrospect, the following motivation is perhaps better.
In~\cite{Iz} and~\cite{Iz2},
Izumi has started an intensive study of finite group actions with
the Rokhlin property, which, to minimize confusion,
we call here the strict Rokhlin property.
The strict Rokhlin property imposes severe restrictions on the
relation between the K-theory of the original algebra,
the action of the group on this K-theory,
and the K-theory of the crossed product.
See especially Section~3 of~\cite{Iz2}.
Since actions with the strict Rokhlin property are so rare,
a less restrictive version of the Rokhlin property is needed.
We give some examples; in them, we write $\Zqn$ for $\Z / n \Z.$
\begin{itemize}
\item
The flip action of $\Zqt$ on $A \otimes A,$ for an \idssuca\  $A,$
often has the tracial Rokhlin property,
but probably almost never has the strict Rokhlin property.
See~\cite{OP3}.
\item
Let $A$ be a simple higher dimensional noncommutative torus,
with standard unitary generators $u_1, u_2, \ldots, u_d.$
Consider the automorphism which sends $u_k$ to
$\exp (2 \pi i / n) u_k,$ and fixes $u_j$ for $j \neq k.$
This automorphism generates an action of $\Zqn$
which has the tracial Rokhlin property,
but for $n > 1$ never has strict Rokhlin property.
The fact that this action has the tracial Rokhlin property
plays a key role in the classification~\cite{PhtRp2} of
simple higher dimensional noncommutative toruses.
\item
Again let $A$ be a simple higher dimensional noncommutative torus,
with standard unitary generators $u_1, u_2, \ldots, u_d.$
The flip automorphism $u_k \mapsto u_k^*$
generates an action of $\Zqt$
which has the tracial Rokhlin property,
but never has the strict Rokhlin property.
See~\cite{ELP},
where this fact is used to prove that the crossed product by the
flip action is always~AF.
\item
The standard actions of $\Zqh,$ $\Zqf,$ and $\Zqs$
on an irrational rotation algebra all have the tracial Rokhlin property,
but never have the strict Rokhlin property.
In~\cite{ELP}, this is used to prove that the crossed
products are always AF~algebras.
\end{itemize}
Of course, one can't expect classification results for actions of the
kind found in~\cite{Iz} and~\cite{Iz2}.

This paper is devoted to the general theory.
In~\cite{PhtRp1b},
we give several useful criteria for the \tRp,
and we give a number of examples of actions of $\Zqt$
on \ca s with tracial rank zero (mostly AF~algebras)
which do and do not have the tracial Rokhlin property,
and are or are not tracially approximately representable.
Further examples,
and results for \ca s with finite but nonzero tracial rank,
will appear in~\cite{OP3}.
The main applications,
already mentioned above, are in \cite{PhtRp2} and~\cite{ELP}.
We also point out that the tracial Rokhlin property has
a generalization to integer actions,
considered in~\cite{OP1} and with applications given in~\cite{OP2}.
Presumably there is a useful generalization to other
countable amenable groups.

This paper replaces Sections 1 through~4
and Section~11 of the unpublished long preprint~\cite{PhW}.
The material of Sections 5 through~7 there will appear in~\cite{PhtRp2},
a greatly improved version of Sections 8 through~10
will appear in~\cite{ELP},
and an improved and expanded version of the material
in Sections 12 and~13
will appear in~\cite{PhtRp1b}.
We give the theory for actions of finite groups,
or of finite abelian groups, as appropriate;
in~\cite{PhW}, only finite cyclic groups were considered.
The definition of the \tRp\  given here differs slightly from
that in~\cite{PhW};
see Remark~\ref{MoreSmallness} and the following discussion for
details.
There is no difference for actions on simple unital \ca s
with tracial rank zero.
The definition of tracial approximate innerness
(Definition~\ref{D:TAI}) has been greatly improved;
more automorphisms satisfy the condition than satisfied the
condition in~\cite{PhW},
and with the new definition the tracially approximately inner
automorphisms form a group.
(See Theorem~\ref{T:TAIIsGp}.)
A variant of the original definition
appears in Definition~\ref{STAInnDfn} here,
where the condition is called strong tracial approximate innerness.
The material in Section~\ref{Sec:TrAppRep} on tracial
approximate representability is mostly new,
although it was motivated by one of the key results in~\cite{PhW}
and its parallel with Lemma~3.8 of~\cite{Iz}.

This paper is organized as follows.
In Section~\ref{Sec:TRP} we introduce the \tRp\  and prove
some basic properties.
Some of the lemmas will be used repeatedly in connection
with the tracial versions of other properties.
Section~\ref{Sec:CrPrd} contains the proofs that
the crossed product of an AF~algebra
by an action of a finite group with the Rokhlin property
is again~AF,
and that the crossed product of an \idssuca\  with tracial rank zero
by an action of a finite group with the tracial Rokhlin property
again has tracial rank zero.
In Section~\ref{Sec:TrAppRep} we treat tracial
approximate representability for actions of finite abelian groups,
and prove the duality between this property and the \tRp,
corresponding to Lemma~3.8 of~\cite{Iz}.
In Section~\ref{Sec:STAI} we introduce
strongly tracially approximately inner automorphism,
and prove that if a strongly tracially approximately inner action
of a
finite cyclic group on an \idssuca\  has the tracial Rokhlin property,
then it is tracially approximately representable.
We also prove the ``nontracial'' analog of this result.
Sections~\ref{Sec:TAI} and~\ref{Sec:PrTAI}
treat tracially approximately inner automorphisms.
We show that they form a group.
We prove that they act trivially
on the tracial state space and on $K_0$ mod infinitesimals,
and give a converse when the algebra is simple with tracial rank zero.
We also prove that if the the algebra is simple with tracial rank zero,
then a tracially approximately inner automorphism of finite order
is necessarily strongly tracially approximately inner.

We use the following notation.
We write $p \precsim q$ to mean that the \pj\  $p$ is \mvnt\  to a
sub\pj\  of $q,$ and $p \sim q$ to mean that
$p$ is \mvnt\  to $q.$
Also, $[a, b]$ denotes the additive commutator $a b - b a.$
If $A$ is a \ca\  and $\af \colon G \to \Aut (A)$ is a group action,
we write $A^{\af}$ for the fixed point algebra.

In most arguments dealing with an arbitrary finite subset $F$
of a \ca,
we will normalize and assume that all elements of $F$
have norm at most~$1.$

We are grateful to Hanfeng Li for valuable comments,
and to Hiroyuki Osaka for a careful reading and suggesting
improvements to several of the proofs.
We are also grateful to Dawn Ashley for catching a number of misprints
and minor mistakes.

\section{The tracial Rokhlin property}\label{Sec:TRP}

\indent
In this section we introduce the \tRp.
We observe several elementary relations and consequences,
and we prove several useful equivalent formulations.
Some of the technical lemmas will be repeatedly used
in connection with other ``tracial'' properties.

We begin with Izumi's definition of the Rokhlin property.
To emphasize
the difference, we call it the strict Rokhlin property here.

\begin{dfn}\label{ERPDfn}
Let $A$ be a separable unital \ca,
and let $\af \colon G \to \Aut (A)$
be an action of a finite group $G$ on $A.$
We say that $\af$ has the
{\emph{strict Rokhlin property}} if for every finite set
$F \S A,$ and every $\ep > 0,$
there are \mops\  $e_g \in A$ for $g \in G$ such that:
\begin{enumerate}
\item\label{ERPDfn:1}
$\| \af_g (e_h) - e_{g h} \| < \ep$ for all $g, h \in G.$
\item\label{ERPDfn:2}
$\| e_g a - a e_g \| < \ep$ for all $g \in G$ and all $a \in F.$
\item\label{ERPDfn:3}
$\sum_{g \in G} e_g = 1.$
\end{enumerate}
\end{dfn}

Izumi's definition (Definition~3.1 of~\cite{Iz})
is actually in terms of central sequences.
Thus, it yields not \mops\  but elements $b_g \in A$
such that $\| b_g b_h - \dt_{g, h} b_g \| < \ep$ for $g, h \in G,$
such that $\| b_g^* - b_g \| < \ep$ for $g \in G,$
and such that $\left\| 1 - \sum_{g \in G} b_g \right\| < \ep.$
However, with $n = \card (G),$ using semiprojectivity of $\C^n$
(see Lemma~14.1.5, Theorem~14.2.1, Theorem~14.1.4,
and Definition~14.1.1 of~\cite{Lr})
and a suitably smaller choice of $\ep,$
one easily sees that the definition above is equivalent to
Definition~3.1 of~\cite{Iz}.

If $\af$ is approximately inner,
requiring $\sum_{g \in G} e_g = 1$ forces $[1_A] \in K_0 (A)$
to be divisible by the order of $G,$
and therefore rules out many \ca s of interest.
In fact, the strict Rokhlin property imposes much more stringent
conditions on the K-theory.
Theorem~3.3 and Lemma~3.2(1) of~\cite{Iz2} show that if
a nontrivial finite group $G$ acts on a \suca\  $A$
in such a way that the induced action on $K_* (A)$ is trivial,
and if one of $K_0 (A)$ and $K_1 (A)$ is a nonzero free abelian
group, then $\af$ does not have the strict Rokhlin property.
Theorem~3.3 and the discussion preceding Theorem~3.4 of~\cite{Iz2}
show that if in addition $G$ is cyclic of order $n,$ then
the strict Rokhlin property implies that
$K_* (A)$ is uniquely $n$-divisible.
It follows that the actions considered in our main applications
(to simple higher dimensional noncommutative toruses~\cite{PhtRp2} and
irrational rotation algebras~\cite{ELP})
never have the strict Rokhlin property.

We now give the definition of the \tRp.
The difference is that we do not require that $\sum_{g \in G} e_g = 1,$
only that $1 - \sum_{g \in G} e_g$ be ``small'' in a tracial sense.
Of course, $\sum_{g \in G} e_g = 1$ is allowed,
in which case Conditions~(\ref{NewTRPDfn:3}) and~(\ref{NewTRPDfn:4})
in the definition are vacuous.

\begin{dfn}\label{NewTRPDfn}
Let $A$ be an \idssuca,
and let $\af \colon G \to \Aut (A)$
be an action of a finite group $G$ on $A.$
We say that $\af$ has the
{\emph{tracial Rokhlin property}} if for every finite set
$F \S A,$ every $\ep > 0,$
and every positive element $x \in A$ with $\| x \| = 1,$
there are \mops\  $e_g \in A$ for $g \in G$ such that:
\begin{enumerate}
\item\label{NewTRPDfn:1}
$\| \af_g (e_h) - e_{g h} \| < \ep$ for all $g, h \in G.$
\item\label{NewTRPDfn:2}
$\| e_g a - a e_g \| < \ep$ for all $g \in G$ and all $a \in F.$
\item\label{NewTRPDfn:3}
With $e = \sum_{g \in G} e_g,$ the \pj\  $1 - e$ is \mvnt\  to a
\pj\  in the \hsa\  of $A$ generated by $x.$
\item\label{NewTRPDfn:4}
With $e$ as in~(\ref{NewTRPDfn:3}), we have $\| e x e \| > 1 - \ep.$
\end{enumerate}
\end{dfn}

\begin{rmk}\label{MoreSmallness}
Our original definition, in~\cite{PhW},
in addition specified a positive integer $N,$
and required the following condition instead of~(\ref{NewTRPDfn:4}):
\begin{itemize}
\item[(4$'$)]
For every $g \in G,$
there are $N$ \mops\  $f_1, f_2, \dots, f_N \leq e_j,$
each of which is \mvnt\  to the \pj\  $1 - e$ of~(\ref{NewTRPDfn:3}).
\end{itemize}
\end{rmk}

As we will see in Lemma~\ref{TRPForFinite} below,
when $A$ is finite,
Condition~(\ref{NewTRPDfn:4})
in Definition~\ref{NewTRPDfn} is unnecessary.
When $A$ has finite tracial topological rank
in the sense of~\cite{LnTTR},
it is not hard to see that
Definition~\ref{NewTRPDfn} implies the definition in Remark~\ref{MoreSmallness}.
In general, the situation is less clear,
and it might be necessary to use both Condition~(\ref{NewTRPDfn:4})
and Condition~(4$'$),
especially for nonsimple \ca s.
We postpone further discussion to Section~4 of~\cite{PhtRp1b}.

\begin{rmk}\label{SRPImpTRP}
If an action of a finite group $G$
on an \idssuca\  $A$ has the
strict Rokhlin property, then it has the \tRp.
\end{rmk}

\begin{lem}\label{TRPImpOuter}
Let $A$ be an \idssuca,
and let $\af \colon G \to \Aut (A)$
be an action of a finite group $G$ on $A$ which has the \tRp.
Then $\af_g$ is outer for every $g \in G \SM \{ 1 \}.$
\end{lem}

\begin{proof}
Let $g \in G \SM \{ 1 \}$ and let $u \in A$ be unitary.
We prove that $\af_g \neq \Ad (u).$
Apply Definition~\ref{NewTRPDfn} with $F = \{ u \},$
with $\ep = \frac{1}{2},$ and with $x = 1.$
Then $e_1$ and $e_g$ are orthogonal nonzero \pj s, so
\[
\| \af_g (e_1) - u e_1 u^* \|
  \geq \| e_g - e_1 \|
     - \| \af_g (e_1) - e_{g} \| - \| u e_1 u^* - e_1 \|
  > 0.
\]
Therefore $\af_g \neq \Ad (u).$
\end{proof}

\begin{cor}\label{CrPrIsSimple}
Let $A$ be an \idssuca,
and let $\af \colon G \to \Aut (A)$
be an action of a finite group $G$ on $A$ which has the \tRp.
Then $C^* (G, A, \af)$ is simple.
\end{cor}

\begin{proof}
In view of Lemma~\ref{TRPImpOuter},
this follows from Theorem~3.1 of~\cite{Ks1}.
\end{proof}

For the \tRp\  to be likely to hold, the \ca\  must have a
reasonable number of \pj s.
For reference, we recall here the definition of the property that
seems most relevant.

\begin{dfn}\label{SPD}
Let $A$ be a \ca.
We say that $A$ has {\emph{Property~(SP)}} if every nonzero
\hsa\  in $A$ contains a nonzero \pj.
\end{dfn}

We state here some results about simple \ca s with Property~(SP)
that will be used repeatedly in this paper.

\begin{lem}\label{L:Switch}
Let $A$ be a \ca, and let $c \in A.$
Then for any projection $p \in {\overline{c A c^*}},$
there exists a \pj\  $q \in {\overline{c^* A c}}$
such that $p \sim q.$
\end{lem}

\begin{proof}
This is essentially in Section~1 of~\cite{Cn1}.
The details can be found in the proof of Lemma~4.1 of~\cite{OP1}.
\end{proof}

The following lemma is essentially Lemma~3.1 of~\cite{LnTAF},
but no proof is given there.

\begin{lem}\label{L:CompSP}
Let $A$ be a simple \ca\  with Property~(SP).
Let $B \S A$ be a nonzero \hsa,
and let $p \in A$ be a nonzero \pj.
Then there is a nonzero \pj\  $q \in B$
such that $q \precsim p.$
\end{lem}

\begin{proof}
Choose a nonzero positive element $a \in B.$
Since $A$ is simple,
there exists $x \in A$ such that $c = a x p$ is nonzero.
Choose a nonzero \pj\  $q \in {\overline{c A c^*}}.$
Then $q \in B,$ and Lemma~4.1 of~\cite{OP1}
(or Lemma~\ref{L:Switch}) provides a \pj\  $e \leq p$
such that $q \sim e.$
\end{proof}

\begin{lem}\label{OrthInSP}
Let $A$ be an \idsuca\  with Property~(SP).
Let $B \S A$ be a nonzero \hsa,
and let $n \in \N.$
Then there exist nonzero \mvnt\  \mops\  $p_1, p_2, \ldots, p_n \in B.$
\end{lem}

\begin{proof}
Since $A$ is unital and infinite dimensional,
it is not isomorphic to the algebra of compact operators
on any Hilbert space.
The lemma is then immediate from Lemma~3.2 of~\cite{LnTAF}.
(The result from~\cite{AS} is on page~61 of that reference.
See above for the proof of Lemma~3.1 of~\cite{LnTAF}.)
\end{proof}

\begin{lem}\label{L:MnSP}
Let $A$ be an \idsuca, and let $n \in \N.$
Then $A$ has Property~(SP) \ifo\  $M_n \otimes A$ has Property~(SP).
Moreover, in this case,
for every nonzero \hsa\  $B \S M_n \otimes A,$
there exists a nonzero \pj\  $p \in A$ such that
$1 \otimes p$ is \mvnt\  to a \pj\  in $B.$
\end{lem}

\begin{proof}
Since $A$ is isomorphic to a \hsa\  in $M_n \otimes A,$
it is obvious that if $M_n \otimes A$ has Property~(SP),
then so does $A.$

For the converse and the last statement,
choose a nonzero element $x \in B.$
Let $( e_{j, k} )_{1 \leq j, k \leq n}$ be a system of matrix units
for $M_n.$
Choose $j$ such that $(e_{j, j} \otimes 1) x \neq 0.$
Then
$C = (e_{j, j} \otimes 1) x (M_n \otimes A) x^* (e_{j, j} \otimes 1)$
is a nonzero \hsa\  in
$(e_{j, j} \otimes 1) (M_n \otimes A) (e_{j, j} \otimes 1) \cong A.$
Because $A$ has Property~(SP),
there exists a nonzero \pj\  $f \in A$
such that $e_{j, j} \otimes f \in C.$
By Lemma~\ref{OrthInSP},
there exist nonzero \mvnt\  \mops\  $f_1, f_2, \ldots, f_n \in A$
with $f_j \leq f$ for all $k.$
Then
\[
1 \otimes f_1
 = \sum_{k = 1}^n e_{k, k} \otimes f_1
 \sim \sum_{k = 1}^n e_{j, j} \otimes f_k
 \leq e_{j, j} \otimes f,
\]
and $e_{j, j} \otimes f$ is \mvnt\  to a \pj\  in $B$ by
Lemma~\ref{L:Switch}.
\end{proof}

The following result implies, in particular,
that if $A$ is an \idsuca\  with Property~(SP),
and if $\af \colon G \to \Aut (A)$
is an action of a finite group $G$ on $A$
such that the crossed product is also simple,
then the crossed product has Property~(SP).

\begin{prp}\label{SPForCrPrd}
Let $A$ be an \idsuca\  with Property~(SP),
and let $\af \colon G \to \Aut (A)$
be an action of a finite group $G$ on $A$
such that $C^* (G, A, \af)$ is also simple.
Let $B \S C^* (G, A, \af)$ be a nonzero \hsa.
Then there exists a nonzero \pj\  $p \in A$
which is \mvnt\  in $C^* (G, A, \af)$ to a \pj\  in $B.$
\end{prp}

\begin{proof}
Set $N = \{ g \in G \colon {\mbox{$\af_g$ is inner}} \}.$
Theorem~4.2 of~\cite{JO} provides a \pj\  %
$q \in C^* (N, A, \af |_N ) \S C^* (G, A, \af)$
which is \mvnt\  in $C^* (G, A, \af)$ to a \pj\  in $B.$
(There is a reference missing in its proof:
one uses the discussion after Proposition~2.8 of~\cite{RfFG}
for the claim, at the top of page~295 of~\cite{JO},
that $C^* (N, A, \af |_N )$ is a direct sum of matrix
algebras over $A.$)
Let $G$ act on $C^* (N, A, \af |_N )$
by conjugation by the standard unitaries $u_g \in C^* (G, A, \af),$
as in the discussion
before Proposition~2.4 of~\cite{RfFG}.
(Since $A$ is simple and unital, an automorphism
is partly inner in the sense of~\cite{RfFG} \ifo\  it is inner.)
Following the discussion after Proposition~2.8 of~\cite{RfFG},
there is a finite dimensional \ca\  $C$ with an action $\gm$ of $G$
such that $C^* (N, A, \af |_N )$ is equivariantly isomorphic
to $A \otimes C$ with the action $g \mapsto \af_g \otimes \gm_g,$
in such a way that the inclusion of $A$ in $C^* (N, A, \af |_N )$
becomes $a \mapsto a \otimes 1.$
Moreover, $C$ is $G$-simple by Propositions~2.5 and~2.10 of~\cite{RfFG}.
Let $e_1, e_2, \ldots, e_n$ be the minimal central \pj s of $C.$
There is some $k$ such that the \pj\  $(1 \otimes e_k) q$ is nonzero.
Lemma~\ref{L:MnSP} provides a nonzero \pj\  $f \in A$
such that $f \otimes e_k \precsim (1 \otimes e_k) q.$
By Lemma~\ref{OrthInSP},
there exist nonzero \mvnt\  \mops\  $f_1, f_2, \ldots, f_n \in A$
with $f_j \leq f$ for all $j.$
Since $C$ is $G$-simple, there exist $g_1, g_2, \ldots, g_n \in G$
such that $\gm_{g_j} (e_j) = e_k$ for $1 \leq j \leq n.$
Use Lemma~\ref{L:CompSP} repeatedly to find a \pj\  $p \in A$
such that $p \precsim \af_{g_j}^{-1} (f_j)$ for all $j.$
Then, with the first relation holding in $C^* (G, A, \af)$
and the second in $C^* (N, A, \af |_N ) = A \otimes C,$
we have
\[
p \otimes e_j
  \sim \af_{g_j} (p) \otimes e_k
  \precsim f_j \otimes e_k.
\]
It follows that, in  $C^* (G, A, \af),$
we have
\[
p \otimes 1
  \precsim \sum_{j = 1}^n f_j \otimes e_k
  \leq f \otimes e_k
  \precsim (1 \otimes e_k) q
  \leq q.
\]
Since $q$ is \mvnt\  in $C^* (G, A, \af)$ to a \pj\  in $B,$
the proof is complete.
\end{proof}

We now return to the main development.

\begin{lem}\label{RImpSP}
Let $A$ be an \idssuca,
and let $\af \colon G \to \Aut (A)$
be an action of a finite group $G$ on $A$ which has the \tRp.
Then $A$ has Property~(SP) or $\af$ has the strict Rokhlin property.
\end{lem}

\begin{proof}
If $A$ does not have Property~(SP),
then there is a nonzero positive element $x \in A$
which generates a \hsa\  which contains no nonzero \pj.
\end{proof}

The following two lemmas will be important for dealing
with Condition~(\ref{NewTRPDfn:4}) of Definition~\ref{NewTRPDfn},
and with similar conditions in other definitions in this paper.
The first comes from an argument
that goes back to Cuntz, in the proof of Lemma~1.7 of~\cite{Cn2}.

\begin{lem}\label{PjAndNorm}
Let $A$ be a \ca\  with Property~(SP),
let $x \in A$ be a positive element with $\| x \| = 1,$
and let $\ep > 0.$
Then there exists a nonzero \pj\  $p \in {\overline{x A x}}$
such that,
whenever $q \leq p$ is a nonzero \pj,
then
\[
\| q x q - q \| < \ep, \,\,\,\,\,\,
\| q x - x q \| < \ep, \andeqn
\| q x q \| > 1 - \ep.
\]
\end{lem}

\begin{proof}
Choose \cfn s $h_1, h_2 \colon [0, 1] \to [0, 1]$ such that
$h_1 (0) = 0,$ $h_1 (t) = 1$ for $t \geq 1 - \tfrac{1}{4} \ep,$
and $| h_1 (t) - t | \leq \tfrac{1}{4} \ep$ for all $t,$
and such that $h_2 (1) = 1$ and $h_1 h_2 = h_2.$
Set $y = h_1 (x)$ and $z = h_2 (x).$
Then $\| x - y \| \leq \tfrac{1}{4} \ep$ and $y z = z.$ 
Furthermore, $z \neq 0$ because $1 \in \spec (x).$
By Property~(SP),
there is a nonzero \pj\  $p \in {\overline{z A z}}.$
Let $q \leq p$ be a nonzero \pj.
Then $q \in {\overline{z A z}},$
so $y q = q y = q.$
Therefore
\[
\| q x - x q \|
 \leq 2 \| x - y \|
 \leq \tfrac{1}{2} \ep
 < \ep.
\]
Furthermore,
\[
\| q x q - q \|
 = \| q x q - q y q \|
 \leq \| x - y \|
 \leq \tfrac{1}{4} \ep
 < \ep,
\]
whence also $\| q x q \| > 1 - \ep.$
\end{proof}

\begin{lem}\label{FPjNorm}
Let $A$ be an infinite dimensional finite unital \ca\  %
with Property~(SP),
let $x \in A$ be a positive element with $\| x \| = 1,$
and let $\ep > 0.$
Then there exists a nonzero \pj\  $q \in {\overline{x A x}}$
such that,
whenever $e \in A$ is a \pj\  such that $1 - e \precsim q,$
then $\| e x e \| > 1 - \ep.$
\end{lem}

\begin{proof}
Apply Lemma~\ref{PjAndNorm} with $x^{1/2}$ in place of $x$
and with $\tfrac{1}{5} \ep$ in place of $\ep,$
obtaining a nonzero \pj\  $p \in {\overline{x^{1/2} A x^{1/2}}}$
such that, whenever $q \leq p$ is a nonzero \pj,
then (combining two of the estimates there)
$\| q x^{1/2} - q \| < \tfrac{2}{5} \ep.$
Note that ${\overline{x^{1/2} A x^{1/2}}} = {\overline{x A x}}.$
By Lemma~\ref{OrthInSP},
there is a nonzero \pj\  $q \leq p$ such that $p - q \neq 0.$
Now suppose $e \in A$ is a \pj\  such that $1 - e \precsim q$
and $\| e x e \| \leq 1 - \ep.$
Then, repeatedly using $\| a^* a \| = \| a a^* \|,$ we get
\[
\| e p e \|
   = \| p e p \|
   < \big\| p x^{1/2} e x^{1/2} p \big\| + \tfrac{4}{5} \ep
   \leq \big\| x^{1/2} e x^{1/2} \big\| + \tfrac{4}{5} \ep
   = \| e x e \| + \tfrac{4}{5} \ep
   < 1 - \tfrac{1}{5} \ep.
\]
Therefore
\[
\| e - e (1 - p) \|
 = \| e p \|
 = \| e p e \|^{1/2}
 < \left( 1 - \tfrac{1}{5} \ep \right)^{1/2}
 < 1.
\]
It follows from Lemma~2.5.2 of~\cite{LnBook}
that $e \precsim 1 - p.$
Since, by assumption, we have $1 - e \precsim q,$
this gives $1 \precsim 1 - (p - q),$
contradicting finiteness of $A.$
\end{proof}

When $A$ is finite,
we do not need Condition~(\ref{NewTRPDfn:4})
of Definition~\ref{NewTRPDfn}:

\begin{lem}\label{TRPForFinite}
Let $A$ be an \idfssuca,
and let $\af \colon G \to \Aut (A)$
be an action of a finite group $G$ on $A.$
Then $\af$ has the
tracial Rokhlin property \ifo\  for every finite set
$F \subset A,$ every $\ep > 0,$
and every nonzero positive element $x \in A,$
there are \mops\  $e_g \in A$ for $g \in G$ such that:
\begin{enumerate}
\item\label{TRPForFinite:1} %
$\| \af_g (e_h) - e_{g h} \| < \ep$ for all $g, h \in G.$
\item\label{TRPForFinite:2} %
$\| e_g a - a e_g \| < \ep$ for all $g \in G$ and all $a \in F.$
\item\label{TRPForFinite:3} %
With $e = \sum_{g \in G} e_g,$ the \pj\  $1 - e$ is \mvnt\  to a
\pj\  in the \hsa\  of $A$ generated by $x.$
\end{enumerate}
\end{lem}

\begin{proof}
It is immediate that the \tRp\  implies the condition in the lemma.
So assume the condition in the lemma holds.

If $A$ does not have Property~(SP),
then $\af$ has the strict Rokhlin property,
by the same proof as for Lemma~\ref{RImpSP}.
Accordingly, we may assume that $A$ has Property~(SP).

Let $F \S A$ be finite,
let $\ep > 0,$
and let $x \in A$ be a positive element with $\| x \| = 1.$
Apply Lemma~\ref{FPjNorm},
obtaining a nonzero \pj\  $q \in {\overline{x A x}}$
such that,
whenever $e \in A$ is a \pj\  such that $1 - e \precsim q,$
then $\| e x e \| > 1 - \ep.$
Apply the hypothesis,
with $F$ and $\ep$ as given,
and with $q$ in place of $x,$
obtaining \pj s $e_g \in A$ for $g \in G.$
Set $e = \sum_{g \in G} e_g.$
We need only prove that $\| e x e \| > 1 - \ep.$
But this is immediate from the choice of $q$
and the relation $1 - e \precsim q.$
\end{proof}

It is convenient to have a formally stronger version of the \tRp,
in which the defect \pj\  is $\af$-invariant.

\begin{lem}\label{StTRPDfn}
Let $A$ be an \idssuca,
and let $\af \colon G \to \Aut (A)$
be an action of a finite group $G$ on $A$ which has the \tRp.
Let $F \S A$ be finite,
let $\ep > 0,$
and let $x \in A$ be a positive element with $\| x \| = 1.$
Then there are \mops\  $e_g \in A$ for $g \in G$ such that:
\begin{enumerate}
\item\label{StTRPDfn:1} %
$\| \af_g (e_h) - e_{g h} \| < \ep$ for all $g, h \in G.$
\item\label{StTRPDfn:2} %
$\| e_g a - a e_g \| < \ep$ for all $g \in G$ and all $a \in F.$
\item\label{StTRPDfn:3} %
The \pj\  $e = \sum_{g \in G} e_g$ is $\af$-invariant.
\item\label{StTRPDfn:4} %
With $e$ as in~(\ref{StTRPDfn:3}), the \pj\  $1 - e$ is \mvnt\  to a
\pj\  in the \hsa\  of $A$ generated by $x.$
\item\label{StTRPDfn:5} %
With $e$ as in~(\ref{StTRPDfn:3}), we have $\| e x e \| > 1 - \ep.$
\end{enumerate}
\end{lem}

\begin{proof}
\Wolog\  $\| a \| \leq 1$ for all $a \in F.$
Set $\ep_0 = \min \left( \frac{1}{5} \ep, \frac{1}{2} \right).$
Choose $\dt > 0$ so small that whenever
$B$ is a \uca, $b \in B$ is selfadjoint, and $p \in B$ is a \pj\  %
such that $\| b - p \| < \card (G) \dt,$
then the functional calculus
$e = \ch_{(1/2, \, \I)} (b)$ is defined and moreover
$\| e - p \|$ is small enough that there exists a unitary
$v \in B$ such that $v p v^* = e$ and $\| v - 1 \| < \ep_0.$ 
We also require $\dt \leq \ep_0.$

Apply Definition~\ref{NewTRPDfn} to $\af,$
with $F$ and $x$ as given,
and with $\dt$ in place of $\ep.$
obtaining \pj s $e_g \in A$ for $g \in G.$
Let $(p_g)_{g \in G}$ be the resulting family of \pj s.
Define $p = \sum_{h \in G} p_h.$
For $g \in G$ we have
\[
\| \af_g (p) - p \|
 \leq \sum_{h \in G} \| \af_g (p_h) - p_{g h} \|
 < \card (G) \dt.
\]
Set
\[
b = \frac{1}{\card (G)} \sum_{g \in G} \af_g (p).
\]
Then $b \in A^{\af}$ and $\| b - p \| < \card (G) \dt.$

By the choice of $\dt,$
there exists a \pj\  $e \in A^{\af}$
and a unitary $v \in A$ such that
$v p v^* = e$ and $\| v - 1 \| < \ep_0.$
Now define $e_g = v p_g v^*$ for $g \in G.$
Clearly $\| e_g - p_g \| < 2 \ep_0.$
So, for $g, h \in G,$
\[
\| \af_g (e_h) - e_{g h} \|
 \leq \| e_h - p_h \| + \| e_{g h} - p_{g h} \|
     + \| \af_g (p_h) - p_{g h} \|
 < 5 \ep_0
 \leq \ep.
\]
For $g \in G$ and $a \in F,$ we similarly
get $\| e_g a - a e_g \| < 5 \ep_0 \leq \ep.$
We have
\[
\| (1 - e) - (1 - p) \| < 2 \ep_0 \leq 1,
\]
so $1 - e \sim 1 - p,$ and is hence \mvnt\  to a
\pj\  in the \hsa\  of $A$ generated by $x.$
Finally,
\[
\| e x e \| \geq \| p x p \| - 2 \| e - p \|
   > 1 - \dt - 2 \ep_0
   \geq 1 - \ep.
\]
This completes the proof.
\end{proof}

\section{Crossed products by actions on C*-algebras with
   tracial rank zero}\label{Sec:CrPrd}

\indent
The main result of this section is that the crossed product of an
\idssuca\  with tracial rank zero by an action with the \tRp\  %
again has tracial rank zero.
We begin by proving the analogous result
in which the word ``tracial'' is omitted everywhere:
the crossed product of a
unital AF~-algebra by an action with the \sRp\  is again~AF.
To our surprise,
we have been unable to find this result in the literature.
Both proofs have been considerably simplified from original versions,
following a suggestion of Osaka.
The main estimates in the proof for actions with the \sRp\  will
be referred to in the proof for actions with the \tRp.
Example~3.1 of~\cite{PhtRp1b} shows that
the hypotheses can't be weakened to require only the \tRp.

We need a lemma.

\begin{lem}\label{SemiProjMn}
Let $n \in \N.$
For every $\ep > 0$ there is $\dt > 0$ such that,
whenever
$(e_{j, k})_{1 \leq j, k \leq n}$ is a
system of matrix units for $M_n,$
whenever $B$ is a unital \ca,
and whenever $w_{j, k},$ for $1 \leq j, k \leq n,$
are elements of $B$ such that
$\| w_{j, k}^* - w_{k, j} \| < \dt$ for $1 \leq j, k \leq n,$
such that
$\| w_{j_1, k_1} w_{j_2, k_2} - \dt_{j_2, k_1} w_{j_1, k_2} \| < \dt$
for $1 \leq j_1, j_2, k_1, k_2 \leq n,$
and such that the $w_{j, j}$ are orthogonal \pj s
with $\sum_{j = 1}^n w_{j, j} = 1,$
then there exists a unital \hm\  $\ph \colon M_n \to B$
such that $\ph (e_{j, j}) = w_{j, j}$ for $1 \leq j \leq n$
and $\| \ph (e_{j, k}) - w_{j, k} \| < \ep$ for $1 \leq j, k \leq n.$
\end{lem}

\begin{proof}
This follows from semiprojectivity of $M_n$
(see Lemma~14.1.5, Theorem~14.2.2, Theorem~14.1.4,
and Definition~14.1.1 of~\cite{Lr}),
and the fact that if two families of $n$ orthogonal \pj s summing to $1$
are close
then there is a unitary close to $1$ which conjugates
one family to the other
(Lemma~2.5.7 of~\cite{LnBook}).
\end{proof}

\begin{thm}\label{SRokhAF}
Let $A$ be a unital AF~algebra.
Let $\af \colon G \to \Aut (A)$
be an action of a finite group $G$ on $A$ which has the \sRp.
Then $C^* (G, A, \af)$ is an AF~algebra.
\end{thm}

\begin{proof}
We prove that for every finite set $S \S C^* (G, A, \af)$
and every $\ep > 0,$
there is an AF~subalgebra $D \S C^* (G, A, \af)$
such that every element of $S$ is within $\ep$ of an element of $D.$
It is then easy to use Theorem~2.2 of~\cite{Brt} to show that
$C^* (G, A, \af)$ is~AF.
It suffices to consider a finite set of the form
$S = F \cup \{ u_g \colon g \in G \},$
where $F$ is a finite subset of the unit ball of $A$ and
$u_g \in C^* (G, A, \af)$
is the canonical unitary implementing the automorphism $\af_g.$
So let $F \S A$ be a finite subset
with $\| a \| \leq 1$ for all $a \in F$ and let $\ep > 0.$

Set $n = \card (G),$
and set $\ep_0 = \ep / (4 n).$
Choose $\dt > 0$ according to Lemma~\ref{SemiProjMn}
for $n$ as given and for $\ep_0$ in place of $\ep.$
Also require $\dt \leq \ep / [2 n (n + 1)].$
Apply the \sRp\  to $\af$ with $F$ as given
and with $\dt$ in place of $\ep,$
obtaining \pj s $e_g \in A$ for $g \in G.$
Define $w_{g, h} = u_{g h^{-1}} e_h$ for $g, h \in G.$

We claim that the $w_{g, h}$ form a $\dt$-approximate system of
$n \times n$ matrix units in $C^* (G, A, \af).$
We estimate:
\[
\| w_{g, h}^* - w_{h, g} \|
  = \| e_h u_{g h^{-1}}^* - u_{h g^{-1}} e_g \|
  = \| u_{g h^{-1}} e_h u_{g h^{-1}}^* - e_g \|
  = \| \af_{g h^{-1}} (e_h)  - e_g \|
  < \dt.
\]
Also, using $e_g e_h = \dt_{g, h} e_h$ at the second step,
\begin{align*}
\| w_{g_1, h_1} w_{g_2, h_2} - \dt_{g_2, h_1} w_{g_1, h_2} \|
& = \big\| u_{g_1 h_1^{-1}} e_{h_1} u_{g_2 h_2^{-1}} e_{h_2}
        - \dt_{g_2, h_1} u_{g_1 h_2^{-1}} e_{h_2} \big\|        \\
& = \big\| u_{g_1 h_1^{-1}} e_{h_1} u_{g_2 h_2^{-1}} e_{h_2}
     - u_{g_1 h_1^{-1} g_2 h_2^{-1}} e_{h_2 g_2^{-1} h_1} e_{h_2} \big\|
                        \\
& = \big\| u_{g_1 h_1^{-1} g_2 h_2^{-1}}
       \big( u_{g_2 h_2^{-1}}^* e_{h_1} u_{g_2 h_2^{-1}}
                                - e_{h_2 g_2^{-1} h_1} \big)
           e_{h_2} \big\|
  < \dt.
\end{align*}
Finally, $\sum_{g \in G} w_{g, g} = \sum_{g \in G} e_{g} = 1.$
This proves the claim.

Let $(v_{g, h})_{g, h \in G}$ be a system of matrix units for $M_n.$
By the choice of $\dt,$
there exists a unital \hm\  $\ph_0 \colon M_n \to C^* (G, A, \af)$
such that
$\| \ph_0 (v_{g, h}) - w_{g, h} \| < \ep_0$ for all $g, h \in G,$
and $\ph_0 (v_{g, g}) = e_g$ for all $g \in G.$
Now define a unital \hm\  %
$\ph \colon M_n \otimes e_1 A e_1 \to C^* (G, A, \af)$
by $\ph (v_{g, h} \otimes a) = \ph_0 (v_{g, 1}) a \ph_0 (v_{1, h})$
for $g, h \in G$ and $a \in e_1 A e_1.$
It is well known that a corner of an AF~algebra is~AF,
and $\ph$ is injective, so $D = \ph (M_n \otimes e_1 A e_1)$
is an AF~subalgebra of $C^* (G, A, \af).$
We complete the proof by
showing that every element of $S$ is within $\ep$ of an element of $D.$

For $g \in G$ we have $\sum_{h \in G} \ph_0 (v_{g h, h}) \in D$
and
\begin{align*}
\left\| u_g - \ssum{h \in G} \ph_0 (v_{g h, h}) \right\|
& \leq \sum_{h \in G} \| u_g e_h - \ph_0 (v_{g h, h}) \|      \\
& = \sum_{h \in G} \| w_{g h, h} - \ph_0 (v_{g h, h}) \|
  < n \ep_0
  \leq \ep.
\end{align*}
Now let $a \in F.$
Set
\[
b = \sum_{g \in G} v_{g, g} \otimes e_1 \af_g^{-1} (a) e_1
       \in M_n \otimes e_1 A e_1.
\]
Using $\| e_g a e_h \| \leq \| [e_g, a] \| + \| a e_g e_h \|,$
we get
\[
\left\| a - \ssum{g \in G} e_g a e_g \right\|
 \leq \sum_{g \neq h} \| e_g a e_h \|
 < n (n - 1) \dt. 
\]
We use this, and the inequalities
\[
\| \ph_0 (v_{g, 1}) e_1 - u_g e_1 \| < \ep_0
\andeqn
\big\| e_1 \af_g^{-1} (a) e_1 - \af_g^{-1} (e_g a e_g) \big\| < 2 \dt,
\]
to get
\begin{align*}
\| a - \ph (b) \|
& = \left\| a - \ssum{g \in G}
       \ph_0 (v_{g, 1}) e_1 \af_g^{-1} (a) e_1 \ph_0 (v_{1, g}) \right\|
                   \\
& < 2 n \ep_0
      + \left\| a - \ssum{g \in G}
               u_g e_1 \af_g^{-1} (a) e_1 u_g^* \right\|
                   \\
& < 2 n \ep_0 + 2 n \dt
      + \left\| a - \ssum{g \in G}
               u_g \af_g^{-1} (e_g a e_g) u_g^* \right\|
                   \\
& < 2 n \ep_0 + 2 n \dt + n (n - 1) \dt
  \leq \ep.
\end{align*}
This completes the proof.
\end{proof}

The following result gives the criterion we use for a
\ssuca\  to have tracial rank zero.
Note that, by Theorem~7.1(a) of~\cite{LnTTR},
tracial rank zero is the same as
tracially AF in the sense of Definition~2.1 of~\cite{LnTAF}.

\begin{prp}\label{TAFCond}
Let $A$ be a \ssuca.
Then $A$ has tracial rank zero in the sense of
Definition~3.1 of~\cite{LnTTR} \ifo\  the following holds.

For every finite set $F \S A,$ every $\ep > 0,$ and
every nonzero positive element $x \in A,$
there is a \pj\  $p \in A$ and a finite
dimensional unital subalgebra $E \S p A p$ (that is,
$p$ is the identity of $E$) such that:
\begin{enumerate}
\item\label{TAFCond:1} %
$\| p a - a p \| < \ep$ for all $a \in F.$
\item\label{TAFCond:2} %
For every $a \in F$ there exists $b \in E$ such that
$\| p a p - b \| < \ep.$
\item\label{TAFCond:3} %
$1 - p$ is \mvnt\  to a \pj\  in ${\overline{x A x}}.$
\end{enumerate}
\end{prp}

\begin{proof}
Theorem~6.13 and Definition~3.4 of~\cite{LnTTR}
give this result, except with unitary equivalence
instead of \mvnc\  in~(\ref{TAFCond:3}).
However, according to Remark~6.12 and Theorem~6.9 of~\cite{LnTTR},
with~(\ref{TAFCond:3}) as it stands,
the resulting condition implies that $A$ has stable rank one.
\end{proof}

The condition of~Proposition~\ref{TAFCond} is given as the definition
in~\cite{LnBook}.
See Definition 3.6.2 there.

For convenient reference,
we recall some properties of \suca s with tracial rank zero.
(Most of them will not be needed until later.)

\begin{dfn}\label{OrdDetD}
Let $A$ be a unital \ca.
We say that the {\emph{order on \pj s over $A$ is determined by traces}}
if whenever $n \in \N$ and $p, q \in M_n (A)$ are \pj s such that
$\ta (p) < \ta (q)$ for all tracial states $\ta$ on $A,$ then
$p \precsim q.$
\end{dfn}

This is just Blackadar's Second Fundamental Comparability Question
for all matrix algebras over $A.$
See 1.3.1 in~\cite{Bl3}.

\begin{thm}\label{TAFProp}
(H.~Lin.)
Let $A$ be a \ssuca\  with tracial rank zero.
Then $A$ has real rank zero and stable rank one.
Moreover, the order on \pj s over $A$ is determined by traces
(Definition~\ref{OrdDetD}).
\end{thm}

\begin{proof}
In view of Theorem~7.1(a) of~\cite{LnTTR},
real rank zero and stable rank one are Theorem~3.4 of \cite{LnTAF}.
That the order is determined by traces is
Corollary~5.7 and Theorems~5.8 and~6.8 of~\cite{LnTTR}.
\end{proof}

\begin{thm}\label{RokhTAF}
Let $A$ be an \idssuca\  with tracial rank zero.
Let $\af \colon G \to \Aut (A)$
be an action of a finite group $G$ on $A$ which has the \tRp.
Then $C^* (G, A, \af)$ has tracial rank zero.
\end{thm}

\begin{proof}
The proof is a modification of that of Theorem~\ref{SRokhAF}.
It suffices to verify the condition of Proposition~\ref{TAFCond},
for a finite set $S$ of the form $S = F \cup \{ u_g \colon g \in G \},$
where $F$ is a finite subset of the unit ball of $A$ and
$u_g \in C^* (G, A, \af)$
is the canonical unitary implementing the automorphism $\af_g.$
So let $F \S A$ be a finite subset
with $\| a \| \leq 1$ for all $a \in F,$ let $\ep > 0,$
and let $x \in C^* (G, A, \af)$ be a nonzero positive element.

The \ca\  $A$ has Property~(SP) by Theorem~\ref{TAFProp}.
So Proposition~\ref{SPForCrPrd} and Corollary~\ref{CrPrIsSimple} provide
a nonzero \pj\  $q \in A$ which is \mvnt\  in
$C^* (G, A, \af)$ to a \pj\  in
${\overline{x C^* (G, A, \af) x}}.$
By Lemma~\ref{OrthInSP},
there are orthogonal nonzero \pj s $q_1, q_2 \in A$
such that $q_1, q_2 \leq q.$

Set $n = \card (G),$
and set $\ep_0 = \ep / (16 n).$
Choose $\dt > 0$ according to Lemma~\ref{SemiProjMn}
for $n$ as given and for $\ep_0$ in place of $\ep.$
Also require $\dt \leq \ep / [8 n (n + 1)].$
Apply Lemma~\ref{StTRPDfn} to $\af,$
with $F$ as given,
with $\dt$ in place of $\ep,$
and with $q_1$ in place of $x,$
obtaining \pj s $e_g \in A$ for $g \in G.$
Set $e = \sum_{g \in G} e_g.$
By construction, $u_g e u_g^* = \af_g (e) = e$
for every $g \in G.$
Also, for $a \in F$ we have
$\| e a - a e \| \leq \sum_{g \in G} \| e_g a - a e_g \| < n \ep_0.$

Define $w_{g, h} = u_{g h^{-1}} e_h$ for $g, h \in G.$
Using the same estimates as in the proof of Theorem~\ref{SRokhAF},
we find an injective unital \hm\  %
$\ph \colon M_n \otimes e_1 A e_1 \to e C^* (G, A, \af) e$
and a finite set $T$ in the closed unit ball
of $M_n \otimes e_1 A e_1$
such that for every $a \in S = F \cup \{ u_g \colon g \in G \},$
there is $b \in T$
such that $\| \ph (b) - e a e \| < \frac{1}{4} \ep.$
Let $e_{1, 1} \in M_n$ denote the usual $(1, 1)$ matrix unit.
Then $\ph$ furthermore has the property that if $a \in e_1 A e_1$
then $\ph (e_{1, 1} \otimes a) = a.$
Use Lemma~\ref{L:CompSP} to choose equivalent nonzero \pj s
$f_1, f_2 \in A$ such that
$f_1 \leq e_1$ and $f_2 \leq q_2.$
It follows from Theorems~3.10 and~3.12(1) of~\cite{LnTAF}
that $M_n \otimes e_1 A e_1$ has tracial
rank zero,
so there is a \pj\  $p_0 \in M_n \otimes e_1 A e_1$ and a finite
dimensional unital subalgebra $E_0 \S p_0 (M_n \otimes e_1 A e_1) p_0$
such that $\| p_0 b - b p_0 \| < \tfrac{1}{4} \ep$ for all $b \in T,$
such that for every $b \in T$ there exists $c \in E_0$ with
$\| p_0 b p_0 - c \| < \tfrac{1}{4} \ep,$
and such that $1 - p_0 \precsim e_{1, 1} \otimes f_1$
in $M_n \otimes e_1 A e_1.$
Set $p = \ph (p_0),$ and set $E = \ph (E_0),$
which is a finite dimensional unital subalgebra
of $p C^* (G, A, \af) p.$

Let $a \in S.$
Choose $b \in T$ such that $\| \ph (b) - e a e \| < \frac{1}{4} \ep.$
Then, using $p e = e p = p,$
\begin{align*}
\| p a - a p \|
& \leq 2 \| e a - a e \| + \| p e a e - e a e p \|      \\
& \leq 2 \| e a - a e \| + 2 \| e a e - \ph (b) \|
              + \| p_0 b - b p_0 \|
 < 2 n \ep_0 + 2 \left( \tfrac{1}{4} \ep \right) + \tfrac{1}{4} \ep
 \leq \ep.
\end{align*}
Further,
choosing $c \in E_0$ such $\| p_0 b p_0 - c \| < \tfrac{1}{4} \ep,$
the element $\ph (c)$ is in $E$ and satisfies
\[
\| p a p - \ph (c) \|
  \leq \| e a e - \ph (b) \| + \| p_0 b p_0 - c \|
  < \tfrac{1}{4} \ep + \tfrac{1}{4} \ep
  \leq \ep.
\]
Finally, in $C^* (G, A, \af)$ we have
\[
1 - p = (1 - e) + (e - p)
  \precsim q_1 + f_2
  \leq q,
\]
and $q$ is \mvnt\  to a \pj\  in ${\overline{x C^* (G, A, \af) x}}.$
\end{proof}

\section{Tracially approximately representable actions
   and duality}\label{Sec:TrAppRep}

In this section, we give the tracial analog of
approximate representability of an action,
Definition~3.6(2) of~\cite{Iz}.
As there, we restrict to abelian groups;
see Remark~3.7 of~\cite{Iz}.
After several elementary properties and reformulations,
we generalize Lemma~3.8 of~\cite{Iz},
showing that an action is tracially approximately representable
\ifo\  the dual action has the \tRp,
and similarly with the action and its dual exchanged.

Definition~3.6(2) of~\cite{Iz} is formulated in terms of
central sequences,
so we give a reformulation with direct estimates.

\begin{lem}\label{AppRep}
Let $A$ be a separable unital \ca,
and let $\af \colon G \to \Aut (A)$
be an action of a finite abelian group $G$ on $A.$
Then $\af$ is approximately representable
(Definition~3.6(2) of~\cite{Iz})
\ifo\  for every finite set $F \S A$ and every $\ep > 0,$
there are unitaries $w_g \in A$ such that:
\begin{enumerate}
\item\label{AppRep:1} %
$\| \af_g (a) - w_g a w_g^* \| < \ep$
for all $a \in F$ and all $g \in G.$
\item\label{AppRep:2} %
$\| w_g w_h - w_{g h} \| < \ep$ for all $g, h \in G.$
\item\label{AppRep:3} %
$\| \af_g (w_h) - w_h \| < \ep$ for all $g, h \in G.$
\end{enumerate}
\end{lem}

\begin{proof}
It is immediate that the condition of the lemma implies
approximate representability.

Now assume $\af$ is approximately representable,
and let $F \S A$ be finite and let $\ep > 0.$
\Wolog\  $\| a \| \leq 1$ for all $a \in F.$
Set $n = \card (G).$
Using semiprojectivity of $C^* (G) \cong \C^n$
(see Lemma~14.1.5, Theorem~14.2.1, Theorem~14.1.4,
and Definition~14.1.1 of~\cite{Lr}),
choose $\dt > 0$ such that whenever $A$ is a unital \ca\  and
elements $x_g \in A,$ for $g \in G,$ satisfy
$\| x_g x_h - x_{g h} \| < \dt$ for all $g, h \in G$
and
$\| x_g^* x_g - 1 \|, \, \| x_g x_g^* - 1 \| < \dt$
for all $g \in G,$
then there exist unitaries $w_g \in A$ such that
$\| w_g w_h - w_{g h} \| < \ep$ for all $g, h \in G$
and $\| w_g - x_g \| <  \frac{1}{3} \ep$ for all $g \in G.$
{}From Definition~3.6(2) of~\cite{Iz},
we get elements $x_g \in A$ such that
the conditions above are satisfied, and also
$\| \af_g (a) - x_g a x_g^* \| < \frac{1}{3} \ep$
for all $a \in F$ and all $g \in G$
and
$\| \af_g (x_h) - x_h \| < \frac{1}{3} \ep$ for all $g, h \in G.$
Moreover, we may clearly require $\| x_g \| \leq 1$ for all $g \in G.$
It is now easy to check that
Conditions~(\ref{AppRep:1}) through~(\ref{AppRep:3}) hold.
\end{proof}

We now give the tracial analog.

\begin{dfn}\label{TrAppRepDfn}
Let $A$ be an \idssuca,
and let $\af \colon G \to \Aut (A)$
be an action of a finite abelian group $G$ on $A.$
We say that $\af$ is
{\emph{tracially approximately representable}}
if for every finite set $F \S A,$ every $\ep > 0,$
and every positive element $x \in A$ with $\| x \| = 1,$
there are a \pj\  $e \in A$
and unitaries $w_g \in e A e$ such that:
\begin{enumerate}
\item\label{TrAppRepDfn:1} %
$\| e a - a e \| < \ep$ for all $a \in F.$
\item\label{TrAppRepDfn:2} %
$\| \af_g (e a e) - w_g e a e w_g^* \| < \ep$
for all $a \in F$ and all $g \in G.$
\item\label{TrAppRepDfn:3} %
$\| w_g w_h - w_{g h} \| < \ep$ for all $g, h \in G.$
\item\label{TrAppRepDfn:4} %
$\| \af_g (w_h) - w_h \| < \ep$ for all $g, h \in G.$
\item\label{TrAppRepDfn:5} %
$1 - e$ is \mvnt\  to a
\pj\  in the \hsa\  of $A$ generated by $x.$
\item\label{TrAppRepDfn:6} %
$\| e x e \| > 1 - \ep.$
\end{enumerate}
\end{dfn}

We next give several elementary properties and reformulations.
In the rest of this section,
we generalize Lemma~3.8 of~\cite{Iz},
showing that an action is tracially approximately representable
\ifo\  the dual action has the \tRp,
and similarly with the action and its dual exchanged.

\begin{lem}\label{AppRepandTrAppRep}
Let $A$ be an \idssuca,
and let $\af \colon G \to \Aut (A)$
be an action of a finite abelian group $G$ on $A.$
If $\af$ is approximately representable in the sense of
Definition~3.6(2)
of~\cite{Iz},
than $\af$ is tracially approximately representable.
If $\af$ is tracially approximately representable,
then $A$ has Property~(SP) or $\af$ is approximately representable.
\end{lem}

\begin{proof}
This is immediate from Lemma~\ref{AppRep}.
\end{proof}

When $A$ is finite,
we do not need Condition~(\ref{TrAppRepDfn:6})
of Definition~\ref{TrAppRepDfn}:

\begin{lem}\label{TARForFinite}
Let $A$ be an \idfssuca,
and let $\af \colon G \to \Aut (A)$
be an action of a finite abelian group $G$ on $A.$
Then $\af$ is tracially approximately representable \ifo\ %
for every finite set $F \S A,$ every $\ep > 0,$
and every positive element $x \in A$ with $\| x \| = 1,$
there are a \pj\  $e \in A$
and unitaries $w_g \in e A e$ such that:
\begin{enumerate}
\item\label{TARForFinite:1} %
$\| e a - a e \| < \ep$ for all $a \in F.$
\item\label{TARForFinite:2} %
$\| \af_g (e a e) - w_g e a e w_g^* \| < \ep$
for all $a \in F$ and all $g \in G.$
\item\label{TARForFinite:3} %
$\| w_g w_h - w_{g h} \| < \ep$ for all $g, h \in G.$
\item\label{TARForFinite:4} %
$\| \af_g (w_h) - w_h \| < \ep$ for all $g, h \in G.$
\item\label{TARForFinite:5} %
$1 - e$ is \mvnt\  to a
\pj\  in the \hsa\  of $A$ generated by $x.$
\end{enumerate}
\end{lem}

\begin{proof}
The proof is the same as for Lemma~\ref{TRPForFinite}.
\end{proof}

In the definition of tracial approximate representability,
we can require invariance instead of approximate invariance,
and we can require that $g \mapsto w_g$ be a \hm.

\begin{lem}\label{STrAppRep}
Let $A$ be an \idssuca,
and let $\af \colon G \to \Aut (A)$
be a tracially approximately representable
action of a finite abelian group $G$ on $A.$
Then for every finite set
$F \S A,$ every $\ep > 0,$
and every positive element $x \in A$ with $\| x \| = 1,$
there are an $\af$-invariant \pj\  $e \in A$
and a \hm\  $g \mapsto w_g$ from $G$ to the unitary group
of $e A e$ such that:
\begin{enumerate}
\item\label{STrAppRep:1} %
$\| e a - a e \| < \ep$ for all $a \in F.$
\item\label{STrAppRep:2} %
$\| \af_g (e a e) - w_g e a e w_g^* \| < \ep$
for all $a \in F$ and all $g \in G.$
\item\label{STrAppRep:3} %
$w_g w_h = w_{g h}$ for all $g, h \in G.$
\item\label{STrAppRep:4} %
$\af_g (w_h) = w_h$ for all $g, h \in G.$
\item\label{STrAppRep:5} %
$1 - e$ is \mvnt\  to a
\pj\  in the \hsa\  of $A$ generated by $x.$
\item\label{STrAppRep:6} %
$\| e x e \| > 1 - \ep.$
\end{enumerate}
\end{lem}

\begin{proof}
Let $F,$ $\ep,$ and $x$ be given.
\Wolog\  $\| a \| \leq 1$ for all $a \in F.$

The group algebra $C^* (G)$ is finite dimensional,
hence semiprojective.
(See Lemma~14.1.5 and Theorems~14.2.1 and~14.2.2 of~\cite{Lr}.)
In particular, by Theorem~14.1.4 of~\cite{Lr},
its standard generators $u_g$ for $g \in G,$
and relations
($u_g$ is unitary and $u_g u_h = u_{g h}$ for $g, h \in G$)
are stable in the sense of Definition~14.1.1 of~\cite{Lr}.
Accordingly, there exists $\ep_0 > 0$ such that whenever
$(x_g)_{g \in G}$ is a collection of elements of a unital \ca\  $B$
such that for all $g, h \in G$ we have
\[
\| x_g x_g^* - 1 \| < 10 \ep_0,
\,\,\,\,\,\,
\| x_g^* x_g - 1 \| < 10 \ep_0,
\andeqn
\| x_g x_h - x_{g h} \| < 10 \ep_0,
\]
then there exists a unital \hm\  $\ph \colon C^* (G) \to B$
such that $\| \ph (u_g) - x_g \| < \frac{1}{13} \ep$ for all $g \in G.$
We also require
$\ep_0 \leq \min \left( \frac{1}{13} \ep, \, 1 \right).$

Next, choose $\dt > 0$ with $\dt \leq \ep_0$
and so small that
whenever $p$ is a \pj\  in a \ca\  $B,$
and $a \in B$ is selfadjoint and satisfies $\| a - p \| < 3 \dt,$
then the \pj\  $q = \ch_{[1/2, \, \I)} (a)$ is defined and
satisfies $\| q - p \| < \ep_0.$

Apply Definition~\ref{TrAppRepDfn}
with $F$ as given, with $\dt$ in place of $\ep,$
and with $x$ as given,
obtaining a \pj\  $f$
and unitaries $y_g \in f A f.$
Note that $\| y_1^2 - y_1 \| < \dt,$ so $\| y_1 - f \| < \dt.$
Therefore $\| \af_g (f) - f \| < 3 \dt$ for all $g \in G.$
Consequently, the element
\[
b = \frac{1}{\card (G)} \sum_{h \in G} \af_h (f) \in A^{\af}
\]
satisfies $\| b - f \| < 3 \dt,$
so that there is a \pj\  $e \in A^{\af}$
such that $\| e - f \| < \ep_0.$

For every $g, h \in G,$ we now have
\[
\| \af_h (e y_g e) - y_g \|
 \leq \| e y_g e - y_g \| + \| \af_h (y_g) - y_g \|
 < 2 \| e - f \| + \dt
 < 3 \ep_0.
\]
Therefore the elements
\[
x_g
 = \frac{1}{\card (G)} \sum_{h \in G} \af_h (e y_g e) \in (e A e)^{\af}
\]
satisfy $\| x_g - y_g \| < 3 \ep_0.$
We now estimate:
\[
\| x_g x_g^* - f \|
 \leq 2 \| x_g - y_g \| + \| y_g y_g^* - f \|
 < 6 \ep_0 + 0
 = 6 \ep_0;
\]
similarly,
$\| x_g^* x_g - 1 \| < 6 \ep_0$;
and also
\[
\| x_g x_h - x_{g h} \|
  \leq \| x_g - y_g \| + \| x_h - y_h \| + \| x_{g h} - y_{g h} \|
    + \| y_g y_h - y_{g h} \|
  < 9 \ep_0 + \dt
  \leq 10 \ep_0.
\]
So there is a unital \hm\  $\ph \colon C^* (G) \to (e A e)^{\af}$
such that $\| \ph (u_g) - x_g \| < \frac{1}{13} \ep$ for all $g \in G.$
Set $w_g = \ph (u_g).$
Then
$\| w_g - y_g \| < \frac{1}{13} \ep + 3 \ep_0 \leq \frac{4}{13} \ep.$

We now have, for $a \in F,$
\[
\| e a - a e \|
  \leq 2 \| e - f \| + \| f a - a f \|
  < 2 \ep_0 + \dt
  < \ep,
\]
and, for $a \in F$ and $g \in G,$
\begin{align*}
\| \af_g (e a e) - w_g e a e w_g^* \|
 & \leq 4 \| e - f \| + 2 \| w_g - y_g \|
            + \| \af_g (f a f) - y_g f a f y_g^* \|
                     \\
 & < 4 \ep_0 + \tfrac{8}{13} \ep + \dt
   \leq \ep.
\end{align*}
This gives Parts (\ref{STrAppRep:1})--(\ref{STrAppRep:4})
of the conclusion.

For Part~(\ref{STrAppRep:5}),
\[
\| e x e \|
 \geq \| f x f \| - 2 \| e - f \|
 > 1 - \dt - 2 \ep_0
 > 1 - \ep.
\]
For Part~(\ref{STrAppRep:6}),
use $\| e - f \| < \ep_0 \leq 1$ to get $e \sim f.$
\end{proof}

We now turn to the proof of the duality relations.
We will use Takai duality to get one part from the other,
and to do so we need to know
that an action of a finite abelian group has the \tRp\  %
\ifo\  the second dual action has the \tRp.
The next few lemmas contain the proof of this fact.

We must identify what happens to the inclusion map
$A \to C^* \big( {\widehat{G}}, \, C^* (G, A, \af),
            \, {\widehat{\af}} \big)$
under Takai duality~\cite{Tk}.
The formula in the next proposition is correct without simplicity,
but assuming simplicity shortens the proof.

\begin{prp}\label{Takai}
Let $A$ be a \suca,
and let $\af \colon G \to \Aut (A)$
be an action of a finite abelian group $G$ on $A.$
In $L (l^2 (G)),$
let $(e_{g, h})_{g, h \in G}$ be the family of
matrix units determined by the requirement that
$e_{g, h}$ send the standard basis vector $\dt_h$
to the standard basis vector $\dt_g,$
and vanish on all other standard basis vectors.
Let
\[
\io \colon A \to C^* (G, A, \af)
\andeqn
\mu \colon C^* (G, A, \af)
 \to C^* \big( {\widehat{G}},
          \, C^* (G, A, \af), \, {\widehat{\af}} \big)
\]
be the inclusions.
Then there exists an isomorphism
$\ph \colon C^* \big( {\widehat{G}},
               \, C^* (G, A, \af), \, {\widehat{\af}} \big)
   \to L (l^2 (G)) \otimes A$
such that
\[
(\ph \circ \mu \circ \io) (a)
 = \sum_{k \in G} e_{k, k} \otimes \af_k^{-1} (a)
\]
for all $a \in A.$
\end{prp}

\begin{proof}
For $g \in G$ let $u_g \in C^* (G, A, \af)$
be the standard implementing unitary,
and for $\ta \in {\widehat{G}}$
let
$v_{\ta} \in C^* \big( {\widehat{G}},
                \, C^* (G, A, \af), \, {\widehat{\af}} \big)$
be the standard implementing unitary.
Identify $A$ and $C^* (G, A, \af)$
with their images in
$C^* \big( {\widehat{G}}, \, C^* (G, A, \af), \, {\widehat{\af}} \big)$
under $\mu \circ \io$ and $\mu.$
Then define $\ph$ on the generators of
$C^* \big( {\widehat{G}}, \, C^* (G, A, \af), \, {\widehat{\af}} \big)$
by
\[
\ph (a) = \sum_{k \in G} e_{k, k} \otimes \af_k^{-1} (a),
\,\,\,\,\,\,
\ph (u_g) = \sum_{k \in G} e_{g k, k} \otimes 1,
\andeqn
\ph (v_{\ta}) = \sum_{k \in G} e_{k, k} \otimes \ta (k)
\]
for $a \in A,$
$g \in G,$ and
$\ta \in {\widehat{G}}.$
One checks that the appropriate relations are satisfied
for this definition to extend to a \hm,
and it is easily seen that the resulting \hm\  is surjective.
It is injective because, by Takai duality~\cite{Tk},
simplicity of $A$ implies
simplicity of
$C^* \big( {\widehat{G}}, \, C^* (G, A, \af), \, {\widehat{\af}} \big).$
\end{proof}

\begin{lem}\label{TRPCorner}
Let $A$ be an \idssuca,
and let $\af \colon G \to \Aut (A)$
be an action of a finite group $G$ on $A$ which has the \tRp.
Let $p \in A$ be a $G$-invariant \pj.
Then the action $g \mapsto \af_g |_{p A p}$ has the \tRp.
\end{lem}

\begin{proof}
Let $F \S p A p$ be finite,
let $\ep > 0,$
and let $x \in p A p$ be a positive element with $\| x \| = 1.$
Set $n = \card (G).$
Set
\[
\ep_0 = \min \left( \frac{1}{n}, \frac{\ep}{4 n + 1} \right).
\]
Using semiprojectivity of $\C^n,$
choose $\dt > 0$
such that whenever $B$ is a unital \ca,
$q_1, \ldots, q_n \in B$ are \mops,
and $p \in B$ is a \pj\  such that $\| p q_j - q_j p \| < \dt$
for $1 \leq j \leq n,$
then there are \mops\  $e_j \in p B p$ such that
$\| e_j - p q_j p \| < \ep_0$ for $1 \leq j \leq n.$
We also require $\dt \leq \ep_0.$

Apply Definition~\ref{NewTRPDfn} to $\af,$
with $F \cup \{ p \}$ in place of $F,$
with $\dt$ in place of $\ep,$
and with $x$ as given,
obtaining \pj s $q_g \in A$ for $g \in G.$
By the choice of $\dt,$
there are \mops\  $e_g \in p A p$ such that
$\| e_g - p q_g p \| < \ep_0$ for $g \in G.$
We now estimate, using $\af_g (p) = p,$
\[
\| \af_g (e_h) - e_{g h} \|
 \leq \| e_h - p q_h p \| + \| e_{g h} - p q_{g h} p \|
           + \| p ( \af_g (q_h) - q_{g h} ) p \|
 < 2 \ep_0 + \dt
 \leq \ep,
\]
and for $a \in F,$ using $p a = a p = a,$
\[
\| e_g a - a e_g \|
  \leq 2 \| e_g - p q_g p \| + \| p (q_g a - a q_g) p \|
 < 2 \ep_0 + \dt
 \leq \ep.
\]
Next, set $e = \sum_{g \in G} e_g$ and $q = \sum_{g \in G} q_g.$
Then $\| e - p q p \| < n \ep_0 \leq 1.$
So
\begin{align*}
\| (1 - q) (p - e) - (p - e) \|
& = \| q ( p - e ) \|
  = \| ( p - e ) q ( p - e ) \|^{1/2}           \\
& = \| ( 1 - e ) p q p ( 1 - e ) \|^{1/2}
  \leq \| e - p q p \|^{1/2}
  < 1.
\end{align*}
It follows from Lemma~2.5.2 of~\cite{LnBook}
that $p - e \precsim 1 - q.$
Since $1 - q$ is \mvnt\  to a
\pj\  in ${\overline{x p A p x}} = {\overline{x A x}},$
so is $p - e.$

Finally, we estimate $\| e x e \|.$
{}From $\| q_g p - p q_g \| < \dt \leq \ep_0$ for all $g,$
we get $\| q p - p q \| < n \ep_0,$
whence $\| e - p q \|, \| e - q p \| < 2 n \ep_0.$
Therefore
\[
\| e x e \| > \| q p x p q \| - 4 n \ep_0
 = \| q x q \| - 4 n \ep_0
 > 1 - \dt - 4 n \ep_0
 \geq 1 - \ep.
\]
This completes the proof.
\end{proof}

\begin{lem}\label{PosMat}
Let $A$ be a \ca, let $n \in \N,$
and let $x = ( x_{j, k})_{1 \leq j, k \leq n} \in M_n (A)$ be positive.
Then there exists $k$
such that $\| x_{k, k} \| \geq n^{-2} \| x \|.$
\end{lem}

\begin{proof}
We first claim that if
\[
y = \left( \begin{array}{cc}
  a     &  b        \\
  b^*   &  c
\end{array} \right)
\in M_2 (A)
\]
is positive, then $\| b \| \leq \tfrac{1}{2} ( \| a \| + \| c \| ).$
\Wolog\  $A \S L (H)$ for some Hilbert space $H,$
and we correspondingly take $M_2 (A) \S L (H^2).$
Let $\xi, \et \in H$ satisfy $\| \xi \| = \| \et \| = 1.$
Choose $t \in \C$ with $| t | = 1$
such that
$t \langle b \et, \, \xi \rangle = - | \langle b \et, \, \xi \rangle |.$
Then one calculates that
\[
0 \leq \langle y (\xi, t \et), \, (\xi, t \et) \rangle
  = \langle a \xi, \xi \rangle + \langle c \et, \et \rangle
      - 2 | \langle b \et, \, \xi \rangle |
  \leq \| a \| + \| c \| - 2 | \langle b \xi, \, \et \rangle |.
\]
Therefore
$| \langle b \et, \, \xi \rangle |
        \leq \tfrac{1}{2} ( \| a \| + \| c \| ).$
The claim follows by taking the supremum over all $\xi, \et \in H$
with $\| \xi \| = \| \et \| = 1.$

By considering $2 \times 2$ submatrices of $x,$
the claim implies that
$\| x_{j, k} \| \leq \tfrac{1}{2} ( \| x_{j, j} \| + \| x_{k, k} \| )$
for $1 \leq j, k \leq n.$
Sum over all $j$ and $k$ to get
\[
\| x \| \leq \sum_{1 \leq j, k \leq n} \| x_{j, k} \|
 \leq n \sum_{k = 1}^n \| x_{k, k} \|.
\]
The statement of the lemma follows.
\end{proof}

\begin{lem}\label{TRPMat}
Let $A$ be an \idssuca,
let $\af \colon G \to \Aut (A)$
be an action of a finite group $G$ on $A,$
and let $g \mapsto v_g$ be a unitary representation of $G$ on $\C^n.$
Then $\af$ has the \tRp\  \ifo\  the action
$g \mapsto \Ad (v_g) \otimes \af_g$ of $G$ on $M_n \otimes A$
has the \tRp.
\end{lem}

\begin{proof}
We first prove that if $\af$ has the \tRp,
then so does $g \mapsto \Ad (v_g) \otimes \af_g.$
If $M_n \otimes A$ does not have Property~(SP),
then, by Lemma~\ref{L:MnSP}, neither does $A,$
so $\af$ has the strict Rokhlin property by Lemma~\ref{RImpSP}.
Then we must show that $g \mapsto \Ad (v_g) \otimes \af_g$
has the strict Rokhlin property.
The proof is similar to but easier than the other case, and is omitted.

So assume that $M_n \otimes A$ has Property~(SP)
and that $\af$ has the \tRp.
Let $F \S M_n \otimes A$ be finite,
let $\ep > 0,$
and let $x \in M_n \otimes A$ be a positive element with $\| x \| = 1.$
Let $(e_{j, k})_{1 \leq j, k \leq n}$ be the standard
system of matrix units for $M_n.$
\Wolog\  we may assume that there is a finite set $F_0 \S A$
such that $\| a \| \leq 1$ for all $a \in F_0$ and such that
\[
F = \{ e_{j, k} \otimes a
 \colon {\mbox{$1 \leq j, k \leq n$ and $a \in F_0$}} \}.
\]

Set
$\ep_0 = \min \left( \ep / 4, \, 1 / (4 n^2) \right).$
Choose $\ep_1 > 0$ with $\ep_1 \leq \ep_0$
and so small that if $B$ is a \ca\  %
and $e, q \in B$ are \pj s such that $\| e q - q e \| < \ep_1,$
then there is a \pj\  $p \in B$ such that $\| e q e - p \| < \ep_0.$
Apply Lemma~\ref{PjAndNorm},
obtaining a nonzero \pj\  $p \in {\overline{x (M_n \otimes A) x}}$
such that $\| p_0 x p_0 - p_0 \| < \ep_0$
for every nonzero \pj\  $p_0 \leq p.$

Write $p = \sum_{j, k = 1}^n e_{j, k} \otimes p_{j, k}.$
By Lemma~\ref{PosMat},
there is $k$ such that $\ld = \| p_{k, k} \| \geq n^{-2}.$
\Wolog\  $k = 1.$
Lemma~\ref{L:MnSP} provides a nonzero \pj\  $d \in A$
such that $1 \otimes d \precsim p.$
Then $q = e_{1, 1} \otimes d$
is a nonzero \pj\  in $\C e_{1, 1} \otimes A$
such that $q \precsim p.$
By Lemma~\ref{OrthInSP},
there exist nonzero \mvnt\  \mops\  $q_1, q_2, \ldots, q_n \leq q.$
Set $y = \ld^{-1} p_{1, 1}.$
Apply Lemma~\ref{PjAndNorm},
obtaining a nonzero \pj\  $z \in {\overline{y A y}}$
such that $\| q_0 y q_0 - q_0 \| < \tfrac{1}{8}$
for every nonzero \pj\  $q_0 \leq z.$
Using Property~(SP) and Lemma~\ref{L:CompSP},
choose a nonzero \pj\  $q_0 \leq z$ such that
$e_{1, 1} \otimes q_0 \precsim q_1.$

Set
\[
\dt = \min \left( \ep, \, \frac{\ep_1}{n^2 \card (G)},
        \, \frac{1}{4 n}, \, \frac{1}{8 \card (G)} \right).
\]
Apply Definition~\ref{NewTRPDfn} to $\af,$
with
$F_0 \cup \{ q_0 \} \cup \{ p_{j, k} \colon 1 \leq j, k \leq n \}$
in place of $F,$
with $\dt$ in place of $\ep,$
and with $q_0$ in place of $x.$
Let $f_g,$ for $g \in G,$ be the resulting \pj s.
As usual, set $f = \sum_{g \in G} f_g.$
Then define $e_g = 1 \otimes f_g.$
Since $\dt \leq \ep,$ it is immediate that
$\| (\Ad (v_g) \otimes \af_g) (e_h) - e_{g h} \| < \ep$
for all $g, h \in G,$
and that
$\| [ e_g, \, (e_{j, k} \otimes a) ] \| < \ep$
for all $g \in G$ and all $a \in F_0.$
Moreover,
with $e = \sum_{g \in G} e_g,$ we have
\[
1 - e
 = \sum_{k = 1}^n e_{k, k} \otimes (1 - f)
 \precsim \sum_{k = 1}^n e_{k, k} \otimes q_0
 \precsim \sum_{k = 1}^n q_k
 \leq q,
\]
which is \mvnt\  to a \pj\  in ${\overline{x (M_n \otimes A) x}}.$

It remains only to show that $\| e x e \| > 1 - \ep.$
We have $\| [f, p_{j, k}] \| < \ep_1 / n^2$
for $1 \leq j, k \leq n,$
so $\| [e, p] \| < \ep_1.$
Therefore there is a \pj\  $r \in M_n \otimes A$ such that
$\| r - e p e \| < \ep_0.$
Now, assuming $r \neq 0$ at the last step,
\[
\| e x e \|
  \geq \| p e x e p \|
  \geq \| e p x p e \| - 2 \ep_1
  \geq \| e p e \| - 2 \ep_1 - \ep_0
  > \| r \| - 2 \ep_1 - 2 \ep_0
  \geq 1 - \ep.
\]
So we need only show that $r \neq 0.$
It suffices to show that $\| e p e \| > \ep_0.$

Using $\| f q_0 f \| > 1 - \dt$ at the seventh step, we have
\begin{align*}
\| e p e \|
& \geq \| f p_{1, 1} f \|
  = \ld \| f y f \|
  \geq \tfrac{1}{n^2} \| f y f \|
  \geq \tfrac{1}{n^2} \| q_0 f y f q_0 \|             \\
& \geq \tfrac{1}{n^2} ( \| f q_0 y q_0 f \| - 2 \card (G) \dt)
  > \tfrac{1}{n^2} \left( \| f q_0 f \|
             - \tfrac{1}{8} - \tfrac{1}{8} \right)     \\
& > \tfrac{1}{n^2} \left( 1 - \dt - \tfrac{1}{4} \right)
  \geq \frac{1}{4 n^2}
  \geq \ep_0.
\end{align*}
This completes the proof that
$g \mapsto \Ad (v_g) \otimes \af_g$ has the \tRp.

Now suppose that $g \mapsto \Ad (v_g) \otimes \af_g$ has the \tRp.
Let $g \mapsto w_g$ be the contragredient representation.
By what we already did,
$g \mapsto \Ad (w_g \otimes v_g) \otimes \af_g,$
which is an action of $G$ on $M_{n^2} \otimes A,$ has the \tRp.
The one dimensional trivial representation of $G$
is a subrepresentation of $g \mapsto w_g \otimes v_g,$
so there is an invariant \pj\  $p \in M_{n^2}$
such that the restriction of
$\Ad (w_g \otimes v_g) \otimes \af_g$
to $(p \otimes 1) (M_{n^2} \otimes A) (p \otimes 1) \cong A$
can be identified with $\af.$
Now Lemma~\ref{TRPCorner} implies that $\af$ has the \tRp.
\end{proof}

\begin{lem}\label{SPAndInv}
Let $A$ be a simple \ca\  with Property~(SP), and
let $\af \colon G \to \Aut (A)$
be an action of a finite group $G$ on $A.$
For each $g \in G,$ let $p_g \in A$ be a nonzero \pj.
Then there exists a nonzero \pj\  $q \in A$ such that
$q \precsim \af_g (p_g)$ for all $g \in G.$
\end{lem}

\begin{proof}
Write $G = \{ g_j \colon 0 \leq j \leq n - 1 \},$ with $g_0 = 1.$
Using Property~(SP) and Lemma~\ref{L:CompSP},
find a nonzero \pj\  $e_1 \leq p_1$
such that $e_1 \precsim \af_{g_1} (p_{g_1}).$
In the same way, find a nonzero \pj\  $e_2 \leq e_1$
such that $e_2 \precsim \af_{g_2} (p_{g_2}).$
Proceed inductively.
Set $q = e_{n - 1}.$
Then $q$ is a nonzero \pj\  such that
$q \precsim \af_{g_j} (p_{g_j})$ for $0 \leq j \leq n - 1,$
whence $q \precsim \af_g (p_g)$ for all $g \in G.$
\end{proof}

The next result is the analog of Lemma~3.8 of~\cite{Iz}.

\begin{thm}\label{TRPDualToTAppRep}
Let $A$ be an \idssuca,
and let $\af \colon G \to \Aut (A)$
be an action of a finite abelian group $G$ on $A$
such that $C^* (G, A, \af)$ is also simple.
Then:
\begin{enumerate}
\item\label{TRPDualToTAppRep:1} %
$\af$ has the \tRp\  \ifo\  %
${\widehat{\af}}$ is tracially approximately representable.
\item\label{TRPDualToTAppRep:2} %
$\af$ is tracially approximately representable \ifo\  %
${\widehat{\af}}$ has the \tRp.
\end{enumerate}
\end{thm}

\begin{proof}
Using Proposition~\ref{SPForCrPrd},
Lemma~\ref{L:MnSP},
and Takai duality~\cite{Tk},
we see that
that either both $A$ and $C^* (G, A, \af)$ have Property~(SP)
or neither does.
If neither does,
by Remark~\ref{SRPImpTRP} and Lemma~\ref{RImpSP},
the statement is equivalent to the
corresponding statement for the strict Rokhlin property
and approximate representability, and follows from
Lemma~3.8 of~\cite{Iz}.
Thus, we assume both have Property~(SP).

We prove~(\ref{TRPDualToTAppRep:2}).
Statement~(\ref{TRPDualToTAppRep:1}) will then follow
from~(\ref{TRPDualToTAppRep:2}) for ${\widehat{\af}},$
by combining Takai duality~\cite{Tk} with Lemma~\ref{TRPMat}.

For $g \in G$ let $u_g \in C^* (G, A, \af)$
be the standard unitary of the crossed product.
Also set $n = \card (G).$

Suppose that $\af$ is tracially approximately representable.
Let $F \S C^* (G, A, \af)$ be finite,
let $\ep > 0,$
and let $x \in C^* (G, A, \af)$
be a positive element with $\| x \| = 1.$
\Wolog\  $F = F_0 \cup \{ u_g \colon g \in G \}$
for some finite subset $F_0 \S A,$
such that $\| a \| \leq 1$ for all $a \in F_0.$

Use Proposition~\ref{SPForCrPrd}
and Lemma~\ref{PjAndNorm} to find
a nonzero \pj\  $q_0$
in the hereditary subalgebra ${\overline{x C^* (G, A, \af) x}}$
such that whenever $p \leq q_0$ is a nonzero \pj,
then $\| p x p \| \geq 1 - \frac{1}{8} \ep.$
By Proposition~\ref{SPForCrPrd} again,
there is a nonzero \pj\  $q \in A$
and a partial isometry $c \in C^* (G, A, \af)$
such that $c c^* = q$ and $c^* c \leq q_0.$
Write $c = \sum_{g \in G} c_g u_g$ with $c_g \in A.$
Note that $\| c_g \| \leq \| c \| = 1.$

Set
$\ep_0 = \ep / (8 n).$
Apply Lemma~\ref{STrAppRep}
with $F_0 \cup \{ c_g \colon g \in G \}$ in place of $F,$
with $\ep_0$ in place of $\ep,$
and with $q$ in place of $x.$
Let $e \in A$ and $w_g \in e A e$ be the resulting \pj\  and unitaries.
For $\sm \in {\widehat{G}},$
define $e_{\sm} \in C^* (G, A, \af)$ by
\[
e_{\sm} = \frac{1}{n} \sum_{g \in G} \sm (g) u_g w_g^*.
\]
We show that the $e_{\sm}$ are \pj s which verify the definition of
the \tRp.

The unitaries $u_g \in C^* (G, A, \af)$ and $w_h \in e A e$
commute for all $g, h \in G,$
because $\af_g (w_h) = w_h.$
It is now easy to check that $e_{\sm}^* = e_{\sm}.$
Next, let $\sm, \ta \in {\widehat{G}}.$
Then, changing the variable $g$ to $g h^{-1}$ at the second step,
\[
e_{\sm} e_{\ta}
  = \frac{1}{n^2}
      \sum_{g, h \in G} \sm (g) \ta (h) u_{g h} w_{g h}^*    \\
  = \frac{1}{n^2}
      \sum_{g \in G} \left( \ssum{h \in G} (\sm^{-1} \ta) (h) \right)
          \sm (g) u_{g} w_{g}^*.
\]
We have $\sum_{h \in G} (\sm^{-1} \ta) (h) = n$ if $\sm = \ta,$
and the sum is zero otherwise, so the $e_{\sm}$ are \mops.
Moreover,
\[
\sum_{\sm \in {\widehat{G}}} e_{\sm}
 = u_1 w_1^*
 = e.
\]
For $\ta \in {\widehat{G}},$
we have ${\widehat{\af}}_{\ta} (w_g) = w_g$
and ${\widehat{\af}}_{\ta} (u_g) = \ta (g) u_g.$
Therefore ${\widehat{\af}}_{\ta} (e_{\sm}) = e_{\ta \sm}$
for $\sm, \ta \in {\widehat{G}}.$

Next, we check the approximate commutation relations.
For $g \in G,$ because $G$ is abelian and $u_g$ commutes
with each $w_h,$
we have $e_{\ta} u_g = u_g e_{\ta}$ for all $\ta \in {\widehat{G}}.$
For $a \in F_0$ and $\ta \in {\widehat{G}},$
we use the relation $u_g e a e u_g^* = \af_g (e a e)$
and the estimate
$\| \af_g (e a e) - w_g e a e w_g^* \| < \ep_0$
to get
\begin{align*}
\| e_{\ta} a - a e_{\ta} \|
& \leq 2 \| e a - a e \| + \| e_{\ta} e a e - e a e e_{\ta} \|    \\
& \leq 2 \| e a - a e \|
    + \frac{1}{n} \sum_{g \in G}
        \big\| \sm (g) u_g w_g^* e a e - e a e \sm (g) u_g w_g^* \big\|
                 \\
& < 2 \ep_0 + \ep_0
  \leq \ep.    
\end{align*}

Next, we observe that $1 - e \precsim q,$
so $1 - e$ is \mvnt\  %
to a \pj\  in ${\overline{x C^* (G, A, \af) x}}.$

Finally, we estimate $\| e x e \|.$
First, since $e$ is $\af$-invariant,
we have $u_g e = e u_g$ for all $g \in G.$
Therefore
\[
\| c e - e c \|
 \leq \sum_{g \in G} \| c_g e - e c_g \|
 < \tfrac{1}{8} \ep,
\]
whence $\| c^* c e - e c^* c \| < \tfrac{1}{4} \ep.$
Using the fact that $c^* c$ is a \pj\  and $c^* c \leq q_0$
at the fourth step,
we get
\begin{align*}
\| e x e \|
& \geq \| c^* c e x e c^* c \|
  \geq \| e c^* c x c^* c e \| - \tfrac{1}{2} \ep
  > \| e c^* c e \| - \tfrac{5}{8} \ep
  = \| (c e) (c e)^* \| - \tfrac{5}{8} \ep      \\
& > \| e c c^* e \| - \tfrac{7}{8} \ep
  = \| e q e \| - \tfrac{7}{8} \ep
  > 1 - \ep_0 - \tfrac{7}{8} \ep
  \geq 1 - \ep.
\end{align*}
This completes the proof that ${\widehat{\af}}$ has the \tRp.

Now assume that ${\widehat{\af}}$ has the \tRp.
Let $F \S A$ be a finite set, let $\ep > 0,$
and let $x \in A$ be a positive element with $\| x \| = 1.$
\Wolog\  $\| a \| \leq 1$ for all $a \in F.$
Apply Lemma~\ref{PjAndNorm},
obtaining a nonzero \pj\  $q \in {\overline{x A x}}$
such that $\| r x r - r \| < \frac{1}{2} \ep$
for every nonzero \pj\  $r \leq q.$
Use Lemma~\ref{OrthInSP} to find a family $(q_g)_{g \in G}$
of nonzero \mops\  in ${\overline{q A q}}.$
By Lemma~\ref{SPAndInv}, there exists a nonzero \pj\  $p \in A$ such that
$p \precsim \af_g (q_g)$ for all $g \in G.$
Therefore there are nonzero \pj s $p_g \in A$
with $p_g \leq q_g$ such that
$p_g \sim \af_g (p_1)$ for all $g \in G.$

Choose $\ep_0 > 0$ so small that if $B$ is a unital \ca\  and
$y \in B$ satisfies
\[
\| y y^* - 1 \| < 5 \ep_0 n
\andeqn
\| y^* y - 1 \| < 5 \ep_0 n,
\]
then there is a unitary $w \in B$
such that $\| w - y \| < \frac{1}{7} \ep.$
We also require that
$9 n^2 \ep_0 \leq \frac{1}{7} \ep.$

Apply Lemma~\ref{StTRPDfn} to ${\widehat{\af}}$
with $F \cup \{ p_1 \} \cup  \{ u_g \colon g \in G \}$ in place of $F,$
with $\ep_0$ in place of $\ep,$
and with $p_1$ in place of $x,$
obtaining \pj s $e_{\ta}$ for $\ta \in {\widehat{G}}.$
Let $e = \sum_{\ta \in {\widehat{G}} } e_{\ta}.$
Note that $e \in A$ because $e$ is ${\widehat{\af}}$-invariant.

Define
\[
x_g = \sum_{\sm \in {\widehat{G}} } \sm (g) e_{\sm} u_g e_{\sm}
\]
for $g \in G.$
Note that $x_g \in e C^* (G, A, \af) e.$
We estimate as follows.
For $g \in G$ and $\ta \in {\widehat{G}},$
\[
\| {\widehat{\af}}_{\ta} (x_g) - x_g \|
  \leq \sum_{\sm \in {\widehat{G}} }
           \big\| \sm (g) \ta (g) {\widehat{\af}}_{\ta} (e_{\sm})
                      u_g {\widehat{\af}}_{\ta} (e_{\sm})
                - (\sm \ta) (g) e_{\ta \sm} u_g e_{\ta \sm} \big\|
  < 2 n \ep_0.
\]
For $g \in G,$
\begin{align*}
\| x_g x_g^* - e \|
& = \left\| \ssum{\sm, \ta \in {\widehat{G}} }
        \sm (g) \ta (g)^{-1} e_{\sm} u_g e_{\sm} e_{\ta} u_g^* e_{\ta}
      - \ssum{\sm \in {\widehat{G}} } e_{\sm} \right\|
                 \\
& \leq \sum_{\sm \in {\widehat{G}} }
           \| e_{\sm} u_g e_{\sm} u_g^* e_{\sm} - e_{\sm} \|
  < n \ep_0.
\end{align*}
Similarly,
\[
\| x_g^* x_g - e \| < n \ep_0.
\]
For $a \in F \cup \{ p_1 \},$
\[
\| e a - a e \|
 \leq \sum_{\ta \in {\widehat{G}} } \| e_{\ta} a - a e_{\ta} \|
 < n \ep_0
 \leq \tfrac{1}{7} \ep.
\]
For $g \in G$ and $a \in F,$
\begin{align*}
\| u_g e a e u_g^* - x_g e a e x_g^* \|
& \leq \sum_{\sm, \ta \in {\widehat{G}} }
 \big\| \sm (g) \ta (g)^{-1} e_{\sm} u_g e_{\sm} a e_{\ta} u_g^* e_{\ta}
               - u_g e_{\sm} a e_{\ta} u_g^* \big\|     \\
& < 9 n^2 \ep_0
    + \sum_{\sm, \ta \in {\widehat{G}} }
       \big\| \sm (g) \ta (g)^{-1} u_g a u_g^* e_{\sm}^2 e_{\ta}^2
               - u_g a u_g^* e_{\sm} e_{\ta} \big\|
                        \\
& = 9 n^2 \ep_0
  \leq \tfrac{1}{7} \ep.
\end{align*}
For $g, h \in G,$
\begin{align*}
\| x_g x_h - x_{g h} \|
& = \left\| \ssum{\sm, \ta \in {\widehat{G}} }
        \sm (g) \ta (h) e_{\sm} u_g e_{\sm} e_{\ta} u_h e_{\ta}
     - \ssum{\sm \in {\widehat{G}} }
        \sm (g h) e_{\sm} u_g u_h e_{\sm} \right\|     \\
& \leq \sum_{\sm \in {\widehat{G}} } \| u_g e_{\sm} - e_{\sm} u_g \|
  < \ep_0 n
 \leq \tfrac{1}{7} \ep.
\end{align*}
For $g, h \in G,$
\begin{align*}
\| u_g x_h u_g^* - x_h \|
& \leq \sum_{\sm \in {\widehat{G}} }
      \| u_g e_{\sm} u_h e_{\sm} u_g^* - e_{\sm} u_h e_{\sm} \|     \\
& \leq 2 n \| u_g e_{\sm} - e_{\sm} u_g \|
  < 2 n \ep_0
  \leq \tfrac{1}{7} \ep.
\end{align*}

Now for $g \in G$ set
\[
y_g = \frac{1}{n}
   \sum_{\ta \in {\widehat{G}} } {\widehat{\af}}_{\ta} (x_g).
\]
Then $y_g \in e C^* (G, A, \af) e$
and $y_g$ is ${\widehat{\af}}$-invariant,
so $y_g \in e A e.$
The first estimate in the previous paragraph implies that
$\| y_g - x_g \| < 2 n \ep_0.$
The next two then imply that
\[
\| y_g y_g^* - e \| < 5 n \ep_0
\andeqn
\| y_g^* y_g - e \| < 5 n \ep_0.
\]
Therefore there are unitaries $w_g \in e A e$ such that
$\| w_g - y_g \| < \tfrac{1}{7} \ep.$
It follows that
$\| w_g - x_g \|
 < \tfrac{1}{7} \ep + 2 n \ep_0 \leq \tfrac{2}{7} \ep.$
The remaining four estimates in the previous paragraph now imply,
in order,
Conditions~(\ref{TrAppRepDfn:1}) through~(\ref{TrAppRepDfn:4})
in Definition~\ref{TrAppRepDfn}.
For the first, this is immediate.
For the second, for $g \in G$ and $a \in F,$
\begin{align*}
\| \af_g (e a e) - w_g e a e w_g^* \|
& = \| u_g e a e u_g^* - w_g e a e w_g^* \|          \\
& \leq 2 \| w_g - x_g \| + \| u_g e a e u_g^* - x_g e a e x_g^* \|
  < \tfrac{5}{7} \ep
  \leq \ep.
\end{align*}
Similarly, the remaining two give
$\| w_g w_h - w_{g h} \| < \ep$
and $\| \af_g (w_h) - w_h \| < \tfrac{5}{7} \ep \leq \ep.$

We next verify Condition~(\ref{TrAppRepDfn:5})
in Definition~\ref{TrAppRepDfn}.
Our choices give $1 - e \precsim p_1$ in $C^* (G, A, \af).$
This is therefore also true in
$C^* \big( {\widehat{G}}, \, C^* (G, A, \af), \, {\widehat{\af}} \big).$
Use the notation of Proposition~\ref{Takai} and the identification there
of
$C^* \big( {\widehat{G}}, \, C^* (G, A, \af), \, {\widehat{\af}} \big)$
with $L (l^2 (G)) \otimes A.$
By Proposition~\ref{Takai},
in $L (l^2 (G)) \otimes A$ we have
\[
\sum_{g \in G} e_{g, g} \otimes \af_g^{-1} (1 - e)
 \precsim \sum_{g \in G} e_{g, g} \otimes \af_g^{-1} (p_1).
\]
Thus, in $L (l^2 (G)) \otimes A$ we have
\[
e_{1, 1} \otimes (1 - e)
 \leq \sum_{g \in G} e_{g, g} \otimes \af_g^{-1} (1 - e)
 \precsim \sum_{g \in G} e_{g, g} \otimes \af_g^{-1} (p_1)
 \sim \sum_{g \in G} e_{1, 1} \otimes p_g
 \leq e_{1, 1} \otimes q.
\]
Therefore $1 - e \precsim q$ in $A,$
which verifies condition~(\ref{TrAppRepDfn:5}).

It remains to verify Condition~(\ref{TrAppRepDfn:6}).
Since $\ep_0 \leq \frac{1}{7} \ep,$
we have $\| e p_1 e \| > 1 - \frac{1}{7} \ep.$
Since $p_1 \leq q$ and is nonzero,
we have $\| p_1 - p_1 x p_1 \| < \frac{1}{2} \ep,$
whence
\begin{align*}
\| e x e \|
 & \geq \| p_1 e x e p_1 \|
   \geq \| e p_1 e \| - 2 \| [ e, p_1 ] \| - \| p_1 - p_1 x p_1 \|
               \\
 & > 1 - \tfrac{1}{7} \ep - \tfrac{2}{7} \ep - \tfrac{1}{2} \ep
   > 1 - \ep.
\end{align*}
This completes the proof.
\end{proof}

\section{Strongly tracially approximately inner
   automorphisms}\label{Sec:STAI}

\indent
In this section, we introduce the notion of a
strongly tracially approximately inner automorphism.
The main result of this section is that if an action
of a finite cyclic group has the \tRp\  and is
generated by a strongly tracially approximately inner automorphism,
then the action is tracially approximately representable
in the sense of Definition~\ref{TrAppRepDfn}.
It will follow that the dual action has the \tRp.

We begin by proving the ``nontracial'' version,
essentially, that for finite cyclic groups,
approximate innerness and the \sRp\  imply
approximate representability.
This result is not particularly useful,
because approximate representability seems usually to be
easier to prove than the \sRp.
(See the discussion in at the beginning of Section~3.2 of~\cite{Iz},
where approximate representability
is seen primarily as a way to get actions with the \sRp\  %
by duality.)
However, the proof of the tracial version
will partially follow this proof.
The tracial version is much more useful,
because the \tRp\  can often be proved directly.

\begin{prp}\label{AppInnRP}
Let $A$ be a separable unital \ca,
let $\af \in \Aut (A)$ be approximately inner
and satisfy $\af^n = \id_A,$
and suppose that the action of $\Zqn$ generated by $\af$
has the \sRp.
Then the action of $\Zqn$ generated by $\af$
is approximately representable.
\end{prp}

\begin{proof}
The proof has two steps.
First, we show that one can choose the unitaries
in the definition of approximate innerness to have order $n,$
and then we show that one can in addition choose them to
be $\af$-invariant.

Step~1: We claim that for every finite subset $F \S A$
and every $\ep > 0,$ there is a unitary $v \in A$ such that
$\| v a v^* - \af (a) \| < \ep$ for all $a \in F,$
and such that $v^n = 1.$

To prove this, \wolog\  $F$ is $\af$-invariant
and $\| a \| \leq 1$ for all $a \in F.$
Choose $\dt > 0$ with
\[
\dt \leq \frac{\ep}{3 (2 n + 7) n},
\]
and also (using semiprojectivity of $\C^n$;
see Lemma~14.1.5, Theorem~14.2.1, Theorem~14.1.4,
and Definition~14.1.1 of~\cite{Lr})
so small that whenever $B$ is a unital \ca\  and $c \in B$
satisfies
\[
\| c^* c - 1 \| < 2 (n - 1) \dt, \,\,\,\,\,\,
\| c c^* - 1 \| < 2 (n - 1) \dt, \andeqn
\| c^n - 1 \| < 2 (n - 1) \dt,
\]
then there is a unitary $v \in B$ such that $v^n = 1$
and $\| v - c \| < \frac{1}{3} \ep.$

Use the \sRp\  to find \mops\  $f_0, f_1, \ldots, f_{n - 1} \in A$
such that $\sum_{j = 0}^{n - 1} f_j = 1,$
such that $\| \af^k (f_j) - f_{j + k} \| < \dt$
for $0 \leq j, k \leq n - 1$ (with the subscripts taken cyclically),
and such that $\| [ e_j, a ] \| < \dt$
for $0 \leq j \leq n - 1$ and $a \in F.$
Then choose a unitary $u \in A$ such that whenever
\[
b \in \{ f_0, f_1, \ldots, f_{n - 1} \}
        \cup \{ f_j a f_j \colon
         {\mbox{$0 \leq j \leq n - 1$ and $a \in F$}} \},
\]
we have $\| u b u^* - \af (b) \| < \dt.$

Set $w_k = f_{k + 1} u f_k$ for $0 \leq k \leq n - 2,$
and set
\[
w_{n - 1} = f_0 u^* f_1 u^* f_2 \cdots f_{n - 2} u^* f_{n - 1}.
\]
Then set $w = \sum_{k = 0}^{n - 1} w_k.$

We estimate $\| w^* w - 1 \|$ and $\| w w^* - 1 \|.$
First, observe that for $0 \leq k \leq n - 2$ we have
\begin{align*}
\| w_k^* w_k - f_k \|
& \leq \| u^* f_{k + 1} u - f_k \|
  = \| f_{k + 1} - u f_k u^* \|              \\
& \leq \| u f_k u^* - \af (f_k) \| + \| \af (f_k) - f_{k + 1} \|
  < 2 \dt
\end{align*}
and similarly $\| w_k^* w_k - f_{k + 1} \| < 2 \dt.$
By downwards induction on $k,$ we also get
\[
\big\| (f_k u^* f_{k + 1} \cdots f_{n - 2} u^* f_{n - 1})
        (f_k u^* f_{k + 1} \cdots f_{n - 2} u^* f_{n - 1})^*
      - f_k \big\|
   < 2 (n - k - 1) \dt,
\]
so that $\| w_{n - 1} w_{n - 1}^* - f_{0} \| < 2 (n - 1) \dt.$
Similarly, $\| w_{n - 1}^* w_{n - 1} - f_{n - 1} \| < 2 (n - 1) \dt.$
Since the $f_k$ are orthogonal, we thus get
\[
\| w w^* - 1 \|
 \leq \max_{0 \leq k \leq n - 1} \| w_k^* w_k - f_k \|
 < 2 (n - 1) \dt.
\]
Similarly, $\| w^* w - 1 \| < 2 (n - 1) \dt.$

Now we estimate $\| w^n - 1 \|.$
Since the $e_k$ are orthogonal, we have
\begin{align*}
& \| w^n - 1 \|                    \\
& \hspace*{2em}
 = \max_{0 \leq k \leq n - 1}
         \big\| w_{k - 1} w_{k - 2} \cdots w_0 w_{n - 1} \cdots w_{k}
                      - f_k \big\|
                                   \\
& \hspace*{2em}
 = \max_{0 \leq k \leq n - 1}
         \big\| [ f_k u f_{k - 1} \cdots f_1 u f_0 u^* f_1 \cdots
                   f_{n - 2} u^* f_{n - 1} u f_{n - 2} \cdots
                   f_{k + 1} u f_k ]
                - f_k \big\|.
\end{align*}
Using the inequalities
\[
\| u f_k u^* - f_{k + 1} \|
 \leq \| u f_k u^* - \af (f_k) \| + \| \af (f_0) - f_{k + 1} \|
 < 2 \dt
\]
for $0 \leq k \leq n - 1$ (indices taken cyclically),
and similarly $\| u^* f_{k} u - f_{k - 1} \| < 2 \dt,$
one gets by induction
\[
\big\| f_k u f_{k - 1} \cdots f_1 u f_0 u^* f_1 \cdots
                   f_{n - 2} u^* f_{n - 1} u f_{n - 2} \cdots
                   f_{k + 1} u f_k 
                - f_k \big\|
         < 2 (n - 1) \dt
\]
for all $k.$
So $\| w^n - 1 \| < 2 (n - 1) \dt.$
Thus, by the choice of $\dt,$
there exists a unitary $v \in A$ such that
$\| v - w \| < \frac{1}{3} \ep$ and $v^n = 1.$

We next estimate $\| w a w^* - \af (a) \|$ for $a \in F.$
Set $b = \sum_{k = 0}^{n - 1} f_k a f_k$
and $c = \sum_{k = 0}^{n - 1} f_k \af (a) f_k.$
Then
\[
\| a - b \|
  \leq \sum_{k = 0}^{n - 1}  \sum_{j \neq k} \| f_j a f_k \|
  \leq \sum_{k = 0}^{n - 1}  \sum_{j \neq k} \| [f_k, a] \|
  < n (n - 1) \dt.
\]
Similarly, recalling that $F$ is $\af$-invariant,
$\| \af (a) - c \| < n (n - 1) \dt.$
For $0 \leq k \leq n - 2,$ we have
\begin{align*}
& \| w_k f_k a f_k w_k^* - f_{k + 1} \af (a) f_{k + 1} \|
               \\
& \mbox{} \hspace*{3em}
  \leq \big\| f_{k + 1} u f_k a f_k u^* f_{k + 1}
         - f_{k + 1} \af (f_k a f_k) f_{k + 1} \big\|
       + 2 \| \af (f_k) - f_{k + 1} \|
  < 3 \dt.
\end{align*}

For $k = n - 1,$ we have to estimate
\begin{align*}
&
\big\| w_{n - 1} f_{n - 1} a f_{n - 1} w_{n - 1}^*
                                   - f_{0} \af (a) f_{0} \big\|
               \\
& \mbox{} \hspace*{3em}
  \leq \big\| f_{n - 1} a f_{n - 1}
        - w_{n - 1}^* f_{0} \af (a) f_{0} w_{n - 1} \big\|
    + 2 \| w_{n - 1} w_{n - 1}^* - f_0 \|.
\end{align*}
The last term was shown above to be less than $4 (n - 1) \dt.$
For the first term, we claim that
\[
\big\| (f_k u f_{k - 1} \cdots f_1 u f_0) \af (a)
            (f_k u f_{k - 1} \cdots f_1 u f_0)^*
         - f_k \af^{k + 1} (a) f_k \big\|
     < 3 k \dt,
\]
and we prove this by induction on $k.$
For $k = 0$ is it trivial.
If it is true for $k,$
then the estimate
\begin{align*}
& \big\| f_{k + 1} u f_k \af^{k + 1} (a) f_k u^* f_{k + 1}
      - f_{k + 1} \af^{k + 2} (a) f_{k + 1} \big\|            \\
& \mbox{} \hspace*{3em}
   \leq 2 \| f_{k + 1} - \af (f_k) \|
           + \big\| u f_k \af^{k + 1} (a) f_k u^*
                 - \af (f_{k} \af^{k + 1} (a) f_{k} ) \big\|
  < 3 \dt.
\end{align*}
This implies that it is true for $k + 1,$ completing the induction.
Since $\af^n = \id_A,$
we get
\[
\big\| w_{n - 1} f_{n - 1} a f_{n - 1} w_{n - 1}^*
                   - f_{0} \af (a) f_{0} \big\|
   < [3 (n - 1) + 4 (n - 1)] \dt.
\]
By orthogonality of the summands,
$\| w b w^* - c \| < 7 (n - 1) \dt,$
so $\| w a w^* - \af (a) \| < (2 n + 7) n \dt.$
It follows that
\[
\| v a v^* - \af (a) \|
 < \tfrac{2}{3} \ep + (2 n + 7) n \dt
 \leq \ep.
\]
This completes the proof of Step~1.

Step~2: We prove the conclusion.
Let $F \S A$ be finite, and let $\ep > 0.$
For the same reasons as in Step~1,
we may choose $\dt > 0$ with
\[
\dt \leq \frac{\ep}{3 (2 n^2 + 2 n + 3)},
\]
and also
so small that whenever $B$ is a unital \ca\  and $c \in B$
satisfies
\[
\| c^* c - 1 \| < (2 n^2 + n) \dt, \,\,\,\,\,\,
\| c c^* - 1 \| < (2 n^2 + n) \dt, \andeqn
\| c^n - 1 \| < (2 n^2 + n) \dt,
\]
then there is a unitary $v \in B$ such that $v^n = 1$
and $\| v - c \| < \frac{1}{3} \ep.$

Apply the result of Step~1
with $F$ as given and with $\dt$ in place of $\ep,$
obtaining a unitary $u.$
Apply the \sRp\  with
$\bigcup_{k = 0}^{n - 1} \af^k (F \cup \{ u \} )$ in place of $F$
and with $\dt$ in place of $\ep,$
obtaining \mops\  $f_0, f_1, \ldots, f_{n - 1} \in A.$

Set $w_k = f_k \af^k (u) f_k$ for $0 \leq k \leq n - 1,$
and set $w = \sum_{k = 0}^{n - 1} w_k.$
We have
\[
\| w_k w_k^* - f_k \| \leq \| [ \af^k (u), \, f_k ] \| < \dt,
\]
and similarly $\| w_k^* w_k - f_k \| < \dt.$
Using orthogonality of the $f_k,$ we then get
\[
\| w w^* - 1 \|
  = \max_{0 \leq k \leq n - 1} \| w_k w_k^* - f_k \| < \dt,
\]
and similarly $\| w^* w - 1 \| < \dt.$
Moreover,
\[
\| w_k^n - f_k \| \leq (n - 1) \| [ \af^k (u), \, f_k ] \|
 < (n - 1) \dt,
\]
so, using orthogonality again,
$\| w^n - 1 \| < (n - 1) \dt.$

Now define $y_k = \af^k (f_0 u f_0)$ for $0 \leq k \leq n - 1,$
and set $y = \sum_{k = 0}^{n - 1} y_k.$
Then $\| y_k - w_k \| \leq 2 \| f_k - \af^k (f_0) \| < 2 \dt,$
so $\| y - w \| < 2 n \dt.$
Therefore
\[
\| y^* y - 1 \| \leq \| w^* w - 1 \| + 2 \| y - w \|
  < (2 n + 1) \dt,
\]
and similarly $\| y y^* - 1 \| < (2 n + 1) \dt$
and $\| y^n - 1 \| < [n \cdot 2 n + (n - 1)] \dt.$
We have $\af (y) = y$ because $\af^n = \id_A.$
Apply the choice of $\dt$ with the fixed point algebra $A^{\af}$
in place of $B,$
obtaining an $\af$-invariant unitary $v \in A$ such that
$\| v - y \| < \frac{1}{3} \ep$ and $v^n = 1.$
Then $\| v - w \| < \frac{1}{3} \ep + 2 n \dt.$

It remains to estimate
$\| v a v^* - \af (a) \|$ for $a \in F.$
Set $b = \sum_{k = 0}^{n - 1} f_k a f_k$ and
$c = \sum_{k = 0}^{n - 1} f_k \af (a) f_k.$
As in the proof of Step~1,
we have $\| a - b \| < n (n - 1) \dt$
and $\| \af (a) - c \| < n (n - 1) \dt.$
Also,
\begin{align*}
\| w_k f_k a f_k w_k^* - f_k \af (a) f_k \|
 & = \big\| f_k \af^k (u) f_k a f_k \af^k (u)^* f_k
                     - f_k \af (a) f_k \big\|
             \\
 & = \big\| f_k \af^k (u \af^{- k} (a) u^*) f_k - f_k \af (a) f_k \big\|
        + 2 \| [ \af^k (u), f_k ] \|
             \\
 & \leq \| u \af^{- k} (a) u^* - \af^{- k + 1} (a) \|
        + 2 \| [ \af^k (u), f_k ] \|
   < 3 \dt.
\end{align*}
Using orthogonality of the $f_k,$ we then get
$\| w b w^* - c \| < 3 \dt,$
whence $\| w a w^* - \af (a) \| < [2 n (n - 1) + 3] \dt.$
Therefore
\[
\| v a v^* - \af (a) \|
  < 2 \left( \tfrac{1}{3} \ep + 2 n \dt \right) + [2 n (n - 1) + 3] \dt
  \leq \ep.
\]
as was to be proved.
\end{proof}

We now define a strongly tracially approximately inner automorphism.
We call the condition ``strong tracial approximate innerness''
because there is a different condition,
Definition~\ref{D:TAI},
which is more appropriately called tracial approximate innerness.
In particular, it is unlikely to be true that the product of two
automorphisms which are strongly tracially approximately inner
as defined here is again strongly tracially approximately inner.
We have not been able to prove the results of this section
with tracial approximate innerness
in place of strong tracial approximate innerness.
However, the tracially approximately inner automorphisms
form a group (Theorem~\ref{T:TAIIsGp}),
and, moreover, if $A$ is an \idssuca\  with tracial rank zero,
then a tracially approximately inner automorphism
of finite order is
necessarily strongly tracially approximately inner
(combine Proposition~\ref{TAIAndInf} and Theorem~\ref{FOTAIOnTAF}).

\begin{dfn}\label{STAInnDfn}
Let $A$ be an \idssuca\  and let $\af \in \Aut (A).$
We say that $\af$ is
{\emph{strongly tracially approximately inner}}
if for every finite set $F \S A,$ every $\ep > 0,$
and every positive element $x \in A$ with $\| x \| = 1,$
there exist a \pj\  $e \in A$ and a unitary $v \in e A e$
such that:
\begin{enumerate}
\item\label{STAInnDfn:1} %
$\| \af (e) - e \| < \ep.$
\item\label{STAInnDfn:2} %
$\| e a - a e \| < \ep$ for all $a \in F.$
\item\label{STAInnDfn:3} %
$\| v e a e v^* - \af (e a e) \| < \ep$ for all $a \in F.$
\item\label{STAInnDfn:4} %
$1 - e$ is \mvnt\  to a \pj\  in ${\overline{x A x}}.$
\item\label{STAInnDfn:5} %
$\| e x e \| > 1 - \ep.$
\end{enumerate}
\end{dfn}

Our original definition,
in~\cite{PhW}, had a condition like that in
Remark~\ref{MoreSmallness} in place of~(\ref{STAInnDfn:5}).

As in Definition~\ref{NewTRPDfn},
we allow $e = 1,$ in which case
conditions~(\ref{STAInnDfn:4}) and~(\ref{STAInnDfn:5}) are vacuous.

\begin{rmk}\label{InnAndTAInn}
Let $A$ be an \idssuca\  and let $\af \in \Aut (A).$
If $\af$ is approximately inner
then $\af$ is strongly tracially approximately inner.
If $\af$ is strongly tracially approximately inner
and $A$ does not have Property~(SP),
then $\af$ is approximately inner.
\end{rmk}

Example~2.9 of~\cite{PhtRp1b}
shows that a strongly tracially approximately inner
automorphism need not be approximately inner,
even on a simple AF~algebra.
In fact, in Theorem~\ref{FOTAIOnTAF},
we give a condition which implies that many automorphisms
which are not inner are nevertheless
strongly tracially approximately inner.

When $A$ is finite,
we do not need Condition~(\ref{STAInnDfn:5})
of Definition~\ref{STAInnDfn}:

\begin{lem}\label{TAIForFinite}
Let $A$ be an \idfssuca\  and let $\af \in \Aut (A).$
Then $\af$ is strongly tracially approximately inner \ifo\  %
for every finite set $F \S A,$ every $\ep > 0,$
and every positive element $x \in A$ with $\| x \| = 1,$
there exist a \pj\  $e \in A$ and a unitary $v \in e A e$
such that:
\begin{enumerate}
\item\label{TAIForFinite:1} %
$\| \af (e) - e \| < \ep.$
\item\label{TAIForFinite:2} %
$\| e a - a e \| < \ep$ for all $a \in F.$
\item\label{TAIForFinite:3} %
$\| v e a e v^* - \af (e a e) \| < \ep$ for all $a \in F.$
\item\label{TAIForFinite:4} %
$1 - e$ is \mvnt\  to a \pj\  in ${\overline{x A x}}.$
\end{enumerate}
\end{lem}

\begin{proof}
The proof is the same as for Lemma~\ref{TRPForFinite}.
\end{proof}

We now prove the tracial analog of Proposition~\ref{AppInnRP}.
We separate the analog of Step~1 of its proof as a separate lemma.

\begin{lem}\label{FinOrdTAInn}
Let $A$ be an \idssuca\  with \PSP,
and let $\af \in \Aut (A)$
be strongly tracially approximately inner and satisfy $\af^n = \id_A.$
Suppose $\af$ generates an action of $\Zqn$ with the
tracial Rokhlin property.
Then for every finite set
$F \S A,$ every $\ep > 0,$
and every positive element $x \in A$ with $\| x \| = 1,$
there exist a \pj\  $e \in A$ and a unitary $v \in e A e$
such that:
\begin{enumerate}
\item\label{FinOrdTAInn:1} %
$\af (e) = e$ and $v^n = e.$
\item\label{FinOrdTAInn:2} %
$\| e a - a e \| < \ep$ for all $a \in F.$
\item\label{FinOrdTAInn:3} %
$\| v e a e v^* - \af (e a e) \| < \ep$ for all $a \in F.$
\item\label{FinOrdTAInn:4} %
$1 - e$ is \mvnt\  to a \pj\  in ${\overline{x A x}}.$
\item\label{FinOrdTAInn:5} %
$\| e x e \| > 1 - \ep.$
\end{enumerate}
\end{lem}

\begin{proof}
Let $F \S A$ be finite,
let $\ep > 0,$
and let $x \in A$ be a positive element with $\| x \| = 1.$
\Wolog\  $\| a \| \leq 1$ for all $a \in F.$

Set
\[
\rh = \min \left( \tfrac{1}{n}, \, \tfrac{1}{5} \right)
\andeqn
\ep_0 = \frac{\ep}{168 n^2}.
\]

Choose $\ep_1 > 0$ with 
\[
\ep_1 \leq \min \left( \tfrac{1}{36}, \, \ep_0 \right),
\]
and so small that the following are all true in any \ca\  $B$:
\begin{itemize}
\item
Whenever $r_1, r_2 \in B$ are \pj s,
and $c \in B$ satisfies $\| c^* c - r_1 \| < 27 \ep_1$
and $\| c c^* - r_2 \| < 27 \ep_1,$
then there is $v \in B$ such that $\| v - c \| < \ep_0,$
$v^* v = r_1,$ and $v v^* = r_2.$
\item
Whenever
$e, q \in B$ are \pj s such that $\| e q - q e \| < \ep_1,$
then there exists a \pj\  $r \leq q$ such that
$\| r - e q \| < \min \left( \tfrac{1}{2}, \tfrac{1}{12} \ep \right).$
\item
Whenever $e_0, e_1, \ldots, e_{n - 1}$
and $f_0, f_1, \ldots, f_{n - 1}$ are two sets of
\mops\  in $B$ such that $\| e_j f_j - e_j \| < 3 \ep_1$
for $0 \leq j \leq n - 1,$
then there exist \mops\  $h_0, h_1, \ldots, h_{n - 1} \in B$
such that $h_j \geq e_j$ and $\| h_j - f_j \| < \rh$
for $0 \leq j \leq n - 1.$
\end{itemize}

Choose $\ep_2 > 0$ so small that if $r_0, \ldots, r_{n - 1}$
are \pj s in a \ca\  $B$ such that $\| r_j r_k \| < 3 \ep_2$
for $j \neq k,$ then there is a \pj\  $z$ in the C*-subalgebra
of $B$ generated by $\sum_{j = 0}^{n - 1} r_j,$
and there are \mops\  $z_j \in B,$
such that $\sum_{j = 0}^{n - 1} z_j = z,$
and such that $\| z_j - r_j \| < \ep_1$ for $0 \leq j \leq n - 1.$
We also require $\ep_2 \leq \min \left( \tfrac{1}{6}, \ep_1 \right).$

Choose $\ep_3 > 0$ so small that whenever $B$ is a \ca\  and
$e, q \in B$ are \pj s such that $\| e q - q e \| < \ep_3,$
then there exists a \pj\  $r \leq q$ such that $\| r - e q \| < \ep_2.$
We also require $\ep_3 \leq \min \left( \ep_2, \, 1 / n^2 \right).$

Apply Lemma~\ref{PjAndNorm},
obtaining a nonzero \pj\  $q \in {\overline{x A x}}$
such that:
\newcounter{TmpEnumi}
\begin{enumerate}
\item\label{Pf:1} %
$\| r x r - r \| < \ep_1$
for every nonzero \pj\  $r \leq q.$
\setcounter{TmpEnumi}{\value{enumi}}
\end{enumerate}
Using Lemma~\ref{OrthInSP},
find nonzero \pj s $q_0, q_1, \ldots, q_{n - 1}, g \in A$
such that $g + \sum_{j = 0}^{n - 1} q_j = q.$
Lemma~\ref{SPAndInv} provides a nonzero \pj\  $g_0 \leq g$
such that $\af^j (g_0) \precsim q_j$
for $0 \leq j \leq n - 1.$

Set $F_1 = \bigcup_{j = 0}^{n - 1} \af^j (F \cup \{ g_0 \}).$
By the \tRp\  (Definition~\ref{NewTRPDfn}),
there exist orthogonal \pj s $f_0, f_1, \ldots, f_{n - 1} \in A$
such that, with $f = \sum_{j = 0}^{n - 1} f_j,$ we have:
\begin{enumerate}
\setcounter{enumi}{\value{TmpEnumi}}
\item\label{Pf:2} %
With indices taken mod $n,$ we have
$\| \af^m (f_j) - f_{j + m} \| < \ep_3$
for $0 \leq j, m \leq n - 1.$
\item\label{Pf:3} %
$\| f_j a - a f_j \| < \ep_3$ for
$0 \leq j \leq n - 1$ and all $a \in F_1.$
\item\label{Pf:4} %
$1 - f$ is \mvnt\  to a \pj\  in $g_0 A g_0.$
\item\label{Pf:5} %
$\| f g_0 f \| > 1 - \ep_3.$
\setcounter{TmpEnumi}{\value{enumi}}
\end{enumerate}

For $j \neq k$ we have $\| f_j g_0 f_k \| < \ep_3,$
so
\begin{align*}
\max_{0 \leq j \leq n - 1} \| f_j g_0 f_j \|
 & = \left\| \Ssum{j = 0}{n - 1} f_j g_0 f_j \right\|        \\
 & > \| f g_0 f \| - n (n - 1) \ep_3
   > 1 - \ep_3 - n (n - 1) \ep_3
   > 1 - n^2 \ep_3.
\end{align*}
Therefore there is $j$ such that $\| f_j g_0 f_j \| > 1 - n^2 \ep_3.$
Since all the other conditions on the $f_j$ are invariant under
cyclic permutation of the indices,
\wolog\  $\| f_0 g_0 f_0 \| > 1 - n^2 \ep_3.$

By the choice of $\ep_3,$ there are \pj s $r_1 \leq g_0$
and $r_2 \leq f_0$ such that $\| r_1 - f_0 g_0 \| < \ep_2$
and $\| r_2 - f_0 g_0 \| < \ep_2.$
In particular:
\begin{enumerate}
\setcounter{enumi}{\value{TmpEnumi}}
\item\label{Pf:6} %
$\| r_1 - r_2 \| < 2 \ep_2.$
\setcounter{TmpEnumi}{\value{enumi}}
\end{enumerate}
Further, using $\| f_0 g_0 - g_0 f_0 \| < \ep_3 \leq \ep_2,$
we get $\| r_1 - f_0 g_0 f_0 \| < 2 \ep_2.$
Since $\| f_0 g_0 f_0 \| > 1 - n^2 \ep_3$
and $3 \ep_2 + n^2 \ep_3 < 1,$
it follows that $r_1 \neq 0.$

Set
\[
F_2 = F_1 \cup \{ r_1 \}
        \cup \{ \af^j (f_k) \colon 1 \leq j, k \leq n - 1 \}.
\]
Apply strong tracial approximate innerness (Definition~\ref{STAInnDfn}),
obtaining a \pj\  $p_0 \in A$ and a partial isometry $w \in A$
such that:
\begin{enumerate}
\setcounter{enumi}{\value{TmpEnumi}}
\item\label{Pf:7} %
$w^* w = p$ and $w w^* = \af (p).$
\item\label{Pf:8} %
$\| p a - a p \| < \ep_3$ for all $a \in F_2.$
\item\label{Pf:9} %
$\| w p a p w^* - \af (p a p) \| < \ep_3$ for all $a \in F_2.$
\item\label{Pf:10} %
$1 - p \precsim r_1.$
\item\label{Pf:11} %
$\| p r_1 p \| > 1 - \ep_3.$
\setcounter{TmpEnumi}{\value{enumi}}
\end{enumerate}

By the choice of $\ep_3,$ there are \pj s $s_1 \leq f_0$
and $s_2 \leq p$
such that:
\begin{enumerate}
\setcounter{enumi}{\value{TmpEnumi}}
\item\label{Pf:12} %
$\| s_1 - f_0 p \| < \ep_2$ and $\| s_2 - f_0 p \| < \ep_2.$
\setcounter{TmpEnumi}{\value{enumi}}
\end{enumerate}
Then also:
\begin{enumerate}
\setcounter{enumi}{\value{TmpEnumi}}
\item\label{Pf:13} %
$\| s_1 - f_0 p f_0 \| < \ep_2.$
\item\label{Pf:14} %
$\| s_2 - p f_0 p \| < \ep_2.$
\item\label{Pf:15} %
$\| s_1 - s_2 \| < 2 \ep_2.$
\setcounter{TmpEnumi}{\value{enumi}}
\end{enumerate}
For $0 \leq j < k \leq n - 1,$
we have, using~(\ref{Pf:13}) at the third step,
$f_0 f_{k - j} = 0$ at the fourth step,
and~(\ref{Pf:2}) and $\ep_3 \leq \ep_2$ at the last step,
\begin{align*}
\| \af^j (s_1) \af^k (s_1) \|
 & = \| s_1 \af^{k - j} (s_1) \|                     \\
 & \leq \| f_0 p f_0 \af^{k - j} (f_0)
              \af^{k - j} (p) \af^{k - j} (f_0) \|
          + 2 \| s_1 - f_0 p f_0 \|                     \\
 & < \| f_0 \af^{k - j} (f_0) \| + 2 \ep_2
   \leq \| f_0 - \af^{k - j} (f_0) \| + 2 \ep_2
   < 3 \ep_2.
\end{align*}
By the choice of $\ep_2,$
there are \pj s $e_0, e_1, \ldots, e_{n - 1} \in A$
such that the sum $e = \sum_{j = 0}^{n - 1} e_j$
is in the \ca\  generated by $\sum_{j = 0}^{n - 1} \af^j (s_1)$
(and, in particular, is $\af$-invariant),
and such that:
\begin{enumerate}
\setcounter{enumi}{\value{TmpEnumi}}
\item\label{Pf:16} %
$\| e_j - \af^j (s_1) \| < \ep_1$ for $0 \leq j \leq n - 1.$
\setcounter{TmpEnumi}{\value{enumi}}
\end{enumerate}
We then get:
\begin{enumerate}
\setcounter{enumi}{\value{TmpEnumi}}
\item\label{Pf:17} %
$\| \af^m (e_j) - e_{j + m} \| < 2 \ep_1,$
with indices taken mod $n,$
for $0 \leq j, k \leq n - 1.$
\setcounter{TmpEnumi}{\value{enumi}}
\end{enumerate}

Define
$c_j = e_{j + 1} \af^j (w) e_j$ for $0 \leq j \leq n - 2.$
(We can't use $e_{j + 1} w e_j,$
imitating the proof of Step~1 of Proposition~\ref{AppInnRP},
because we don't know $e_j$ and $e_{j + 1}$ are even
approximately dominated by $p$ and $\af (p).$)
We claim that
\[
\| c_j^* c_j - e_j \| < 27 \ep_1
\andeqn
\| c_j c_j^* - e_{j + 1} \| < 27 \ep_1
\]
for $0 \leq j \leq n - 2.$
We start by observing
that~(\ref{Pf:16}), (\ref{Pf:15}), and~(\ref{Pf:14}) imply
\[
\| e_0 - p f_0 p \|
 \leq \| e_0 - s_1 \| + \| s_1 - s_2 \| + \| s_2 - p f_0 p \|
 < \ep_1 + 2 \ep_2 + \ep_2
 \leq 4 \ep_1.
\]
So~(\ref{Pf:9}) implies
\[
\| w e_0 w^* - \af (e_0) \|
 < 8 \ep_1 + \| w p f_0 p w^* - \af (p f_0 p) \|
 < 8 \ep_1 + \ep_3
 \leq 9 \ep_1.
\]
Combining this with $\| \af (e_0) - e_1 \| < 2 \ep_1$
(from~(\ref{Pf:17})),
we get:
\begin{enumerate}
\setcounter{enumi}{\value{TmpEnumi}}
\item\label{Pf:18} %
$\| w e_0 w^* - e_1 \| < 11 \ep_1.$
\setcounter{TmpEnumi}{\value{enumi}}
\end{enumerate}
Therefore also $\| e_1 w^* e_0 w e_1 - e_1 \| < 11 \ep_1.$
Furthermore, using $w^* w = p$ at the first step,
$s_2 \leq p$ at the second step,
and~(\ref{Pf:15}) and~(\ref{Pf:16}) at the third step,
\begin{align*}
\| w^* e_1 w - e_0 \|
 & \leq \| w^* \| \cdot \| e_1 - w e_0 w^* \| \cdot \| w \|
              + \| p e_0 p - e_0 \|                     \\
 & \leq \| e_1 - w e_0 w^* \|
              + 2 ( \| s_2 - s_1 \| + \| s_1 - e_0 \|)
                                   \\
 & < 11 \ep_1 + 4 \ep_2 + 2 \ep_1
   \leq 17 \ep_1.
\end{align*}
So $\| e_0 w^* e_1 w e_0 - e_0 \| < 17 \ep_1.$
This does the case $j = 0$ of the claim,
with $17 \ep_1$ in place of $27 \ep_1.$
For the general case, use~(\ref{Pf:17}) and the definition of $c_j$
to get:
\begin{enumerate}
\setcounter{enumi}{\value{TmpEnumi}}
\item\label{Pf:19} %
$\| c_j - \af^j (c_0) \| < 4 \ep_1.$
\setcounter{TmpEnumi}{\value{enumi}}
\end{enumerate}
Use this and~(\ref{Pf:17}) again to get
\[
\| c_j^* c_j - e_j \|
 < 10 \ep_1 + \| c_0^* c_0 - e_0 \|
 < 27 \ep_1
\]
and,
similarly, $\| c_j c_j^* - e_{j + 1} \| < 27 \ep_1.$
This proves the claim.

By the choice of $\ep_1,$
there are partial isometries
$v_j \in B$ such that:
\begin{enumerate}
\setcounter{enumi}{\value{TmpEnumi}}
\item\label{Pf:20} %
$\| v_j - c_j \| < \ep_0,$
$v_j^* v_j = e_j,$ and $v_j v_j^* = e_{j + 1}$
for $0 \leq j \leq n - 2.$
\setcounter{TmpEnumi}{\value{enumi}}
\end{enumerate}
Define
$v = v_0^* v_1^* \cdots v_{n - 2}^* + \sum_{j = 0}^{n - 2} v_j.$
Then it is easily checked that $v$ is a unitary in $e A e$
which satisfies $v^n = 1.$

We now have $e$ and $v,$
and Part~(\ref{FinOrdTAInn:1}) of the conclusion has been verified.

For Part~(\ref{FinOrdTAInn:2}), if $a \in F$ and $0 \leq j \leq n - 1,$
then $\af^{-j} (a) \in F_1 \subset F_2,$
so, using~(\ref{Pf:17}),
(\ref{Pf:16}), (\ref{Pf:12}), (\ref{Pf:8}), and~(\ref{Pf:3}),
\begin{align*}
\| e_j a - a e_j \|
 & \leq 2 \| e_j - \af^j (e_0) \|
      + \| e_0 \af^{-j} (a) - \af^{-j} (a) e_0 \|       \\
 & \leq 2 \| e_j - \af^j (e_0) \| + 2 \| e_0 - s_1 \|          \\
 & \hspace*{3em} \mbox{}
           + 2 \| s_1 - f_0 p \|
           + \| p a - a p \| + \| f_0 a - a f_0 \|       \\
 & < 4 \ep_1 + 2 \ep_1 + 2 \ep_2 + \ep_3 + \ep_3
   \leq 10 \ep_1.
\end{align*}
Since $e = \sum_{j = 0}^{n - 1} e_j,$
it follows that
$\| e a - a e \| < 10 n \ep_1 \leq \ep.$
This is Part~(\ref{FinOrdTAInn:2}) of the conclusion.
The relation $\| e_j a - a e_j \| < 10 \ep_1$ furthermore implies that:
\begin{enumerate}
\setcounter{enumi}{\value{TmpEnumi}}
\item\label{Pf:21} %
$\left\| e a e - \sum_{j = 0}^{n - 1} e_j a e_j \right\|
  < 10 n (n - 1) \ep_1.$
\setcounter{TmpEnumi}{\value{enumi}}
\end{enumerate}

To start~(\ref{FinOrdTAInn:3}),
we estimate $\| v_j e_j a e_j v_j^* - \af (e_j a e_j) \|$
for $a \in F$ and $0 \leq j \leq n - 2.$
We begin with $j = 0$
and $a \in \bigcup_{j = 0}^{n - 1} \af^{-j} (F) \S F_2.$
Use~(\ref{Pf:16}), (\ref{Pf:15}), and $s_2 \leq p$ to get
\[
\| e_0 p - e_0 \|
 \leq 2 \| e_0 - s_1 \| + 2 \| s_1 - s_2 \|
 < 2 \ep_1 + 4 \ep_2
 \leq 6 \ep_1,
\]
and then use $w^* w = p$ and~(\ref{Pf:18}) to get
\[
\| e_1 w - w e_0 \|
 \leq \| e_1 - w e_0 w^* \| \cdot \| w \|
         + \| e_0 p - e_0 \|
 < 11 \ep_1 + 6 \ep_1
 = 17 \ep_1.
\]
Next, use the definition of $c_0$ in the first step
and~(\ref{Pf:17}) and~(\ref{Pf:9}) in the last step to get
\begin{align*}
\| c_0 e_0 a e_0 c_0^* - \af (e_0 a e_0) \|
 & \leq 4 \| e_0 - e_0 p \|
    + \| e_1 w e_0 p a p e_0 w^* e_1 - \af (e_0 p a p e_0) \|     \\
 & \leq 4 \| e_0 - e_0 p \|
          + 2 \| e_1 w e_0 - w e_0 \|          \\
 & \hspace*{3em} \mbox{}
          + 2 \| \af (e_0) - e_1 \|
          + \| e_1 [w p a p w^* - \af (p a p)] e_1 \|     \\
 & < 24 \ep_1 + 34 \ep_1 + 4 \ep_1 + \ep_3
   \leq 63 \ep_1.
\end{align*}
Using~(\ref{Pf:20}) and $\ep_1 \leq \ep_0,$ we get
$\| v_0 e_0 a e_0 v_0^* - \af (e_0 a e_0) \| < 65 \ep_0.$

Now, for $a \in F$ and $0 \leq j \leq n - 2,$
combine
\[
\| c_0 e_0 \af^{-j} (a) e_0 c_0^* - \af (e_0 \af^{-j} (a) e_0) \|
        < 63 \ep_1
\]
with~(\ref{Pf:17}) and~(\ref{Pf:19}) to get
\[
\| c_j e_j a e_j c_j^* - \af (e_j a e_j) \|
 < 63 \ep_1 + 2 \| c_j - \af^j (c_0) \|
          + 4 \| e_j - \af^j (e_0) \|
 < 79 \ep_1.
\]
Using~(\ref{Pf:20}) and $\ep_1 \leq \ep_0,$ we then get
\begin{enumerate}
\setcounter{enumi}{\value{TmpEnumi}}
\item\label{Pf:22} %
$\| v_j e_j a e_j v_j^* - \af (e_j a e_j) \| < 81 \ep_0$
for $0 \leq j \leq n - 2$
and $a \in \bigcup_{j = 0}^{n - 1} \af^j (F).$
\setcounter{TmpEnumi}{\value{enumi}}
\end{enumerate}

Set
\[
y_k = \af^{k - 1} (v_0) \af^{k - 2} (v_0) \cdots \af (v_0) v_0.
\]
An induction argument gives
$\| \af^k (e_0 a e_0) - y_k e_0 a e_0 y_k^* \| < 81 k \ep_0$
for $k \geq 0$ and $a \in \bigcup_{j = 0}^{n - 1} \af^j (F).$
The inequalities~(\ref{Pf:19}) and~(\ref{Pf:20})
give $\| \af^j (v_0) - v_j \| < 4 \ep_1 + 2 \ep_0 \leq 6 \ep_0.$
Therefore, with $z = v_0^* v_1^* \cdots v_{n - 2}^*,$ we get
\[
\| \af^{n - 1} (e_0 a e_0) - z^* e_0 a e_0 z \|
  < 81 (n - 1) \ep_0 + 6 (n - 1) \ep_0
  = 87 (n - 1) \ep_0.
\]
Putting $\af (a)$ in place of $a,$ using $\af^n = \id_A,$
and rearranging,
we get
\[
\| z \af^{n - 1} (e_0) a \af^{n - 1} (e_0) z^*
    - \af (\af^{n - 1} (e_0) a \af^{n - 1} (e_0)) \|
  < 87 (n - 1) \ep_0.
\]
{}From~(\ref{Pf:17}) and $\ep_1 \leq \ep_0,$ it follows that
\[
\| z e_{n - 1} a e_{n - 1} z^* - \af (e_{n - 1} a e_{n - 1}) \|
  < [87 (n - 1) + 8] \ep_0.
\]

Using the definition of $v$ and~(\ref{Pf:22}),
we now get
\[
\left\| v \left( \Ssum{j = 0}{n - 1} e_j a e_j \right) v^*
         - \af \left( \Ssum{j = 0}{n - 1} e_j a e_j \right) \right\|
   < 81 (n - 1) \ep_0 + [87 (n - 1) + 8] \ep_0
   \leq 168 n \ep_0.
\]
{}From~(\ref{Pf:21}) it now follows that
\[
\| v e a e v^* - \af (e a e) \|
  < 168 n \ep_0 + 20 n (n - 1) \ep_1
  \leq 168 n^2 \ep_0
  \leq \ep.
\]
This completes the proof of Part~(\ref{FinOrdTAInn:3})
of the conclusion.

For Part~(\ref{FinOrdTAInn:4}),
begin by using~(\ref{Pf:16}), (\ref{Pf:2}), and $s_1 f_0 = s_1$ to
estimate
\[
\| e_j f_j - e_j \|
    \leq 2 \| e_j - \af^j (s_1) \| + \| f_j - \af^j (f_0) \|
    < 2 \ep_1 + \ep_3
    \leq 3 \ep_1.
\]
The conditions on $\ep_1$ provide \mops\  %
$h_j \geq e_j$ such that $\| h_j - f_j \| < \rh.$
Set $h = \sum_{j = 0}^{n - 1} h_j.$
Then
$1 - e = 1 - h + \sum_{j = 0}^{n - 1} (h_j - e_j).$
Since $\| h - f \| < n \rh \leq 1,$
it follows (using~(\ref{Pf:4}) at the second step) that
\[
1 - h \sim 1 - f \precsim g_0 \leq g.
\]

The next step is to show that $h_0 - e_0 \precsim 1 - p.$
We have, using~(\ref{Pf:12}), (\ref{Pf:15}), (\ref{Pf:16}), and $s_2 \leq p$
at the third step,
\begin{align*}
\| (1 - p) (h_0 - e_0) - (h_0 - e_0) \|
  &  = \| p h_0 - p e_0 \|                    \\
  &  \leq \| h_0 - f_0 \| + \| p f_0 - p s_2 \| + \| s_2 - s_1 \|
                 + \| s_1 - e_0 \|                    \\
  & < \rh + \ep_2 + 2 \ep_2 + \ep_1
    < 1.
\end{align*}
So $h_0 - e_0 \precsim 1 - p$ by Lemma~2.5.2 of~\cite{LnBook}.

Next, using~(\ref{Pf:2}), for $0 \leq j \leq n - 1$ we have
\[
\| h_j - \af^j (h_0) \|
  \leq \| h_j - f_j \| + \| h_0 - f_0 \| + \| f_j - \af^j (f_0) \|
  < \rh + \rh + \ep_3,
\]
and so, using~(\ref{Pf:17}) and $\ep_3 \leq \ep_1,$
\[
\| (h_j - e_j) - \af^j (h_0 - e_0) \|
  < 2 \rh + 3 \ep_1
  \leq 1.
\]
Thus, using~(\ref{Pf:10}) at the third step,
\[
h_j - e_j \sim \af^j (h_0 - e_0)
   \precsim \af^j (1 - p)
   \precsim \af^j (r_1)
   \leq \af^j (g_0)
   \precsim q_j.
\]
We conclude that
\[
1 - e = 1 - h + \sum_{j = 0}^{n - 1} (h_j - e_j)
     \precsim g + \sum_{j = 0}^{n - 1} q_j
     \leq q
     \in {\overline{x A x}}.
\]
This proves Part~(\ref{FinOrdTAInn:4}).

It remains to prove Part~(\ref{FinOrdTAInn:5}) of the conclusion.
Combine~(\ref{Pf:12}) and~(\ref{Pf:16}) to get
$\| e_0 - f_0 p \| < \ep_2 + \ep_1 \leq 2 \ep_1.$
Use~(\ref{Pf:6}) at the second and third steps, $r_2 \leq f_0$ at
the second step, and~(\ref{Pf:11}) at the fourth step, to get
\begin{align*}
\| e_0 r_1 e_0 \|
 & > \| p f_0 r_1 f_0 p \| - 4 \ep_1
   > \| p r_2 p \| - 2 \ep_2 - 4 \ep_1
   > \| p r_1 p \| - 4 \ep_2 - 4 \ep_1             \\
 & > 1 - \ep_3 - 4 \ep_2 - 4 \ep_1
   \geq 1 - 9 \ep_1.
\end{align*}
Since $1 - 9 \ep_1 \geq \tfrac{1}{4},$
we get in particular $\| e_0 r_1 \| > \tfrac{1}{2}.$

Next, we estimate $\| e_0 r_1 - r_1 e_0 \|.$
We have $\| p r_1 - r_1 p \| < \ep_3$ by~(\ref{Pf:8}),
and $\| f_0 r_1 - r_1 f_0 \| \leq 2 \| r_1 - r_2 \| < 4 \ep_2$
by~(\ref{Pf:6}) and because $r_2 \leq f_0.$
So
\[
\| e_0 r_1 - r_1 e_0 \|
   \leq 2 \| e_0 - f_0 p \| + \| f_0 r_1 - r_1 f_0 \|
                + \| p r_1 - r_1 p \|
   < 4 \ep_1 + 4 \ep_2 + \ep_3
   \leq 9 \ep_1.
\]
It follows from the choice of $\ep_1$ that there is a \pj\  $r \leq e_0$
such that
$\| r - e_0 r_1 \|
 < \min \left( \tfrac{1}{2}, \tfrac{1}{12} \ep \right).$
Since $\| e_0 r_1 \| > \tfrac{1}{2},$ we have $r \neq 0.$

Now use $r \leq e_0 \leq e$ at the first step,
and $r_1 \leq g_0 \leq q$ and~(\ref{Pf:1}) at the third step,
to get
\[
\| e x e \|
  \geq \| r x r \|
  > \| e_0 r_1 x r_1 e_0 \| - \tfrac{1}{6} \ep
  > \| e_0 r_1 e_0 \| - \ep_1 - \tfrac{1}{6} \ep
  > 1 - 9 \ep_1 - \ep_1 - \tfrac{1}{6} \ep
  \geq 1 - \ep.
\]
This completes the proof.
\end{proof}

\begin{thm}\label{FinOrdTAInn2}
Let $A$ be an \idssuca\  and let $\af \in \Aut (A)$
be strongly tracially approximately inner and satisfy $\af^n = \id_A.$
Suppose $\af$ generates an action of $\Zqn$ with the
tracial Rokhlin property.
Then this action is tracially approximately representable
in the sense of Definition~\ref{TrAppRepDfn}.
\end{thm}

\begin{proof}
If $A$ does not have Property~(SP),
then $\af$ is approximately inner (by Remark~\ref{InnAndTAInn})
and the action of
$\Zqn$ it generates has the \sRp\  (by Lemma~\ref{RImpSP}),
so the result is Proposition~\ref{AppInnRP}.
Accordingly, we assume that $A$ has Property~(SP).
Let $F \S A$ be finite, let $\ep > 0,$ and let $x \in A$
be a positive element with $\| x \| = 1.$

Choose $\ep_1 > 0$ with
\[
\ep_1 \leq \min \left(
       \frac{1}{n}, \, \frac{\ep}{3 n (2 n^2 + 2 n + 3)},
        \, \frac{\ep}{2 (n + 5)} \right),
\]
and also
(using semiprojectivity of $\C^n,$
as at the beginning of the proof of Proposition~\ref{AppInnRP})
so small that whenever $B$ is a unital \ca\  and $c \in B$
satisfies
\[
\| c^* c - 1 \| < 2 (n - 1) \ep_1, \,\,\,\,\,\,
\| c c^* - 1 \| < 2 (n - 1) \ep_1, \andeqn
\| c^n - 1 \| < 2 (n - 1) \ep_1,
\]
then there is a unitary $v \in B$ such that $v^n = 1$
and $\| v - c \| < \ep / (3 n).$
Choose $\ep_2 > 0$ with $\ep_2 \leq \frac{1}{8} \ep_1,$
and also so small that whenever $B$ is a unital \ca,
$C \S B$ is a subalgebra, $p \in B$ is a \pj,
and $b \in C$ satisfies $\| b - p \| < 3 n^2 \ep_2,$
then there is a \pj\  $q \in C$ and a unitary $z \in B$
such that $z p z^* = q$ and $\| z - 1 \| < \frac{1}{8} \ep_1.$
Choose $\ep_3 > 0$ with $\ep_3 \leq \min \left( \frac{1}{2}, \ep_2 \right),$
and also so small that whenever $B$ is a unital \ca,
$f_0, f_1, \ldots, f_{n - 1} \in B$ are \mops,
and $e \in B$ is a \pj\  such that $\| [ f_j, e ] \| < \ep_3$
for $0 \leq j \leq n - 1,$
then there exist \mops\  $g_0, g_1, \ldots, g_{n - 1} \leq e$
such that $\| g_j - f_j e \| < \ep_2$ for $0 \leq j \leq n - 1.$

By considering polynomial approximations
to the continuous functional calculus,
choose $\dt > 0$ with $\dt \leq \tfrac{1}{6} \ep,$
and also so small that whenever $B$ is a unital \ca,
$x \in B$ satisfies $0 \leq x \leq 1,$
and $q \in B$ is a \pj\  such that $\| [x, q] \| < \dt,$
then $\| [x^{1/2}, \, q] \| < \tfrac{1}{6} \ep.$
Apply Lemma~\ref{PjAndNorm}
to find a nonzero \pj\  $q_0 \in {\overline{x A x}}$
such that whenever $q \leq q_0$ is a nonzero \pj,
then $\| q x q - q \| < \dt$ and $\| [x, q] \| < \dt.$
Use Lemma~\ref{OrthInSP} to choose orthogonal nonzero \pj s
$q_1, q_2 \leq q.$

Apply Lemma~\ref{FinOrdTAInn}
with $\bigcup_{k = 0}^{n - 1} \af^k (F)$ in place of $F,$
with $\ep_3$ in place of $\ep,$
and with $q_1$ in place of $x.$
Call the resulting \pj\  $e,$
and let $u \in e A e$ be the resulting unitary.

Set $\ld = \| e q_1 e \| > 1 - \ep_3,$
and set $b = \ld^{-1} e q_1 e.$
Thus $\| b - e q_1 e \| = \ld^{-1} - 1 < 2 \ep_3.$
Apply Lemma~\ref{PjAndNorm}
to find a nonzero \pj\  $r \in {\overline{e q_1 e A e q_1 e}}$
such that whenever $p \leq r$ is a nonzero \pj,
then $\| p b p \| > 1 - \ep_3,$ $\| p b p - p \| < \ep_3,$
and $\| p b - b p \| < \ep_3.$
By Lemma~\ref{L:CompSP}, we may require $r \precsim q_2.$
It follows that, for such $p,$ we have
$\| p e q_1 e p \| > 1 - 3 \ep_3$
and $\| p e q_1 e p - p \| < \ep_3.$

Apply the \tRp\  with
$\{ e \} \cup \bigcup_{k = 0}^{n - 1} \af^k (F \cup \{ u \})$
in place of $F,$
with $\ep_3$ in place of $\ep,$
and with $r$ in place of $x.$
Let $f_0, f_1, \ldots, f_{n - 1}$ be the resulting \pj s,
and set $f = \sum_{k = 0}^{n - 1} f_k.$

Since $\| [ f_k, e ] \| < \ep_3$
for $0 \leq k \leq n - 1,$
there exist \mops\  $g_0, g_1, \ldots, g_{n - 1} \leq e$
such that $\| g_k - f_k e \| < \ep_2$ for $0 \leq k \leq n - 1.$
Set $g = \sum_{k = 0}^{n - 1} g_k.$
We have
\[
\| \af (g_k) - g_{k + 1} \|
 < 2 \ep_2 + \| \af (f_k) - f_{k + 1} \|
 < 3 \ep_2,
\]
so that $\| \af (g) - g \| < 3 n \ep_2$
and $\| \af^k (g) - g \| < 3 k n \ep_2 \leq 3 n^2 \ep_2.$
With $b = \frac{1}{n} \sum_{k = 0}^{n - 1} \af^k (g),$
we thus get $\| b - g \| < 3 n^2 \ep_2,$
so the choice of $\ep_2$ provides an $\af$-invariant \pj\  %
$h \in e A e$ and a unitary $z \in A$ such that $z g z^* = h$
and $\| z - 1 \| < \frac{1}{8} \ep_1.$
Set $h_k = z g_k z^*$ for $0 \leq k \leq n - 1.$
Then $\| h_k - g_k \| < \frac{2}{8} \ep_1,$
so
\[
\| h_k - f_k e \| < \tfrac{2}{8} \ep_1 + \ep_2 \leq \tfrac{3}{8} \ep_1.
\]
Since $\af (e) = e,$ we get
\[
\| \af (h_k) - h_{k + 1} \|
  \leq \| h_k - f_k e \| + \| h_{k + 1} - f_{k + 1} e \|
        + \| \af (f_k) - f_{k + 1} \|
  < 2 \left( \tfrac{3}{8} \ep_1 \right) + \ep_3
  < \ep_1
\]
for $0 \leq k \leq n - 1,$
\[
\| [h_j, \af^k (u) ] \|
 \leq 2 \| h_k - f_k e \| + \| [f_j, \af^k (u) ] \|
  < 2 \left( \tfrac{3}{8} \ep_1 \right) + \ep_3
  < \ep_1
\]
for $0 \leq k \leq n - 1,$ and
\[
\| [h_j, \af^k (a) ] \|
 \leq 2 \| h_k - f_k e \| + \| [f_j, \af^k (a) ] \|
            + \| [e, \af^k (a) ] \|
  < 2 \left( \tfrac{3}{8} \ep_1 \right) + 2 \ep_3
  < \ep_1
\]
for $a \in F$ and $0 \leq k \leq n - 1.$

We now follow Step~2 of the proof of Proposition~\ref{AppInnRP}.
Set $w_k = h_k \af^k (u) h_k$ for $0 \leq k \leq n - 1,$
and set $w = \sum_{k = 0}^{n - 1} w_k.$
The same estimates as there give
\[
\| w w^* - h \| < \ep_1, \,\,\,\,\,\,
\| w^* w - h \| < \ep_1, \andeqn
\| w^n - h \| < (n - 1) \ep_1.
\]
As there, $y = \sum_{k = 0}^{n - 1} \af^k (h_0 u h_0) \in h A^{\af} h$
satisfies $\| y - w \| < 2 n \ep_1,$
so there is a unitary $v \in h A^{\af} h$ such that
$v^n = h$ and
$\| v - w \| < 2 n \dt + \ep / (3 n).$
Moreover, also using the same reasoning as there, we have
$\| v h a h v^* - \af (h a h) \| < \ep / n$ for $a \in F.$
Setting $v_k = v^k$ for $0 \leq k \leq n - 1,$
and iterating the last estimate as required,
we now have
Conditions~(\ref{TrAppRepDfn:2}), (\ref{TrAppRepDfn:3}),
and~(\ref{TrAppRepDfn:4}) of Definition~\ref{TrAppRepDfn}.

For Condition~(\ref{TrAppRepDfn:1}),
we estimate
\[
\| [h, a] \| \leq \sum_{k = 0}^{n - 1} \| h_k, a ] \|
 < n \ep_1
 \leq \ep.
\]

For Condition~(\ref{TrAppRepDfn:5}),
we have $1 - e \precsim q_1$ by the choice of $e.$
Also,
\[
\| (1 - f) (e - h) - (e - h) \|
  = \| f (e - h) \|
  \leq \| f e - h \|
  \leq \sum_{k = 1}^{n - 1} \| h_k - f_k e \|
  < \tfrac{3}{8} n \ep_1
  < 1.
\]
It follows from Lemma~2.5.2 of~\cite{LnBook} and the choice of $f$
that $e - h \precsim 1 - f \precsim r \precsim q_2.$
So $1 - h \precsim q_1 + q_2 \in {\overline{x A x}},$
as desired.

It remains only to prove Condition~(\ref{TrAppRepDfn:6}).
We observed in the previous paragraph that
$\| f e - h \| < \tfrac{3}{8} n \ep_1 < \tfrac{1}{2} n \ep_1.$
The choice of $r$ implies $\| r e q_1 e - r \| < 2 \ep_3$
and $\| e q_1 e r - r \| < 2 \ep_3.$
Therefore
\begin{align*}
\| h q_1 h \|
 & > \| f e q_1 e f \| - n \ep_1
   = \| q_1 e f e q_1 \| - n \ep_1
   \geq \| r q_1 e f e q_1 r \| - n \ep_1
             \\
 &
   > \| r f r \| - 4 \ep_3 - n \ep_1
   = \|f  r f \| - 4 \ep_3 - n \ep_1
   > 1 - 5 \ep_3 - n \ep_1
   \geq 1 - \tfrac{1}{2} \ep.
\end{align*}
So
\begin{align*}
\| h x h \|
 & = \| x^{1/2} h x^{1/2} \|
   \geq \| q_1 x^{1/2} h x^{1/2} q_1 \|
   = \| h x^{1/2} q_1 x^{1/2} h \|
             \\
 &
   \geq \| h q_1 x q_1 h \| - 2 \| [ x^{1/2}, q_1] \|
   \geq \| h q_1 h \| - \| q_1 x q_1 - q_1 \| - 2 \| [ x^{1/2}, q_1] \|
             \\
 &
   > \left( 1 - \tfrac{1}{2} \ep \right)
      - \dt - 2 \left( \tfrac{1}{6} \ep \right)
   \geq 1 - \ep,
\end{align*}
as desired.
\end{proof}

We finish this section by giving one way in which
the three main results so far can be combined.
This theorem is what is really needed for the
classification
of higher dimensional noncommutative toruses~\cite{PhtRp2}.
We first collect some information from~\cite{Rs}.

\begin{prp}\label{FPIsSimple}
Let $A$ be a \ca, let $G$ be a compact group, and let
$\af \colon G \to \Aut (A)$ be a \ct\  action of $G$ on $A.$
Suppose $C^* (G, A, \af)$ is simple.
Then the fixed point algebra $A^{\af}$ is simple,
is isomorphic to a full \hsa\  of $C^* (G, A, \af),$ and is
strongly Morita equivalent to $C^* (G, A, \af).$
\end{prp}

\begin{proof}
See the Proposition, Corollary, and proof of the Corollary in~\cite{Rs}.
\end{proof}

\begin{thm}\label{P:FixPtAlg}
Let $A$ be an \idssuca, and let $\af \in \Aut (A)$ be a
strongly tracially approximately inner automorphism which
satisfies $\af^n = \id_A.$
Suppose $\af$ generates an action of $\Zqn$ with the
tracial Rokhlin property.
Then $A$ has tracial rank zero \ifo\  %
the fixed point algebra $A^{\af}$ has tracial rank zero.
\end{thm}

\begin{proof}
Assume that $A$ has tracial rank zero.
Theorem~\ref{RokhTAF} implies that $\CZnAa$
has tracial rank zero.
Corollary~\ref{CrPrIsSimple} implies that $C^* (\Zqn, A, \af)$
is simple.
It therefore follows from Proposition~\ref{FPIsSimple} that
$A^{\af}$ is isomorphic to a \hsa\  in $\CZnAa.$
So Theorem~3.12(1) of~\cite{LnTAF} implies that
$A^{\af}$ has tracial rank zero.

Now assume that $A^{\af}$ has tracial rank zero.
By the reverse of the reasoning in the first part of the proof,
$\CZnAa$ has tracial rank zero.
It follows from Theorem~\ref{FinOrdTAInn2}
that the action of $\Zqn$ on $A$
is tracially approximately representable,
and then from Theorem~\ref{TRPDualToTAppRep} that its dual action,
which by abuse of notation we call ${\widehat{\af}},$ has the \tRp.
So Theorem~\ref{RokhTAF} implies that
$C^* \big( {\widehat{\Zqn}}, \, \CZnAa, {\widehat{\af}} \big)
      \cong M_n \otimes A$
has tracial rank zero.
Then Theorem~3.12(1) of~\cite{LnTAF} implies that
$A$ has tracial rank zero.
\end{proof}

\section{Tracially approximately inner automorphisms}\label{Sec:TAI}

In this section, we introduce tracially approximately inner automorphisms,
and we prove that they form a group
(Theorem~\ref{T:TAIIsGp}).
(We do not know whether the composition of
two strongly tracially approximately inner automorphisms
is again strongly tracially approximately inner;
in fact, we suspect that this is false.)
The following definition is an improvement over,
and uses a weaker condition than,
the definition of tracial approximate innerness given in~\cite{PhW}.
Again, we do not know, and suspect that it is not true,
that the tracially approximately inner automorphisms
as defined in~\cite{PhW} form a group.

\begin{dfn}\label{D:TAI}
Let $A$ be an \idssuca\  and let $\af \in \Aut (A).$
We say that $\af$ is
{\emph{tracially approximately inner}} if for every finite set
$F \S A,$ every $\ep > 0,$
and every positive element $x \in A$ with $\| x \| = 1,$
there exist $p_1, p_2, v \in A$
such that:
\begin{enumerate}
\item\label{D:TAI:1} %
$p_1$ and $p_2$ are \pj s,
$v^* v = p_1,$ and $v v^* = \af (p_2).$
\item\label{D:TAI:2} %
$\| p_1 a - a p_1 \| < \ep$ and $\| p_2 a - a p_2 \| < \ep$
for all $a \in F.$
\item\label{D:TAI:3} %
$\| v p_1 a p_1 v^* - \af (p_2 a p_2) \| < \ep$ for all $a \in F.$
\item\label{D:TAI:4} %
$1 - p_1$ and $1 - p_2$ are \mvnt\  to \pj s in ${\overline{x A x}}.$
\item\label{D:TAI:5} %
$\| p_1 x p_1 \| > 1 - \ep$ and $\| p_2 x p_2 \| > 1 - \ep.$
\end{enumerate}
\end{dfn}

\begin{lem}\label{L:TAIAndSP}
Let $A$ be an \idssuca,
and let $\af \in \Aut (A).$
If $\af$ is approximately inner,
then $\af$ is tracially approximately inner.
If $\af$ is tracially approximately inner
and $A$ does not have Property~(SP),
then $\af$ is approximately inner.
\end{lem}

\begin{proof}
The first statement is obvious.
For the second,
applying the definition with an element $x$ such that
${\overline{x A x}}$ has no nontrivial \pj s,
one forces $p_1 = p_2 = 1,$
so that $v$ is unitary.
\end{proof}

\begin{lem}\label{L:ImpTAI}
Let $A$ be an \idssuca,
and let $\af \in \Aut (A).$
If $\af$ is strongly tracially approximately inner
(Definition~\ref{STAInnDfn}),
then $\af$ is tracially approximately inner.
\end{lem}

\begin{proof}
This is obvious.
\end{proof}

\begin{lem}\label{L:TAIForFin}
Let $A$ be an \idssuca,
and let $\af \in \Aut (A).$
Suppose that $A$ is finite.
Then $\af$ is tracially approximately inner \ifo\  %
$\af$ satisfies the conditions of Definition~\ref{D:TAI}
without Part~(\ref{D:TAI:5}).
\end{lem}

\begin{proof}
Let $\ep > 0,$ let $F \S A$ be finite,
and let $x \in A$ be positive with $\| x \| = 1.$
Apply Lemma~\ref{FPjNorm},
obtaining a \pj\  $q.$
Apply the hypotheses, with $\ep$ and $F$ as given,
and with $q$ in place of $x.$
Then the relations $1 - p_1 \precsim q$ and $1 - p_2 \precsim q$
imply $\| p_1 x p_1 \| > 1 - \ep$ and $\| p_2 x p_2 \| > 1 - \ep.$
\end{proof}

\begin{lem}\label{L:ManyElt}
Let $A$ be an \idssuca,
and let $\af \in \Aut (A).$
Then $\af$ is tracially approximately inner \ifo\  %
for every finite $F \S A,$ every $\ep > 0,$
and any three nonzero positive elements $x_0, x_1, x_2 \in A$
with $\| x_1 \| = \| x_2 \| = 1,$
there are $p_1, p_2, v \in A$ such that
Conditions~(\ref{D:TAI:1}), (\ref{D:TAI:2}), and~(\ref{D:TAI:3})
of Definition~\ref{D:TAI} hold,
and, in addition,
\begin{itemize}
\item[(4$'$)]
$1 - p_1$ and $1 - p_2$ are \mvnt\  to \pj s
in ${\overline{x_0 A x_0}}.$
\item[(5$'$)]
$\| p_1 x_1 p_1 \| > 1 - \ep$ and $\| p_2 x_2 p_2 \| > 1 - \ep.$
\end{itemize}
\end{lem}

\begin{proof}
It is trivial that the condition of the lemma implies
tracial approximate innerness.
For the reverse direction,
if $A$ does not have \PSP,
then $\af$ is approximately inner by Lemma~\ref{L:TAIAndSP},
and the condition of the lemma is immediate.
Accordingly, assume that $A$ has \PSP,
let $F \S A$ be finite, let $\ep > 0,$
and let $x_0, x_1, x_2 \in A$ be nonzero positive
elements with $\| x_1 \| = \| x_2 \| = 1.$

Set $\dt = \frac{1}{4} \ep.$
Use Lemma~\ref{PjAndNorm} to choose nonzero \pj s
$q_1 \in {\overline{x_1 A x_1}}$ and $q_2 \in {\overline{x_2 A x_2}}$
such that, whenever $e \leq q_j$ is a nonzero \pj,
then $\| e x_j e - e \| < \dt.$
Use Lemma~\ref{L:CompSP}
to find a nonzero \pj\  $e \in {\overline{x_0 A x_0}}$
and partial isometries $w_1, w_2 \in A$ such that
\[
w_1^* w_1 = w_2^* w_2 = e, \,\,\,\,\,\,
w_1 w_1^* \leq q_1, \andeqn w_2 w_2^* \leq q_2.
\]
Apply Definition~\ref{D:TAI}
with $F \cup \{ w_1, w_1^*, w_2, w_2^* \}$ in place of $F,$
with $\dt$ in place of $\ep,$
and with $e$ in place of $x,$
obtaining $p_1,$ $p_2,$ and $v.$
It is immediate that 
Conditions~(\ref{D:TAI:1}), (\ref{D:TAI:2}), and~(\ref{D:TAI:3})
of Definition~\ref{D:TAI}
and Condition~(4$'$) of the present lemma hold.
It remains to prove Condition~(5$'$).

Using
\[
\| (w_1 w_1^*) x_1 (w_1 w_1^*) - w_1 w_1^* \| < \dt, \,\,\,\,\,\,
w_1^* w_1 w_1^* = w_1^*, \andeqn
(w_1^* w_1)^2 = e,
\]
we get $\| w_1^* x_1 w_1 - e \| < \dt.$
Therefore
\begin{align*}
\| p_1 x_1 p_1 \|
 & \geq \| w_1^* p_1 x_1 p_1 w_1 \|
   \geq \| p_1 w_1^* x_1 w_1 p_1 \|
               - \| [p_1, w_1] \| - \| [p_1, w_1^*] \|
            \\
 &
   >  \| p_1 e p_1 \| - \dt 
               - \| [p_1, w_1] \| - \| [p_1, w_1^*] \|
   > 1 - 4 \dt = 1 - \ep.
\end{align*}
The same reasoning gives $\| p_2 x_2 p_2 \| > 1 - \ep.$
\end{proof}

For the proof of the next result,
it is convenient to observe that the relations in the definition
of tracial approximate innerness need only hold approximately.

\begin{lem}\label{L:ATAI}
Let $A$ be an \idssuca,
and let $\af \in \Aut (A).$
Suppose that for every $\ep > 0,$ every finite $F \S A,$
and any three nonzero positive elements $x_0, x_1, x_2 \in A$
with $\| x_1 \| = \| x_2 \| = 1,$
there are $e_1, e_2, y \in A$ such that:
\begin{enumerate}
\item\label{L:ATAI:0}
$\| e_1 \|, \| e_2 \|, \| y \| \leq 1.$
\item\label{L:ATAI:1} %
$\| e_1^2 - e_1 \|,$ $\| e_1^* - e_1 \|,$
$\| e_2^2 - e_2 \|,$ $\| e_2^* - e_2 \|,$
$\| y^* y - e_1 \|,$ and $\| y y^* - \af (e_2) \|$
are all less than~$\ep.$
\item\label{L:ATAI:2} %
$\| e_1 a - a e_1 \| < \ep$ and $\| e_2 a - a e_2 \| < \ep$
for all $a \in F.$
\item\label{L:ATAI:3} %
$\| y e_1 a e_1 y^* - \af (e_2 a e_2) \| < \ep$ for all $a \in F.$
\item\label{L:ATAI:4} %
$\frac{1}{2} \not\in \spec (e_1^* e_1) \cup \spec (e_2^* e_2),$
and
$1 - \ch_{[1/2, \, \I)} (e_1^* e_1)$
and $1 - \ch_{[1/2, \, \I)} (e_2^* e_2)$
are \mvnt\  to \pj s in ${\overline{x_0 A x_0}}.$
\item\label{L:ATAI:5} %
$\| e_1 x_1 e_1 \| > 1 - \ep$ and $\| e_2 x_2 e_2 \| > 1 - \ep.$
\end{enumerate}
Then $\af$ is tracially approximately inner.
\end{lem}

\begin{proof}
Let $F,$ $\ep,$ $x_0,$ $x_1,$ $x_2$ be as in the hypotheses of
Lemma~\ref{L:ManyElt}.
\Wolog\  $\| a \| \leq 1$ for all $a \in F.$
Choose $\ep_0 > 0$ with $\ep_0 \leq \frac{1}{7} \ep$
and so small that whenever $p_1, p_2 \in A$ are \pj s,
and $y \in A$ satisfies $\| y^* y - p_1 \| < 2 \ep_0$
and  $\| y y^* - p_2 \| < 2 \ep_0,$
then there is a partial isometry $v \in A$ such that
\[
\| v - y \| < \tfrac{1}{7} \ep,
\,\,\,\,\,\,
v^* v = p_1, \andeqn v v^* = p_2.
\]
Choose $\dt > 0$ with $\dt \leq \ep_0$
and also so small that whenever $e \in A$ satisfies
$\| e^2 - e \| < \dt$ and $\| e^* - e \| < \dt,$
then $p = \ch_{[1/2, \, \infty)} (e^* e)$ is defined and
satisfies $\| e - p \| < \ep_0.$

Apply the conditions of the lemma with $\dt$ in place of $\ep,$
obtaining $e_1,$ $e_2,$ and~$y.$
Set $p_1 = \ch_{[1/2, \, \infty)} (e_1^* e_1)$ and
$p_2 = \ch_{[1/2, \, \infty)} (e_2^* e_2),$
giving $\| p_1 - e_1 \| < \ep_0$ and $\| p_2 - e_2 \| < \ep_0.$
Then $\| y^* y - p_1 \| < \dt + \ep_0 \leq 2 \ep_0,$
and similarly $\| y y^* - \af ( p_2) \| < 2 \ep_0.$
Therefore there is a partial isometry $v \in A$
such that
\[
\| v - y \| < \tfrac{1}{7} \ep,
\,\,\,\,\,\,
v^* v = p_1, \andeqn v v^* = \af (p_2).
\]
Since also
\[
\| p_1 - e_1 \| < \tfrac{1}{7} \ep,
\,\,\,\,\,\,
\| p_2 - e_2 \| < \tfrac{1}{7} \ep,
\,\,\,\,\,\,
\| e_1 x e_1 \| > 1 - \tfrac{1}{7} \ep,
\andeqn
\| e_2 x e_2 \| > 1 - \tfrac{1}{7} \ep,
\]
we easily obtain the conditions of Lemma~\ref{L:ManyElt}.
\end{proof}

\begin{prp}\label{P:CompTAI}
Let $A$ be an \idssuca,
and let $\af, \bt \in \Aut (A)$ be tracially approximately inner.
Then $\af \circ \bt$ is tracially approximately inner.
\end{prp}

\begin{proof}
If $A$ does not have \PSP, then $\af$ and $\bt$
are approximately inner by Lemma~\ref{L:TAIAndSP},
so $\af \circ \bt$ is approximately inner.
Thus, we assume that $A$ has \PSP.

We verify the hypotheses of Lemma~\ref{L:ATAI}.
Let $F \S A$ be finite, let $\ep > 0,$
and let $x_0, x_1, x_2 \in A$ be nonzero positive
elements with $\| x_1 \| = \| x_2 \| = 1.$
\Wolog\  $\| a \| \leq 1$ for all $a \in F.$

Choose a nonzero \pj\  $r \in {\overline{x_0 A x_0}},$
and use Lemma~\ref{OrthInSP} to choose orthogonal nonzero \pj s
$r_1^{(0)}, r_2 \leq r.$
Then use Lemma~\ref{L:CompSP} to choose a nonzero \pj\  $r_1 \leq r_1^{(0)}$
such that $r_1 \precsim \bt (r_1^{(0)}).$

Choose $\dt > 0$
with $\dt \leq \min \left( \frac{1}{6} \ep, \frac{1}{40} \right),$
and also so small that whenever $B$ is a \ca\  and $p, q \in B$
are \pj s such that $\| [ p, q ] \| < \dt,$
then there exists a \pj\  $g \in B$ which commutes with $p$
and satisfies $\| g - q \| < \frac{1}{4}.$
Apply the condition of Lemma~\ref{L:ManyElt} to $\bt,$
with $F$ as given,
with $\dt$ in place of $\ep,$
with $r_1$ in place of $x_0,$
and with $x_1$ and $x_2$ as given.
Let $q_1, q_2, w \in A$ be the resulting \pj s and partial isometry.

Set
\[
E = F \cup \bt (F) \cup \{ \bt (q_2 a q_2) \colon a \in F \}
         \cup \{ w, w^*, q_1, q_2, \bt (q_2) \}.
\]
Set $\ld_1 = \| q_1 x_1 q_1 \|$ and $\ld_2 = \| q_2 x_2 q_2 \|.$
Apply the condition of Lemma~\ref{L:ManyElt} to $\af,$
with $E$ in place of $F,$
with $\dt$ in place of $\ep,$
with $r_2$ in place of $x_0,$
with $\ld_1^{-1} q_1 x_1 q_1$ in place of $x_1,$
and with $\ld_2^{-1} \bt ( q_2 x_2 q_2)$ in place of $x_2.$
Let $p_1, p_2, v \in A$ be the resulting \pj s and partial isometry.

Set
\[
e_1 = p_1 q_1, \,\,\,\,\,\,
e_2 = \bt^{-1} (p_2) q_2,
\andeqn y = v w.
\]
Condition~(\ref{L:ATAI:0}) of Lemma~\ref{L:ATAI} is obvious.
We have
$\| e_1^2 - e_1 \| \leq \| [p_1, q_1] \| < \dt \leq \ep,$
and the same estimate applies to $\| e_1^* - e_1 \|.$
Also,
\[
\| e_2^2 - e_2 \|
 \leq \| [\bt^{-1} (p_2), \, q_2] \|
 = \| [p_2, \bt (q_2)] \|
 < \dt
 \leq \ep,
\]
and again the same estimate applies to $\| e_2^* - e_2 \|.$
Next, since $v^* v = p_1$ and $w^* w = q_1,$ we get
\[
\| y^* y - e_1 \|
  = \| w^* p_1 w - p_1 q_1 \|
  \leq \| [w^*, p_1] \|
 < \dt
 \leq \ep.
\]
Furthermore, by the choice of $p_1$, $p_2,$ and $v,$
and by the definition of $e_2,$ we have
\[
\| y y^* - (\af \circ \bt) (e_2 e_2^* ) \|
 = \| v p_1 \bt (q_2) p_1 v^* - \af (p_2 \bt (q_2) p_2) \|
 < \dt.
\]
Since
$\| e_2 e_2^* - e_2 \|
   \leq \| e_2^* - e_2 \| + \| e_2^2 - e_2 \| < 2 \dt,$
we get $\| y y^* - (\af \circ \bt) (e_2) \| < 3 \dt \leq \ep.$
This verifies~(\ref{L:ATAI:1}) of Lemma~\ref{L:ATAI}.

For~(\ref{L:ATAI:2}), for $a \in F$ we estimate
\[
\| [ p_1 q_1, \, a ] \|
 \leq \| p_1 \| \cdot \| [q_1, a] \| + \| [p_1, a] \| \cdot \| q_1 \|
 < 2 \dt \leq \ep
\]
and
\[
\| [ \bt^{-1} (p_2) q_2, \, a ] \|
   \leq \| \bt^{-1} (p_2) \| \cdot \| [q_2, a] \|
            + \| [p_2, \,\bt ( a)] \| \cdot \| q_2 \|
   < 2 \dt
   \leq \ep.
\]

We verify~(\ref{L:ATAI:3}).
We have
\begin{align*}
& \| y e_1 a e_1 y^* - (\af \circ \bt) (e_2 a e_2) \|
                   \\
& \mbox{} \hspace*{3em}
    = \| v w p_1 q_1 a p_1 q_1 w^* v^*
         - (\af \circ \bt) (\bt^{-1} (p_2) q_2 a \bt^{-1} (p_2) q_2) \|
                   \\
& \mbox{} \hspace*{3em}
    \leq 2 \| [w, p_1] \| + \| [p_1, q_1] \|
    + \| [\bt^{-1} (p_2), \, q_2] \|
                    \\
& \mbox{} \hspace*{6em}
    + \| v p_1 w q_1 a q_1 w^* p_1 v^*
              - \af (p_2 \bt (q_2 a q_2) p_2) \|
                   \\
& \mbox{} \hspace*{3em}
    \leq 2 \| [w, p_1] \| + \| [p_1, q_1] \|
    + \| [p_2, \, \bt (q_2)] \|
    + \| w q_1 a q_1 w^* - \bt (q_2 a q_2) \|
                    \\
& \mbox{} \hspace*{6em}
    + \| v p_1 \bt (q_2 a q_2) p_1 v^*
               - \af (p_2 \bt (q_2 a q_2) p_2) \|
                   \\
& \mbox{} \hspace*{3em}
    < 6 \dt
   \leq \ep,
\end{align*}
as desired.

Next, we prove (\ref{L:ATAI:4})~of Lemma~\ref{L:ATAI}.
We have
\[
\| (e_1^* e_1)^2 - e_1^* e_1 \|
 \leq 2 \| e_1^* - e_1 \| + 2 \| e_1^2 - e_1 \|
 < 4 \dt
 \leq \tfrac{1}{10},
\]
so that $\frac{1}{2} \not\in \spec (e_1^* e_1),$
the \pj\  $f_1 = \ch_{[1/2, \, \I)} (e_1^* e_1)$ is defined,
and it satisfies $\| e_1^* e_1 - f_1 \| < \tfrac{1}{5}.$
So
\[
\| f_1 - e_1 \|
 \leq \| f_1 - e_1^* e_1 \| + \| e_1^* - e_1 \| + \| e_1^2 - e_1 \|
 < \tfrac{1}{5} + \dt + \dt
 \leq \tfrac{1}{4}.
\]
Moreover, by the choice of $\dt$ and because $\| [p_1, q_1] \| < \dt,$
there is a \pj\  $g_1 \in A$ which commutes with $p_1$
and such that $\| g_1 - q_1 \| < \tfrac{1}{4}.$
It follows that
\[
\| (1 - f) - [ (1 - p_1) + p_1 (1 - g_1) ] \|
 = \| p_1 g_1 - f_1 \|
 \leq \| g_1 - q_1 \| + \| e_1 - f_1 \|
 < \tfrac{1}{2}.
\]
Therefore
\begin{align*}
1 - f_1
 & \sim (1 - p_1) + p_1 (1 - g_1)
            \\
 &
   \precsim (1 - p_1) \oplus (1 - g_1)
   \sim (1 - p_1) \oplus (1 - q_1)
   \precsim r_2 + r_1
   \leq r
   \in {\overline{x_0 A x_0}}.
\end{align*}
Applying the same argument with $\bt^{-1} (p_2)$ in place of $p_1$
and with $q_2$ in place of $q_1,$
it follows that
the \pj\  $f_2 = \ch_{[1/2, \, \I)} (e_2^* e_2)$ is defined,
it satisfies $\| e_2^* e_2 - f_2 \| < \tfrac{1}{5},$
and
\[
1 - f_2
 \precsim r_2 \oplus \bt^{-1} (r_1)
 \precsim r_2 + r_1^{(0)}
 \leq r
   \in {\overline{x_0 A x_0}}.
\]

It remains only to prove~(\ref{L:ATAI:5}).
Since $\frac{1}{2} < \ld_1 \leq 1,$ we have
\[
1 \leq \ld_1^{-1} < 1 + 2 (1 - \ld_1) < 1 + 2 \dt.
\]
Therefore
\[
\| \ld_1^{-1} p_1 q_1 x_1 q_1 p_1 - p_1 q_1 x_1 q_1 p_1 \|
 < 2 \dt,
\]
so that, using $\| [p_1, q_1] \| < \dt$ at the first step,
\[
\| e_1 x_1 e_1 \|
 > \| p_1 q_1 x_1 q_1 p_1 \| - \dt
 > \big\| \ld_1^{-1} p_1 q_1 x_1 q_1 p_1 \big\| - 3 \dt
 > 1 - 4 \dt
 \geq 1 - \ep.
\]
Similar reasoning gives
\[
\| e_2 x_2 e_2 \|
 > \big\| \bt^{-1} (p_2) q_2 x_2 q_2 \bt^{-1} (p_2) \big\| - \dt
 > \big\| \ld_2^{-1} p_2 \bt ( q_2 x_2 q_2 ) p_2 \big\| - 3 \dt
 > 1 - \ep.
\]
This completes the proof.
\end{proof}

\begin{thm}\label{T:TAIIsGp}
Let $A$ be an \idssuca.
Then the set of tracially approximately inner automorphisms of $A$
is a group.
\end{thm}

\begin{proof}
Obviously $\id_A$ is tracially approximately inner.
In view of Proposition~\ref{P:CompTAI},
we need only prove that if $\af$
is tracially approximately inner,
then so is $\af^{-1}.$
We use the condition given in Lemma~\ref{L:ManyElt}.

Let $F \S A$ be finite, let $\ep > 0,$
and let $x_0, x_1, x_2 \in A$ be nonzero positive
elements with $\| x_1 \| = \| x_2 \| = 1.$
Apply the condition of Lemma~\ref{L:ManyElt}
with $F,$ $\ep,$ and $x_0$ as given,
with $x_2$ in place of $x_1,$ and with $x_1$ in place of $x_2.$
Let $q_1, q_2, w \in A$ be the resulting \pj s and partial isometry.
Then set $p_1 = q_2,$ $p_2 = q_1,$ and $v = \af^{-1} (w^*).$
These satisfy the conditions in Lemma~\ref{L:ManyElt}.
\end{proof}

\section{Properties of tracially approximately inner
                         automorphisms}\label{Sec:PrTAI}

In this section, we prove that
a tracially approximately inner automorphism is necessarily
trivial on the tracial state space and on $K_0$ mod infinitesimals.
On an \idssuca\  with tracial rank zero,
an automorphism which is
trivial on $K_0$ mod infinitesimals
is necessarily tracially approximately inner;
this is the analog of the fact that an automorphism of
an AF~algebra which is trivial on $K_0$
is necessarily approximately inner.
We also prove that if a tracially approximately inner automorphism
of an \idssuca\  with tracial rank zero
has finite order,
then it is strongly tracially approximately inner
in the sense of Definition~\ref{STAInnDfn}.
Thus, the results of Section~\ref{Sec:STAI}
apply to such automorphisms.
(This does not help with our main application of those
results~\cite{PhtRp2},
because we use them as part of the proof that the algebra
involved has tracial rank zero.)

\begin{prp}\label{P:TrAndTAI}
Let $A$ be an \idssuca,
and let $\af \in \Aut (A)$ be tracially approximately inner.
Then $\ta \circ \af = \ta$
for every tracial state $\ta$ on $A.$
\end{prp}

\begin{proof}
If $\af$ is approximately inner, the result is immediate.
We may therefore assume that $A$ has \PSP.
Let $\ta$ be a tracial state on $A,$
let $\ep > 0,$ and let $a \in A.$
We prove that $| \ta (\af (a)) - \ta (a) | < \ep.$
\Wolog\  $\| a \| = 1.$

Choose $n$ such that $\frac{1}{n} < \frac{1}{3} \ep.$
By page~61 of~\cite{AS},
there exists a selfadjoint $b \in A$ such that $\spec (a) = [0, 1].$
Choose \cfn s
\[
f_1, f_2, \ldots, f_n, g_1, g_2, \ldots, g_n \colon [0, 1] \to [0, 1]
\]
such that $\| g_k \| = 1$ and $f_k g_k = g_k$ for $1 \leq k \leq n,$
and such that the supports of $f_1, f_2, \ldots, f_n$ are disjoint.
Then $\sum_{k = 1}^{n} \ta (f_k (b)) \leq 1,$
so there exists $k$ such that $\ta (f_k (b)) \leq \frac{1}{n}.$
Choose a nonzero \pj\  $q_0 \in {\overline{g_k (b) A g_k (b)}},$
and use Lemma~\ref{L:CompSP} to find a nonzero \pj\  $q \leq q_0$
such that $q \precsim \af^{-1} (q_0).$
Since $f_k (b) q_0 = q_0,$
we have $\ta (q_0) \leq \ta (f_k (b)) \leq \frac{1}{n},$
so that $\ta (q) \leq \frac{1}{n}$ and $\ta (\af (q)) \leq \frac{1}{n}.$

Apply Definition~\ref{D:TAI}
with $F = \{ a \},$
with $\frac{1}{3} \ep$ in place of $\ep,$
and with $x = q,$
obtaining $p_1,$ $p_2,$ and $v.$
Since $\| a \| \leq 1,$ we have
\[
| \ta (a) - \ta (p_1 a p_1) | \leq \ta (1 - p_1) \leq \tfrac{1}{n}
\andeqn
| \ta (\af (a)) - \ta (\af (p_2 a p_2)) |
 \leq \ta (\af (1 - p_2)) \leq \tfrac{1}{n}.
\]
Also,
\[
\ta (v p_1 a p_1 v^*) = \ta (v^* v p_1 a p_1) = \ta (p_1 a p_1),
\]
so
\begin{align*}
| \ta (\af (a)) - \ta (a) |
 & \leq | \ta (\af (a)) - \ta (\af (p_2 a p_2)) |
                   \\
& \mbox{} \hspace*{3em}
       + \| \af (p_2 a p_2) - v p_1 a p_1 v^* \|
       + | \ta (p_1 a p_1) - \ta (a) |
          \\
 & < \tfrac{1}{n} + \tfrac{1}{n} + \tfrac{1}{3} \ep
   < \ep.
\end{align*}
This completes the proof.
\end{proof}

Recall that an element $\et$ of a partially ordered group $G$
with order unit $u \in G_+ \SM \{ 0 \}$ is
{\emph{infinitesimal}} if $- m u \leq n \et \leq m u$ for all
$m, \, n \in \N$ with $m > 0.$
See Definition~1.10 of~\cite{GPS}, where this definition is given for
simple dimension groups.
Clearly we need only consider $m = 1.$
By Proposition~4.7 of~\cite{Gd0}, an equivalent condition is that all
states on $(G, u)$ vanish on $\et.$

\begin{prp}\label{TAIAndInf}
Let $A$ be an \idssuca,
and let $\af \in \Aut (A)$ be tracially approximately inner.
Then $\af_* (\et) - \et$ is infinitesimal for every $\et \in K_0 (A).$
\end{prp}

\begin{proof}
If $A$ does not have Property~(SP),
then $\af$ is approximately inner by Lemma~\ref{L:TAIAndSP},
so $\af_* (\et) - \et = 0$ for every $\et \in K_0 (A).$
So assume that $A$ has Property~(SP).

We prove that for every $\et \in K_0 (A)$
we have $- [1_A] \leq \af_* (\et) - \et \leq [1_A].$
This implies the result, because replacing $\et$ by $n \et$
gives $- [1_A] \leq n [ \af_* (\et) - \et ] \leq [1_A].$

Accordingly, let $\et \in K_0 (A),$
and choose $n \in \N$ and \pj s $p, \, r \in M_n (A)$
such that $\et = [q] - [r].$
Let $\ch \colon \R \setminus \left\{ \frac{1}{2} \right\} \to \R$
be the characteristic function of $\left( \frac{1}{2}, \I \right).$
Choose $\ep > 0$ so small that $n^2 \ep < \frac{1}{6},$ and
also so small that whenever $C$ is a \ca\  and
$f, \, p \in C$ are \pj s such that $\| f p - p f \| < n^2 \ep,$
then $\frac{1}{2}$ is not in the spectrum of either
$f p f$ or $(1 - f) p (1 - f),$ and moreover
the \pj s $p_0 = \ch (f p f)$ and $p_1 = \ch ((1 - f) p (1 - f))$
satisfy
$\| p_0 + p_1 - p \| < \frac{1}{6}.$
Use Lemma~\ref{OrthInSP} to choose $2 n$ nonzero \mops\   %
$g_1, g_2, \ldots, g_{2 n} \in A,$
and use Lemma~\ref{L:CompSP} to choose a nonzero \pj\  $h \leq g_1$
such that $\af (h) \leq g_1.$

Apply Definition~\ref{D:TAI} with
$F = \{ q_{j, k}, \, r_{j, k} \colon 1 \leq j, \, k \leq n \},$
the set of all matrix entries of $q$ and $r,$
with $\ep$ as just chosen,
and with $x = h.$
Let $p_1, p_2, v \in e A e$
be the resulting \pj s and partial isometry.

We have
\[
\| (1 \otimes p_1) q - q (1 \otimes p_1) \|
  \leq \sum_{j, \, k = 1}^n \| p_1 q_{j, k} - q_{j, k} p_1 \| < n^2 \ep.
\]
By the choice of $\ep$ the \pj s
\[
q_1 = \ch ((1 \otimes p_1) q (1 \otimes p_1) ) \in M_n (p_1 A p_1)
\]
and
\[
e_1 = \ch ((1 - 1 \otimes p_1) q (1 - 1 \otimes p_1) )
  \in M_n ((1 - p_1) A (1 - p_1))
\]
are defined and satisfy $\| q_1 + e_1 - q \| < \frac{1}{6}.$
Similarly,
\[
q_2 = \ch ((1 \otimes p_2) q (1 \otimes p_2) ) \in M_n (p_2 A p_2)
\]
and
\[
e_2 = \ch ((1 - 1 \otimes p_2) q (1 - 1 \otimes p_2) )
  \in M_n ((1 - p_2) A (1 - p_2))
\]
are defined and satisfy $\| q_2 + e_2 - q \| < \frac{1}{6}.$
Note in particular that
\[
[q] = [q_1] + [e_1] = [q_2] + [e_2]
\]
in $K_0 (A).$

We now claim that $[q_2] = [\af (q_1)]$ in $K_0 (A).$
First,
\[
\| (1 \otimes p_1) q (1 \otimes p_1) - q_1 \|
   = \| (1 \otimes p_1) [q - (q_1 + e_1)] (1 \otimes p_1) \|
   < {\textstyle{ \frac{1}{6} }}.
\]
Similarly,
\[
\| (1 \otimes p_2) q (1 \otimes p_2) - q_2 \|
   < {\textstyle{ \frac{1}{6} }}.
\]
Finally,
\begin{align*}
& \| (1 \otimes v) (1 \otimes p_1) q (1 \otimes p_1) (1 \otimes v)^*
       - (\id \otimes \af) ((1 \otimes p_2) q (1 \otimes p_2)) \|   \\
& \hspace*{10em}
  \leq \sum_{j, \, k = 1}^n
      \| v p_1 q_{j, k} p_1 v^* - \af (p_2 q_{j, k} p_2) \|
  < n^2 \ep.
\end{align*}
Putting these estimates together gives
\[
\| (1 \otimes v) q_2 (1 \otimes v)^* - \af (q_1) \|
     < {\textstyle{ \frac{1}{6} }} + n^2 \ep
         + {\textstyle{ \frac{1}{6} }}
     < {\textstyle{ \frac{1}{2} }}.
\]
The claim follows.

Repeating the argument of the last two paragraphs with $r$ in place
of $q,$ we find \pj s
\[
r_1 \in M_n (p_1 A p_1), \,\,\,\,\,\,
f_1 \in M_n ((1 - p_1) A (1 - p_1))
\]
and
\[
r_2 \in M_n (p_2 A p_2), \,\,\,\,\,\,
f_2 \in M_n ((1 - p_2) A (1 - p_2))
\]
such that
\[
[r] = [r_1] + [f_1] = [r_2] + [f_2]
  \andeqn [r_2] = [\af (r_1)]
\]
in $K_0 (A).$

We have
\[
\af_* (\et) - \et
  = \af_* ([q]) - [q] - \af_* ([r]) + [r]
  = \af_* ([e_2]) - [e_1] - \af_* ([f_2]) + [f_1].
\]
Since
\[
e_1 \in M_n ((1 - p_1) A (1 - p_1))
\andeqn
1 - p_1 \precsim h \leq g_1,
\]
and since
\[
f_2 \in M_n ((1 - p_2) A (1 - p_2))
\andeqn
1 - \af (p_2) \precsim \af (h) \precsim g_1,
\]
in $K_0 (A)$ we have
\[
[e_1] + \af_* ([f_2])
  \leq n [1 - p_1] + n [1 - \af (p_2)]
  \leq n [h] + n [\af (h)]
  \leq 2 n [g_1]
  \leq [1_A].
\]
Similarly,
$\af_* ([e_2]) + [f_1] \leq [1_A].$
Therefore
\[
- [1_A] \leq - [e_1] - \af_* ([f_2])
        \leq \af_* (\et) - \et
        \leq \af_* ([e_2]) + [f_1]
        \leq [1_A].
\]
This completes the proof.
\end{proof}

We now prove that if $A$ has tracial rank zero,
then there is a converse to Proposition~\ref{TAIAndInf}.
We give a preliminary lemma.

\begin{lem}\label{L:InfInTRZ}
Let $A$ be a \suca\  with Property~(SP),
with stable rank one,
and such that the order on \pj s over $A$ is determined by traces
(Definition~\ref{OrdDetD}).
Let $p, q \in A$ be \pj s such that $[p] - [q]$
is infinitesimal in $K_0 (A).$
Then for every nonzero positive element $x \in A$ there
exist \pj s $p_0 \leq p$ and $q_0 \leq q$ such that
$p_0 \sim q_0$ and such that
$p - p_0$ and $q - q_0$ are \mvnt\  to \pj s in ${\overline{x A x}}.$
\end{lem}

\begin{proof}
For every tracial state $\ta$ on $A,$
we have
$- \ta_* ([1_A]) \leq n \ta_* ([p] - [q]) \leq \ta_* ([1_A])$
for all $n \in \N.$
Therefore $\ta (p) = \ta (q).$
Since $A$ is simple, if $p \in \{ 0, 1 \}$ then $q = p,$
and there is nothing to prove.
Accordingly, we may assume $p \not\in \{ 0, 1 \}.$

Lemma~\ref{OrthInSP} gives orthogonal nonzero \pj s
$r, s \in {\overline{x A x}}.$
Use Lemma~\ref{L:CompSP} to choose nonzero \pj s
$e \leq p$ and $e_0$ such that $e \sim e_0 \leq r,$
and $f \leq 1 - p$ and $f_0$ such that $f \sim f_0 \leq s.$
Set $p_0 = p - e.$
For every tracial state $\ta$ on $A,$
we have
$\ta (q) - \ta (p_0) = \ta (e) > 0.$
Since the order on \pj s over $A$ is determined by traces,
there exists a \pj\  $q_0 \leq q$ such that $q_0 \sim p_0.$
Similarly (also using stable rank one),
there exists a \pj\  $q_1 \geq q$ such that
$p + f \sim q_1.$

Clearly
\[
p - p_0 = e \precsim r \leq r + s \in {\overline{x A x}}.
\]
Also,
\[
q - q_0 \leq q_1 - q_0 \sim e + f \precsim r + s \in {\overline{x A x}}.
\]
This completes the proof.
\end{proof}

\begin{thm}\label{TAIOnTAF}
Let $A$ be an \idssuca\  with tracial rank zero.
Let $\af \in \Aut (A)$ be an automorphism such that
$\af_* (\et) - \et$ is infinitesimal for every $\et \in K_0 (A).$
Then $\af$ is tracially approximately inner.
\end{thm}

\begin{proof}
Let $F \S A$ be finite, let $\ep > 0,$
and let $x \in A$ be a positive element with $\| x \| = 1.$
Set $\ep_0 = \frac{1}{6} \ep.$
Choose a nonzero \pj\  $r \in {\overline{x A x}}.$
By Lemma~\ref{L:CompSP},
there exists a nonzero \pj\  $r \leq r_0$ such that $r - r_0 \neq 0.$
Use Proposition~\ref{TAFCond} to choose
a \pj\  $p \in A$ and a finite
dimensional unital subalgebra $E \S p A p$ such that:
\begin{enumerate}
\item\label{TAIOnTAF:1} %
$\| p a - a p \| < \ep_0$ for all $a \in F.$
\item\label{TAIOnTAF:2} %
For every $a \in F$ there exists $b \in E$ such that
$\| p a p - b \| < \ep_0.$
\item\label{TAIOnTAF:3} %
$1 - p \precsim r_0.$
\end{enumerate}
Write $E = \bigoplus_{l = 1}^m E_l$ with $E_l \cong M_{n (l)}$
for $1 \leq l \leq m.$
Let $( e_{j, k}^{(l)} )_{1 \leq j, k \leq n (l)}$
be a system of matrix units for $E_l.$

Choose (Lemma~\ref{OrthInSP})
orthogonal nonzero equivalent \pj s $r_{l, j} \leq r - r_0$
for $1 \leq l \leq m$ and $1 \leq j \leq n (l).$
By Lemma~\ref{L:CompSP},
there exist equivalent nonzero \pj s $s_{l, j} \leq r_{l, j}$
such that $\af^{-1} (s_{l, j}) \precsim r_{l, j}.$
For $1 \leq l \leq m$ use
Lemma~\ref{L:InfInTRZ} and Theorem~\ref{TAFProp}
to find \pj s $f_l \leq e_{1, 1}^{(l)}$
and $g_l \leq \af (e_{1, 1}^{(l)}),$
and $v_l \in A,$
such that
\[
v_l^* v_l = f_l, \,\,\,\,\,\,
v_l v_l^* = g_l, \,\,\,\,\,\,
e_{1, 1}^{(l)} - f_l \precsim s_{l, 1}, \andeqn
\af (e_{1, 1}^{(l)}) - g_l \precsim s_{l, 1}.
\]

Define
\[
f = \sum_{l = 1}^m \sum_{j = 1}^{n (l)}
        e_{j, 1}^{(l)} f_l e_{1, j}^{(l)},
\,\,\,\,
g = \sum_{l = 1}^m \sum_{j = 1}^{n (l)}
        \af (e_{j, 1}^{(l)}) g_l \af (e_{1, j}^{(l)}),
\,\,\,\, {\mbox{and}} \,\,\,\,
v = \sum_{l = 1}^m \sum_{j = 1}^{n (l)}
        \af (e_{j, 1}^{(l)}) v_l e_{1, j}^{(l)}.
\]
Then it is easily checked that $f$ and $g$ are \pj s,
and that $v^* v = f$ and $v v^* = g.$
So Condition~(\ref{D:TAI:1}) of Definition~\ref{D:TAI} holds
with $f$ in place of $p_1$ and with $\af^{-1} (g)$ in place of $p_2.$
We verify conditions~(\ref{D:TAI:2}) through~(\ref{D:TAI:4}).
This will verify the hypotheses of Lemma~\ref{L:TAIForFin},
completing the proof.

For~(\ref{D:TAI:2}),
let $a \in A$ and choose $b \in E$ such that
$\| p a p - b \| < \ep_0.$
Then
\[
\| a - [b + (1 - p) a (1 - p)] \|
  \leq \| p a p - b \| + 2 \| p a - a p \|
  < 3 \ep_0.
\]
It is easily checked that $f$ commutes with every element of $E,$
and $f$ commutes with $(1 - p) a (1 - p)$ because $f \leq p,$
so
$\| f a - a f \| < 6 \ep_0 \leq \ep.$
Similarly, we get $\| \af^{-1} (g) a - a \af^{-1} (g) \|< \ep.$

For~(\ref{D:TAI:3}),
let $1 \leq l \leq m$ and $1 \leq j, k \leq n (l).$
Then
\begin{align*}
v f e_{j, k}^{(l)} f v^*
 & = \left[ \af \big( e_{j, 1}^{(l)} \big) v_l e_{1, j}^{(l)} \right]
       \left[e_{j, 1}^{(l)} f_l e_{1, j}^{(l)} \right]  e_{j, k}^{(l)}
       \left[e_{k, 1}^{(l)} f_l e_{1, k}^{(l)} \right]
       \left[e_{k, 1}^{(l)} v_l^* \af \big( e_{1, k}^{(l)} \big) \right]
                         \\
 & = \af \big( e_{j, 1}^{(l)} \big)
              v_l f_l v_l^* \af \big( e_{1, k}^{(l)} \big)
   = g \af \big( e_{j, k}^{(l)} \big) g
   = \af \left( \af^{-1} (g) e_{j, k}^{(l)} \af^{-1} (g) \right).
\end{align*}
Since this is true for all $l,$ $j,$ and $k,$
it follows that
$v f b f v^* = \af \big( \af^{-1} (g) b \af^{-1} (g) \big)$
for every $b \in E.$
Now let $a \in A$ and choose $b \in E$ such that
$\| p a p - b \| < \ep_0.$
Then, using $f \leq p$ and $g \leq \af (p),$ we get
\[
\big\| v f a f v^* - \af \big(\af^{-1} (g) a \af^{-1} (g) \big) \big\|
  \leq \| f a f - f b f \| + \| g \af (a) g - g \af (b) g \|
  < 2 \ep_0
  \leq \ep,
\]
as desired.

Finally, we verify~(\ref{D:TAI:4}).
We have
\[
1 - f
  = 1 - p + \sum_{l = 1}^m \sum_{j = 1}^{n (l)}
      \left( e_{j, j}^{(l)} - e_{j, 1}^{(l)} f_l e_{1, j}^{(l)} \right)
  \precsim r_0 + \sum_{l = 1}^m \sum_{j = 1}^{n (l)} s_{l, j}
  \leq r,
\]
and similarly 
\[
1 - \af^{-1} (g)
  = 1 - p + \sum_{l = 1}^m \sum_{j = 1}^{n (l)}
    \left(
       e_{j, j}^{(l)} - e_{j, 1}^{(l)} \af^{-1} (g_l) e_{1, j}^{(l)}
                              \right)
  \precsim r_0 + \sum_{l = 1}^m \sum_{j = 1}^{n (l)} \af^{-1} (s_{l, j})
  \precsim r.
\]
This completes the proof of~(\ref{D:TAI:4}), and of the theorem.
\end{proof}

When the automorphism has finite order,
we get strong tracial approximate innerness.
Again, we need a lemma.

\begin{lem}\label{L:InvK0}
Let $A$ be a \suca\  with real rank zero,
with stable rank one,
and such that the order on \pj s over $A$ is determined by traces
(Definition~\ref{OrdDetD}).
Let $\af \in \Aut (A)$ be an automorphism such that
$\af_* \colon K_0 (A) \to K_0 (A)$ has finite order and such that
$\af_* (\et) - \et$ is infinitesimal for every $\et \in K_0 (A).$
Let $p \in A$ be a \pj.
Then for every nonzero positive element $x \in A$ there
exists a \pj\  $p_0 \leq p$ such that
$\af (p_0) \sim p_0$ and such that
$p - p_0$ is \mvnt\  to a \pj\  in ${\overline{x A x}}.$
\end{lem}

\begin{proof}
\Wolog\  $p \neq 0.$
By Lemma~\ref{L:CompSP},
there is a nonzero \pj\  $r \in {\overline{x A x}}$
such that $r \precsim p.$
Set $\ep = \inf_{\ta \in T (A)} \ta (r) > 0.$
Choose $N \in \N$ with $\frac{1}{N} < \frac{1}{2} \ep.$
Let $n$ be the order of $\af_*.$
Use Lemma~2.3 of~\cite{OP1}
to find \pj s  $f, q_0, q_1, \ldots, q_{n N} \in A$
such that
\[
f + \sum_{k = 0}^{n N} q_k = p
\andeqn
f \precsim q_0 \sim q_1 \sim \cdots \sim q_{n N}.
\]

For every $\ta \in T (A)$ and $\et \in K_0 (A),$
we have
$- \ta_* ([1_A]) \leq n \ta_* (\af_* (\et) - \et) \leq \ta_* ([1_A])$
for all $n \in \N.$
Therefore $\ta_* \circ \af_* = \ta_*.$
Set
\[
\et = \sum_{k = 1}^{n N} [ \af^k (q_k)].
\]
Then $\af_* (\et) = \et$
and $\ta_* (\et) = n N \ta (q_1) < \ta (p).$
(We get strict inequality because $f + q_0 \neq 0.$)
Since the order on \pj s is determined by traces,
there exists a \pj\  $p_0 \leq p$
such that $[p_0] = \et.$
Moreover, since $A$ has stable rank one, $\af (p_0) \sim p_0.$

For every $\ta \in T (A),$ we have
\[
(n N + 1) \ta (q_1) \leq \ta (p) \leq 1
\andeqn
\ta (f) + \ta (q_0) \leq 2 \ta (q_1).
\]
So
\[
\ta (p - p_0) = \ta (f) + \ta (q_0)
  \leq \frac{2}{n N + 1}
  < \ep
  \leq \ta (r).
\]
It follows that $p - p_0 \precsim r.$
\end{proof}

\begin{thm}\label{FOTAIOnTAF}
Let $A$ be an \idssuca\  with tracial rank zero.
Let $\af \in \Aut (A)$ be an automorphism of finite order such that
$\af_* (\et) - \et$ is infinitesimal for every $\et \in K_0 (A).$
Then $\af$ is strongly tracially approximately inner.
\end{thm}

\begin{proof}
The proof is similar to that of Theorem~\ref{TAIOnTAF}.
Let $F \S A$ be finite, let $\ep > 0,$
and let $x \in A$ be a positive element with $\| x \| = 1.$
Let $\ep_0,$ $r,$ $r_0,$ $p,$ and
$E = \bigoplus_{l = 1}^m E_l$ with $E_l \cong M_{n (l)}$
and matrix units $( e_{j, k}^{(l)} )_{1 \leq j, k \leq n (l)}$
be as there,
and further let $r_{l, j}$ and $s_{l, j}$ be as there.
For $1 \leq l \leq m$ use Lemma~\ref{L:InvK0} and Theorem~\ref{TAFProp}
to find \pj s $f_l \leq e_{1, 1}^{(l)}$
and partial isometries $v_l \in A,$
such that
\[
v_l^* v_l = f_l, \,\,\,\,\,\,
v_l v_l^* = \af (f_l), \andeqn
e_{1, 1}^{(l)} - f_l \precsim s_{l, 1}.
\]

Define
\[
f = \sum_{l = 1}^m \sum_{j = 1}^{n (l)}
        e_{j, 1}^{(l)} f_l e_{1, j}^{(l)}
\andeqn
v = \sum_{l = 1}^m \sum_{j = 1}^{n (l)}
        \af (e_{j, 1}^{(l)}) v_l e_{1, j}^{(l)}.
\]
Then it is easily checked that $f$ is a \pj,
and that $v^* v = f$ and $v v^* = \af (f).$
So Condition~(\ref{STAInnDfn:1}) of Definition~\ref{STAInnDfn} holds
with $f$ in place of $e.$
The rest of the proof is essentially the same as the rest of
the proof of Theorem~\ref{TAIOnTAF},
using Lemma~\ref{TAIForFinite} in place of Lemma~\ref{L:TAIForFin}.
\end{proof}

In Theorem~\ref{FOTAIOnTAF},
the hypothesis that $\af$ have finite order may
be replaced by the assumption
that $\af_* \colon K_0 (A) \to K_0 (A)$
have finite order,
or by the assumption that the fixed points of $\af_*$
on $K_0 (A)$ have dense image in the real affine functions
on $T (A).$
One might hope to weaken it to the requirement that
for every nonzero positive element $x \in A,$
there be a nonzero \pj\  $p \in {\overline{x A x}}$
such that $\af_* ([p]) = [p].$
We do not know if this is possible.

\end{document}